\documentclass[12pt,english]{article}  
\usepackage{geometry}
\usepackage{setspace}
\usepackage{lineno}
\usepackage{graphicx}
\newenvironment{keywords}{\noindent\textbf{Keywords:} }{\par}
\usepackage{amsmath, amssymb, amsthm}

\usepackage{authblk}       
\usepackage{times}  
\usepackage{lmodern}
\usepackage{dcolumn}
\usepackage{bm}
\usepackage[utf8]{inputenc}
\usepackage[T1]{fontenc}
\usepackage{mathptmx}
\usepackage{etoolbox}
\usepackage{subcaption}
\usepackage{array}
\usepackage{float}
\usepackage{hyperref}
\usepackage{multicol}
\usepackage[table]{xcolor}
\usepackage{booktabs}

\title{The First Mathematical Model for Elk–Wolf Interaction in Yellowstone National Park Using the E-SINDy Algorithm}

\author{
Anurag Singh\(^{1}\),
Nitu Kumari\(^{1}\),
and Arun Kumar\(^{1}\)\\
\(^1\)School of Mathematical \& Statistical Sciences, Indian Institute of Technology Mandi, Himachal Pradesh, India-175005\\
Author(s) Email: anuragdu2203@gmail.com; nitu@iitmandi.ac.in; arunanuj94@gmail.com
}
\date{}

\begin{document}

\maketitle

 
\section*{Abstract}
The ecological dynamics between elk and wolves in Northern Yellowstone have been a focal point of long-term research, particularly following the reintroduction of wolves to the region. Although numerous studies have explored this prey-predator interaction from ecological and behavioral perspectives, there remains a lack of comprehensive analysis using mathematical modeling approaches capable of uncovering underlying dynamical patterns and system-level insights.  
In this study, we investigate the prey-predator dynamics of the elk–wolf system in northern Yellowstone National Park, USA, using a data-driven modeling approach. We used yearly population data for elk and wolves from $1995$ to $2022$ ($28$ years) to construct a mathematical model using a sparse regression modeling framework. To the best of our knowledge, no previous work has applied this framework to capture elk–wolf interactions over this time period. Our modeling pipeline integrates Gaussian process regression for data smoothing, sparse identification of nonlinear dynamics for model discovery, and model selection techniques to identify the most suitable mathematical representation. The resulting model is analyzed for its non-linear dynamics with ecologically meaningful parameters. Stability and bifurcation analyzes are then performed to understand the system's qualitative behavior. A saddle-node bifurcation identifies parameter ranges where both species can coexist, while regions outside this range may lead to the extinction of one or both populations. Hopf and saddle-node bifurcations together delineate zones of stable co-existence, periodic oscillations, and extinction scenarios. Furthermore, co-dimension two bifurcations, including Bogdanov–Takens and cusp bifurcations, are explored by varying two parameters simultaneously. Ecologically, these bifurcations reflect the complex interplay between wolf pressure and elk defence mechanisms, such as grouping or herd behavior. They suggested that small changes in ecological parameters can lead to sudden shifts in population outcomes ranging from stable co-existence to extinction or oscillatory cycles.

\begin{keywords}
Gaussian Process Regression, SINDy, Elk-Wolf, Ecological Modeling, Nonlinear Dynamics.
\end{keywords}

\section{Introduction} \label{sec:1}
Scientists are increasingly using mathematical models to understand how animal populations change over time. These models are important tools in ecology because they help researchers understand how different species interact and how changes in the environment, like habitat loss or climate change affect ecosystems. As human activities continue to damage nature and the climate changes rapidly, these models are becoming essential for planning how to protect the environment. They not only help to explain how ecosystems work, but also guide governments and organizations in making better decisions about conservation and sustainability \cite{kendall2015some}.

\begin{figure}[H]
	\centering
	\includegraphics[width=\textwidth]{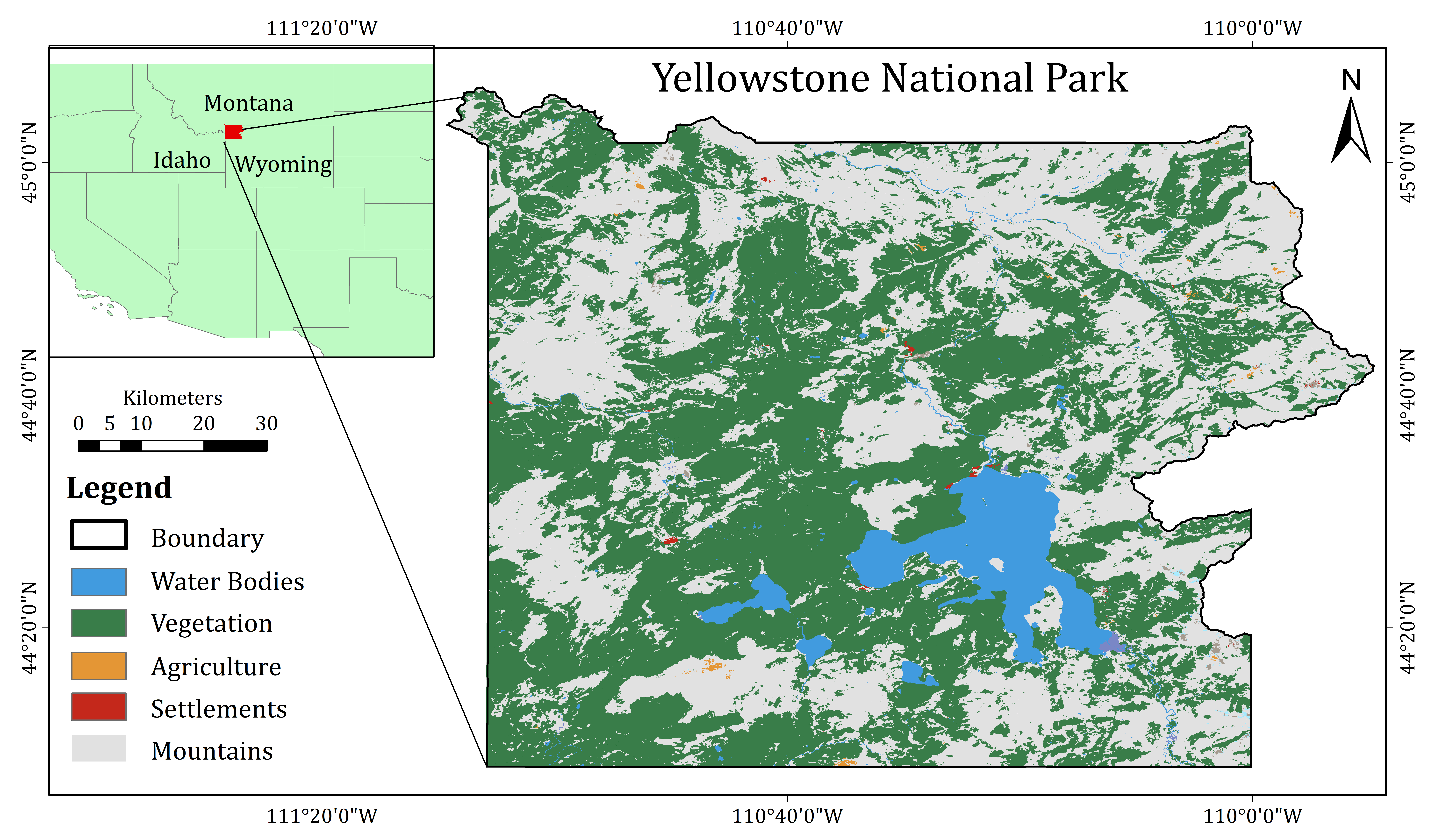}
	\caption{Study area map} 
	\label{Fig:1}
\end{figure}

Yellowstone National Park, established in $1872$, is one of the most studied ecosystems in the world. Yellowstone National Park, located primarily in the state of Wyoming in the USA and extending to Montana and Idaho (Fig. \ref{Fig:1}), is the world’s first National Park and one of the most ecologically important protected areas in North America. 

The northern Yellowstone range, a lower elevation area with relatively milder winters, supports some of the highest wildlife concentrations in the park and has been the focal point of intensive research on prey-predator interactions, especially involving gray wolves (Canis lupus) and elk (Cervus canadensis). The northern range is the hub of wildlife activity in Yellowstone. It covers only $10\%$ of the park, but it harbors the largest wintering elk herd in the park and is classified as one of the highest density areas for carnivores in all of North America. In Yellowstone, the wolves (Canis lupus) were extirpated by mid-$1920s$ and did not return for nearly seven decades until reintroduction in the winters of $1995$ and $1996$ \cite{smith2009reintroduction, macnulty2020population}. This sequence of removal and return presents an uncommon long-term natural experiment to examine cascading effects on the food web. Without wolves, elk have a far-reaching influence on vegetation, soil quality, and habitat for other species. Before the reintroduction of wolves, woody vegetation, such as aspen and willows, was unable to mature in the canopy forest in the winter range of the northern park, the only exception being areas where fencing offered protection \cite{ripple2012trophic, nps_wolf_yellowstone}.

Several studies on the prey-predator relationship between elk and wolf are being carried out by several scientists, including ecologists and biologists \cite{macnulty2020population, nps_wolf_yellowstone, boyce1992wolf, macnulty2012nonlinear, passoni2024investigating}. These studies include field studies, theoretical ecology studies, and statistical studies on how reintroduction of wolves affects the population of elk and other animals in northern areas. It also includes the effect of predation on the ecosystem and vegetation of the park. During the last three decades, this system has served as a natural laboratory to examine complex ecological interactions, including the influences of climate variability, human hunting, and multipredator effects on elk populations.

In northern Yellowstone, wolf packs typically average around $11-12$ individuals. In summer, elk continue to be the main food source, comprising approximately $85\%$ of kills, with deer accounting for around $14\%$, and bison making up a small fraction (less than $1\%$). Wolves in northern Yellowstone are also known to engage in aggressive encounters with other wolves and large carnivores like coyotes and cougars, usually in disputes over territory or access to carcasses. Disease occasionally causes mortality among both wolf pups and older individuals. Notably, outbreaks of canine distemper were recorded in $2005, 2008, \text{ and } 2009$ \cite{nps_wolf_yellowstone}. The largest elk herd in Yellowstone winters along the park's northern boundary, particularly in the Lamar and Yellowstone river valleys, where milder temperatures and lower snowfall support large numbers of elk. Today, most of this northern herd migrates outside the park into the Custer Gallatin National Forest and adjacent private lands. Historically, concerns focused on overgrazing due to high elk numbers, but more recently, attention has shifted to the herd’s declining size \cite{houston1982northern}. The decline in elk numbers has been linked to the recovery of predators (wolves, cougars, bears), human hunting, and climate factors like drought, which affect reproduction and survival \cite{NPS_Elk_Yellowstone}.

Mathematical modeling is a tool in ecology that helps researchers understand the ecosystem related to the dynamics of prey and predators \cite{may2007theoretical}. Ecological modeling began in the early $20^{th}$ century as scientists sought to describe population interactions and dynamics in ecosystems. One of the first milestones in the $1920s$ was the Lotka-Volterra model \cite{lotka1925elements}, where prey-predator interactions were characterized with differential equations. During the mid $20^{th}$ century, models developed to incorporate competition, resource dynamics, and the spread of diseases. With the development of computing in the second half of the century, ecological modeling became more sophisticated, in which spatial structures, climatic influences, and human impacts were considered \cite{murray2007mathematical}. Few mathematical modeling studies are performed on the Yellowstone data, which is based on parameter estimation \cite{menden2019computationally, hatter2019assessment}. These studies have been performed on a presumed models but in real life there are many factors which influences ecology of the park which is not always possible through traditional mathematical modeling \cite{sugihara2012detecting}. This highlights the need for more flexible, data-adaptive approaches that can reveal hidden patterns using ecological time series.  

In today’s world data-driven techniques has emerged as a powerful tool for modeling due to the emergence of high-volume data. Sparse identification of nonlinear dynamics (SINDy) \cite{brunton2016discovering} is one of the data-driven method used to find the governing equations of a variety of dynamical systems. While SINDy has transformed fields like fluid dynamics and engineering, its ecological applications remain limited  \cite{fasel2022ensemble, dos2023new, kumari2025data, singh2025modeling}. Our study bridges this gap by developing a comprehensive framework that combines: Gaussian process regression, sparse identification, information criterion, stability, and bifurcation theory.

In this article, we have adopted a modeling framework which comprises Gaussian process regression (GPR) for data smoothing, sparse regression for mathematical model discovery, and model selection techniques to select the best mathematical model. After getting the best mathematical model with fits smoothed data we performed nonlinear study of the model to find some interesting dynamics in elk-wolf system.
The main contributions of this study are as follows:
\begin{itemize}
	\item We develop a mathematical model to characterize the long-term dynamics between elk and wolf populations, grounded in empirical data from northern Yellowstone National Park.
	\item By leveraging a multi-decadal prey-predator dataset, we provide a quantitative framework to investigate the ecological mechanisms driving species interactions over nearly thirty years.
	\item We analyze the system’s stability and bifurcation structure with respect to key ecological parameters, enabling a deeper understanding of threshold behaviors and nonlinear responses in elk–wolf interactions.
\end{itemize}
By integrating cutting edge computational methods with short-term ecological data, this work advances both theoretical ecology and wildlife management practice. Our approach offers a template for studying other complex ecosystems where traditional modeling approaches have proven inadequate.

\section{Methodology}  \label{sec:2}
\begin{figure}[H]
	\centering
	\includegraphics[width=\textwidth]{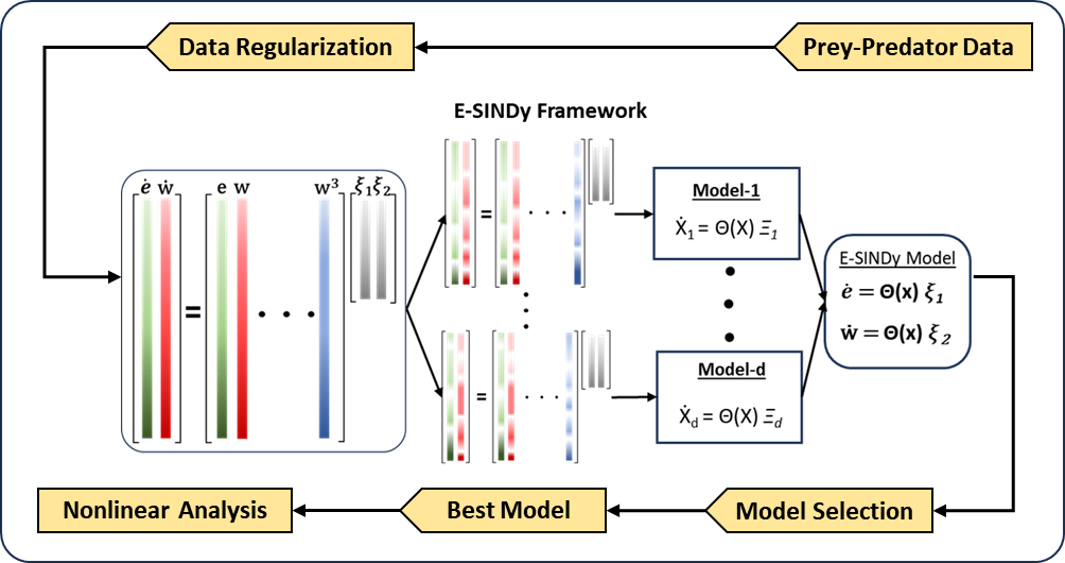}
	\caption{Schematic diagram of the modeling and analysis framework} 
	\label{Fig:3}
\end{figure}

In this work, we develop a framework to study and interpret the prey-predator interactions between elk and wolves. As shown in Fig. \ref{Fig:3}, our approach combines modern techniques such as data regularization, sparse model discovery, and model selection with classical nonlinear analysis to uncover the critical parameters driving the elk–wolf dynamics. It starts with time series observations that record population fluctuations of both species. Since ecological data are often noisy and uncertain, we first apply data regularization methods to reduce variability and obtain smoother, more reliable trajectories for further study.
We have used this smooth time series data to discover ordinary differential equations using E-SINDy method \cite{fasel2022ensemble}. The system identification approach is applied under multiple experimental setups, using varied parameter choices and candidate functions. Each setup has the potential to yield a different model. Further, we used information criterion to find best model among these models. Finally, stability and bifurcation study has been carried out to find interesting dynamics between elk and wolves.
All the simulations has been carried out using open source PySINDy package \cite{Kaptanoglu2022}, scikit-learn Machine learning library, MATLAB, and Jupyter notebook.

\subsection{Time Series Data}  \label{sec:2.1}
\begin{figure}[H]
	\centering
	\includegraphics[width=\textwidth]{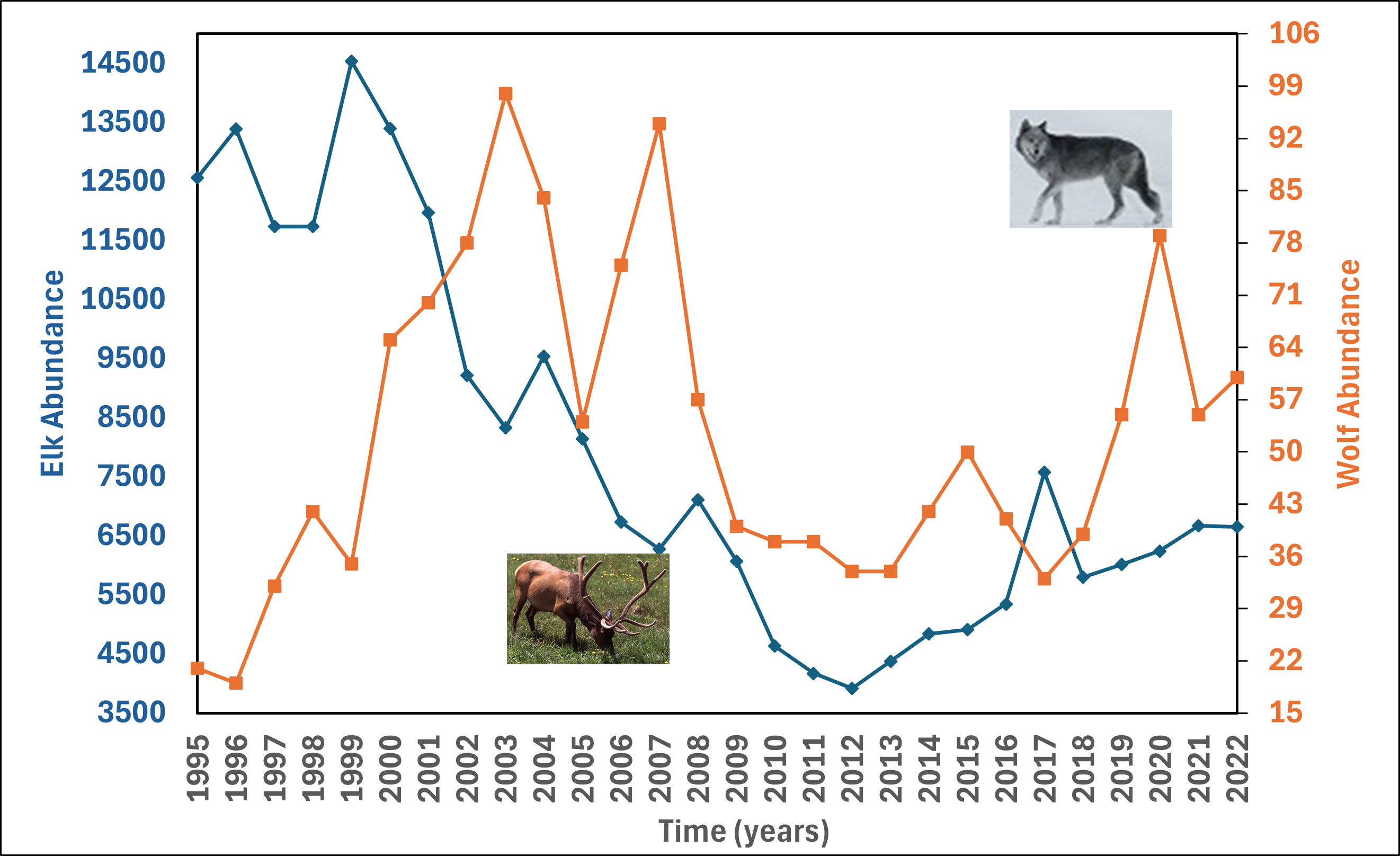}
	\caption{ Time series plot of the Northern Yellowstone National Park yearly elk and wolf population dataset. Elk and wolf photos are from Yellowstone Park \cite{NPS_YellowstonePhotography}.} 
	\label{Fig:2}
\end{figure}

The dataset utilized in this study is obtained through direct communication with the Yellowstone Centre for Resources, U.S. National Park Service, providing annual population estimates for northern Yellowstone of elk (Cervus canadensis) \cite{NPS_Elk_Yellowstone} and wolves (Canis lupus) \cite{cassidy2025yellowstone} from $1995$ to $2022$. Fig. \ref{Fig:2} presents these data, illustrating the prey-predator system composed of the wolves preying on the elk. In the elk population dataset, population estimates for the years $1995$, $2005$, and $2013$ were missing. To address this, different imputation methods are applied based on data availability. The missing value for the year $1995$ is estimated by taking the average of the elk population recorded in $1996$ and $1997$. For $2005$ and $2013$, the missing values are interpolated by computing the mean of the population data from the preceding and succeeding years. Analysis of this time-series data reveals significant fluctuations in both populations, indicative of complex prey-predator dynamics. Over the study period, the elk population ranged from a minimum of $3,915$ to a maximum of $14,539$, with a mean of $7,926.71$ and a standard deviation of $3,141.41$. The wolf population exhibited variation between $19$ and $98$, with an average of $52.21$ and a standard deviation of $20.88$.

The time series in this study is denoted as $
D_i^{(j)} = D^{(j)}(t_i), \text{ for } i = 1, \dots, 28 \text{ and } j \in \{e, w\},
$
where \( e \) represents the elk population (prey) and \( w \) denotes the wolf population (predator).

\subsection{Data Regularization} \label{sec:2.2}  

Ecological data are often characterized by large variability and uncertainty levels, thereby making it extremely challenging to model the underlying dynamics of such a complex system in mathematical terms. The critical characteristics of such complex systems have to be captured by strong techniques that smooth data in addition to addressing their uncertainty. Gaussian Process Regression (GPR) \cite{williams1995gaussian, libotte2022framework} has shown itself to be a very effective technique in this respect. A nonparametric flexible technique is provided by GPR to build smooth approximations from noisy data along with confidence intervals quantifying uncertainty. GPR is used as a Bayesian nonlinear regression technique, where a Gaussian process prior is assumed, and the resulting predictive distribution is a multivariate Gaussian.

For a given new time \( t^* \) and the corresponding variable \( D^* \), the conditional PDF is given by:  
\[
p(D^*|t^*, t, D) = \mathcal{N} (\eta^*, \Sigma^*),
\]  
where the mean function is given as:  
\[
\eta^* = K(t^*, t) (K(t, t) + \sigma^2 I)^{-1} D,
\]  
and the covariance matrix is defined by:  
\[
\Sigma^* = K(t^*, t^*) - K(t^*, t) (K(t, t) + \sigma^2 I)^{-1} K(t, t^*),
\]  
for $D = D^j, j \in \{e, w\}$. The resulting regularized data are respectively defined as $\eta^{*j}, j \in \{e, w\}$.
Here, \( K(\cdot, \cdot) \) represents the kernel function. The kernel function quantifies the similarity between pairs of input points and serves as a crucial component for capturing and representing the underlying structure of the observed data.

Since the elk and wolf populations exhibit distinct temporal patterns, different kernel functions are employed for their modeling. For the elk data, we utilize the Matérn kernel \cite{williams2006gaussian}, which provides an effective compromise between smoothness and flexibility, making it suitable for capturing gradual ecological variations in population dynamics. In contrast, the wolf data is modeled using the Radial Basis Function (RBF) kernel \cite{williams2006gaussian}, as its capacity to handle rapid changes makes it well-suited for representing the pronounced fluctuations in wolf population trends. The GPR models are implemented using the Scikit-learn library, where kernel hyperparameters are optimized based on initial estimates. To mitigate the risk of convergence to local minima, 500 restarts are performed during optimization. Specifically, for the elk data, the Matérn kernel is initialized with a smoothness parameter of $\nu = 2.5$, while for the wolf data, the RBF kernel is initialized with a length-scale parameter constrained within the range of 5 to 100.

\subsection{Mathematical models} \label{sec:2.3}
Mathematical models serve as simplified representations of real-world systems. They are designed to capture the most critical ecological components relevant to the study at hand. By abstracting complex biological processes into a set of equations or rules, these models help researchers analyze, predict, and understand the behavior of ecological systems under various conditions.

One effective method for developing such models is through data-driven techniques, particularly the Ensemble Sparse Identification of Nonlinear Dynamics (E-SINDy) approach. This framework enables the discovery of governing equations directly from noisy and less amount of data, minimizing assumptions about underlying mechanisms. Using E-SINDy, we can construct interpretable and parsimonious models that accurately reflect the key dynamical features of the elk-wolf system, guided by empirical observations rather than purely theoretical formulations.

\subsubsection{Data-driven models}  \label{sec:2.4}
In this study, we use data-driven models, derived through system identification applied to the regularized elk and wolf dataset. Recent advances in data-driven approaches \cite{brunton2016discovering, fasel2022ensemble} have demonstrated their effectiveness in uncovering nonlinear dynamical systems from noisy time-series data. A large body of this work builds upon the Sparse Identification of Nonlinear Dynamics (SINDy) framework, with applications reported in fields such as epidemiology \cite{horrocks2020algorithmic, jiang2021modeling} and ecology \cite{dos2023new}. In our analysis, we adopt the E-SINDy method for model discovery, which enhances the robustness of the standard SINDy approach by leveraging ensemble techniques, making it particularly well-suited for scenarios involving limited and noisy ecological data.

The E-SINDy framework operates by utilizing the data together with a predefined library of candidate functions, with the objective of identifying the minimal set of functions required to accurately capture the underlying system dynamics. For a detailed exposition of the E-SINDy methodology, the reader is referred to \cite{fasel2022ensemble}. In this study, we specifically employ the library b(r)agging strategy as part of the ensemble procedure. Therefore we are dropping two candidate terms from library on each ensemble. 

Since different inputs can lead to the identification of distinct models, we formulate several experimental settings by varying the regularization parameter $\alpha$, SINDy threshold parameter $\lambda$, number of models for ensemble, and number of data points going into ensembeling. In particular, We have taken $\alpha \in \{0.001,~0.01,~0.1,~1\}$, $\lambda \in \{0.01,~0.02,~0.03,~0.04,~0.05\}$, number of models $\in \{20,~50,~100\}$, and number of data points going into ensembeling $\in \{70\%,~80\%,~90\%\}$. Also, we have choosen the set of candidate functions, $\Theta_3(X)$ composed by cubic polynomial functions. The choice of $\Theta_3(X)$ is given by experiments. First we used the library, $\Theta_2(X)$ consists of second order polynomial functions. But second order polynomial library is unable to find the mathematical model which best fit the data. Then we have used $\Theta_3(X)$, with the above hyper parameters and we obtained a mathematical model which fits the data. Also we take another library which consists higher order polynomial terms which fits the data same as $\Theta_3(X)$ but we have not considered those models because as we add more polynomial terms to the library, the model obtained by E-SINDy method will be more complex and hard to interpret. In traditional ecological models, we have mostly seen terms upto third degree that is why we have finally choose our library as $\Theta_3(X)$.

\subsection{Model selection}  \label{sec:2.5}
We then used two established model selection criteria, Akaike Information Criterion (AIC) and Bayesian Information Criterion (BIC) \cite{mangan2017model}. These are extensively used in ecological studies \cite{JOHNSON2004101}. Such measures are especially useful when more than one model is capable of explaining a given dataset, as in the case in this study. Although all two have a similar goodness-of-fit component, they vary in model complexity penalty to reduce bias. While AIC based measures are generally advised to be used in the case of medicine, biology, and social sciences \cite{kp1998model}, a single criterion is never best for all situations \cite{brewer2016relative}.

\subsection{Nonlinear analysis}  \label{sec:2.7}
After getting the final model, we perform nonlinear analysis of the model. We find the equilibrium points of the model and check the stability of the model around these equilibrium points. We also perform the bifurcation analysis of the model, which revels interesting dynamics of the model.

\section{Results}
In this section, we present and discuss the results obtained through the application of the methodologies described earlier. These results serve to validate the effectiveness of the proposed approaches and provide insights into their practical implications.

\subsection{Data Regularization} \label{sec:2.8}
Before applying the GPR, we have employed Z-score normalization (also called standardization)

\begin{equation} \label{eq:1}
	z_i = \frac{x_i-\mu}{\sigma}
\end{equation}

\begin{figure}[H]
	\centering
	\includegraphics[width=\textwidth, height=0.28\textheight]{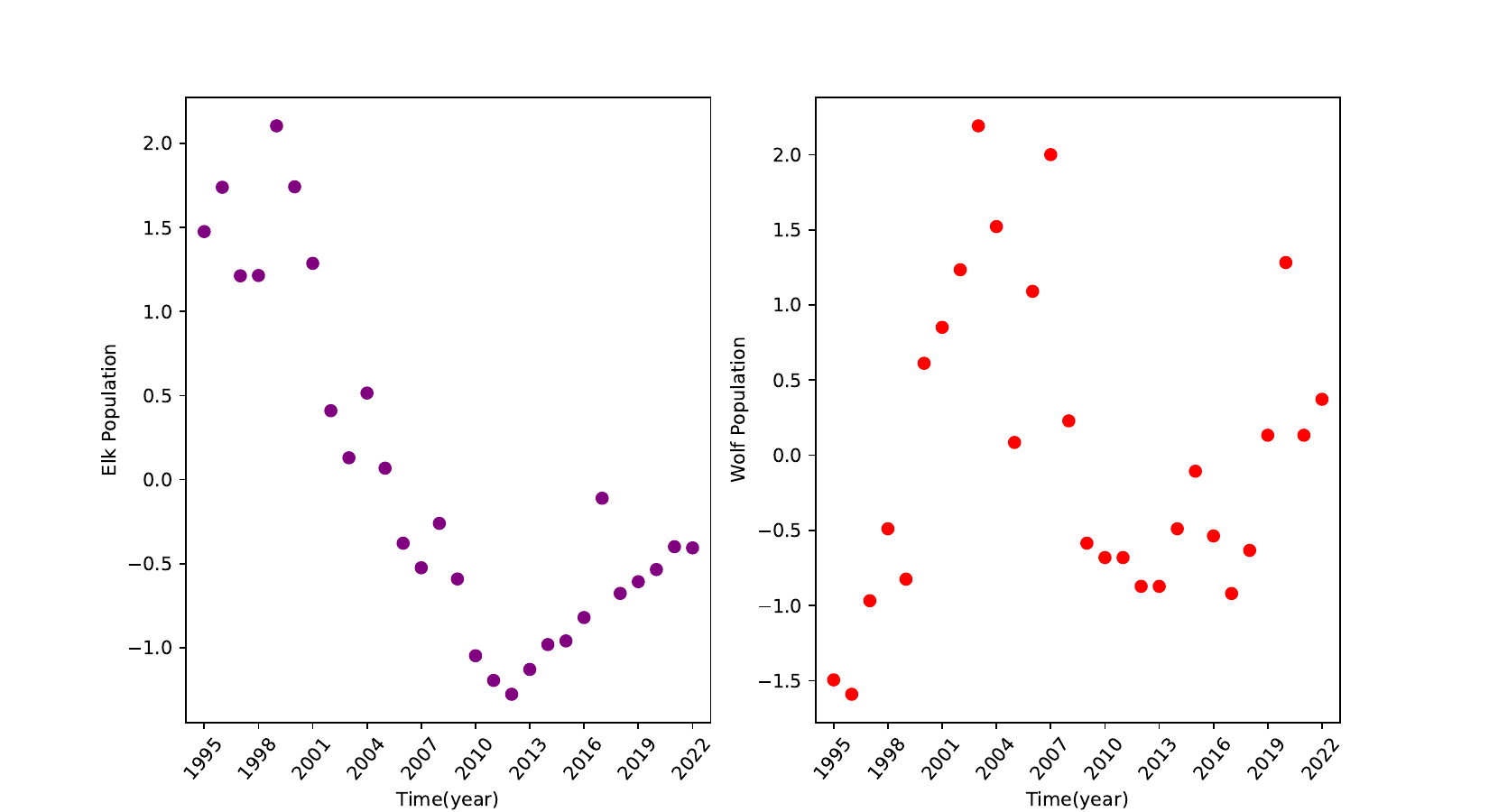}
	\caption{Elk and Wolf populations after Z-score normalization of the original dataset} 
	\label{Fig:4}
\end{figure}

(where $z_i$ = normalized value, $x_i$ = original value, $\mu$ = mean of data, and $\sigma$  = standard deviation of data) to transform the original data so that it has a mean of $0$ and a standard deviation of $1$ (Fig. \ref{Fig:4}). This process is important for several reasons: (1) brings different scales to a common scale, (2) improves performance of GPR algorithm, and (3) handles outliers better than other normalization techniques, such as min-max scaling.

\begin{figure}[H]
	\begin{subfigure}{0.5\textwidth} 
		\centering
		\includegraphics[width= \textwidth]{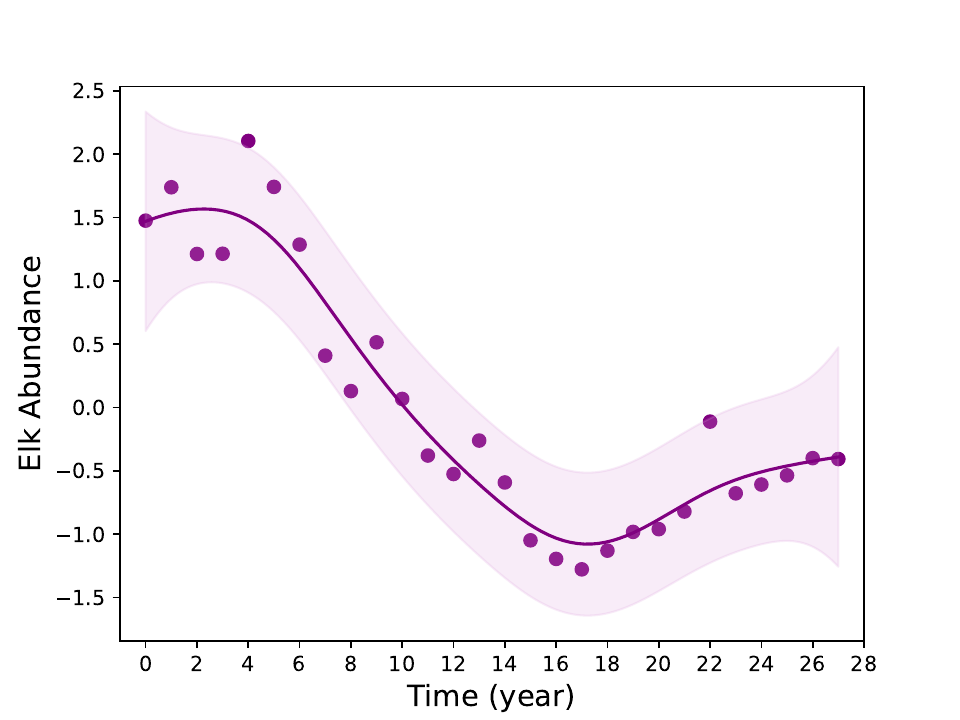} 
		\caption{\centering}\label{Fig:i}
	\end{subfigure} 
	\begin{subfigure}{0.5\textwidth}
		\centering
		\includegraphics[width= \textwidth]{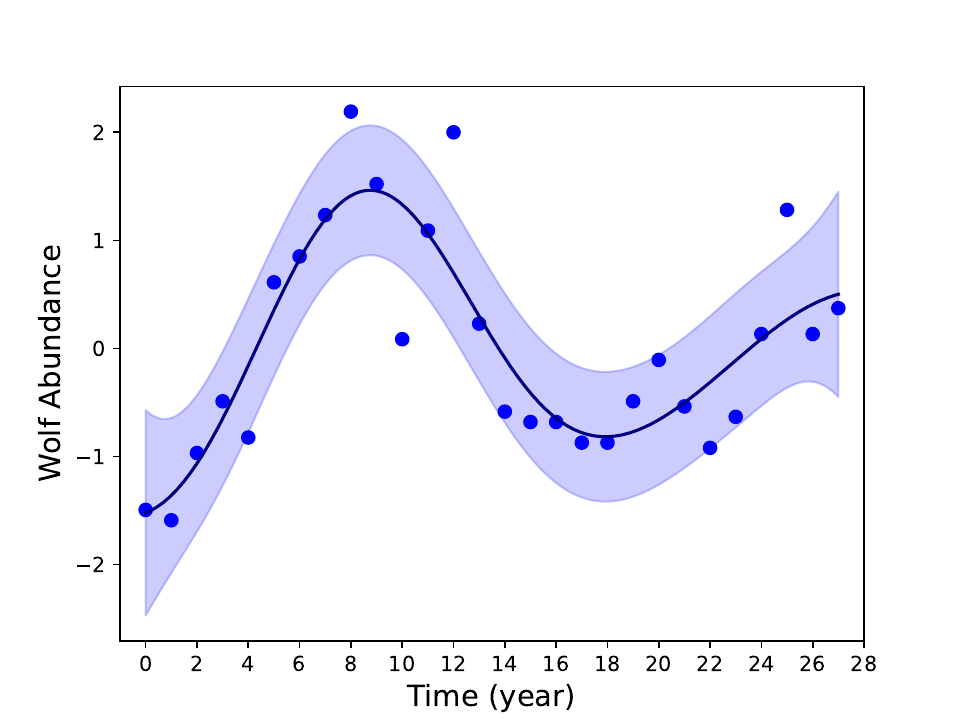} 
		\caption{\centering}\label{Fig:ii}
	\end{subfigure}
	\caption{Regularisation of (a) elk, and (b) wolf data. The solid lines are the mean values, and the shaded areas are the $95\%$ confidence interval of the GPR. The dots are the original time series data.}
	\label{Fig:5}   
\end{figure}

Fig. \ref{Fig:5} depicts regularized data for elk (Fig. \ref{Fig:i}) and wolf (Fig. \ref{Fig:ii}) populations, represented by purple and blue solid curves, respectively. The GPR model also captured almost all observed data within its uncertainty band. We employed regularized data $\eta^*$ rather than the original data set but evaluated it at $200$ equally spaced time instances. We further rescaled the dataset using eq. (\ref{eq:1}) because we want all values positive, which can be seen in Fig. \ref{Fig:6}.

\subsection{Mathematical model}
The rescaled data (Fig. \ref{Fig:6}) has different scales, so we divide all the data points by its respective standard deviation. It only changes the scale of the data not the pattern of the data. For the further analyses, we have used the dataset shown in Fig. \ref{Fig:7}.

\noindent
\begin{minipage}[t]{0.48\textwidth}
	\begin{figure}[H]
		\centering
		\includegraphics[width=\textwidth, height=0.2\textheight]{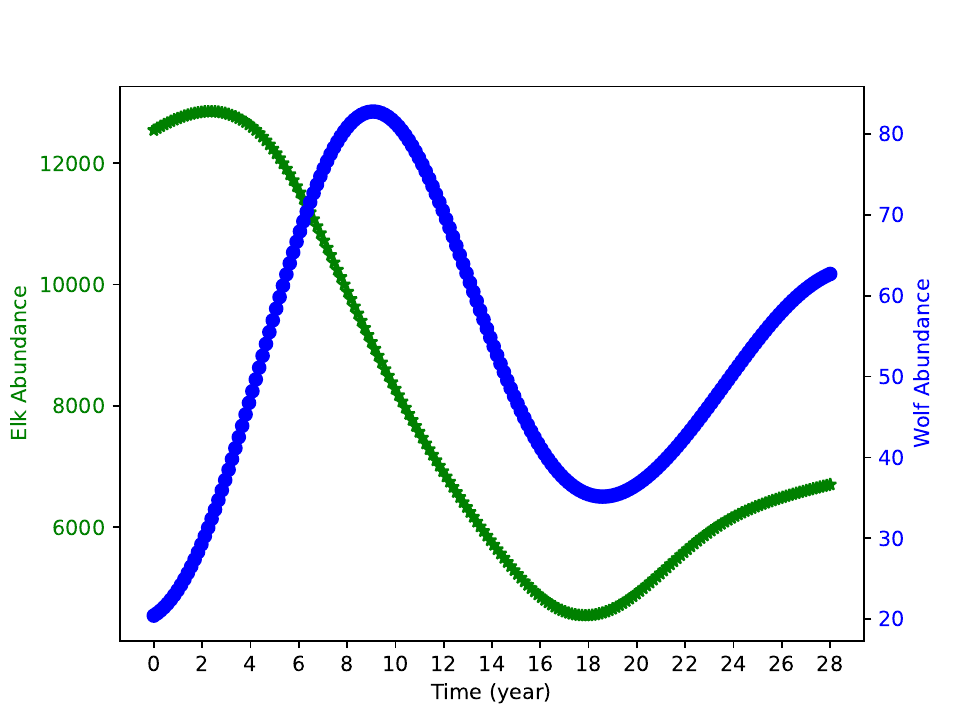}
		\caption{Rescaled data after regularization containing 200 points.} 
		\label{Fig:6}
	\end{figure}
\end{minipage}\hfill
\begin{minipage}[t]{0.48\textwidth}
	\begin{figure}[H]
		\centering
		\includegraphics[width=\textwidth, height=0.2\textheight]{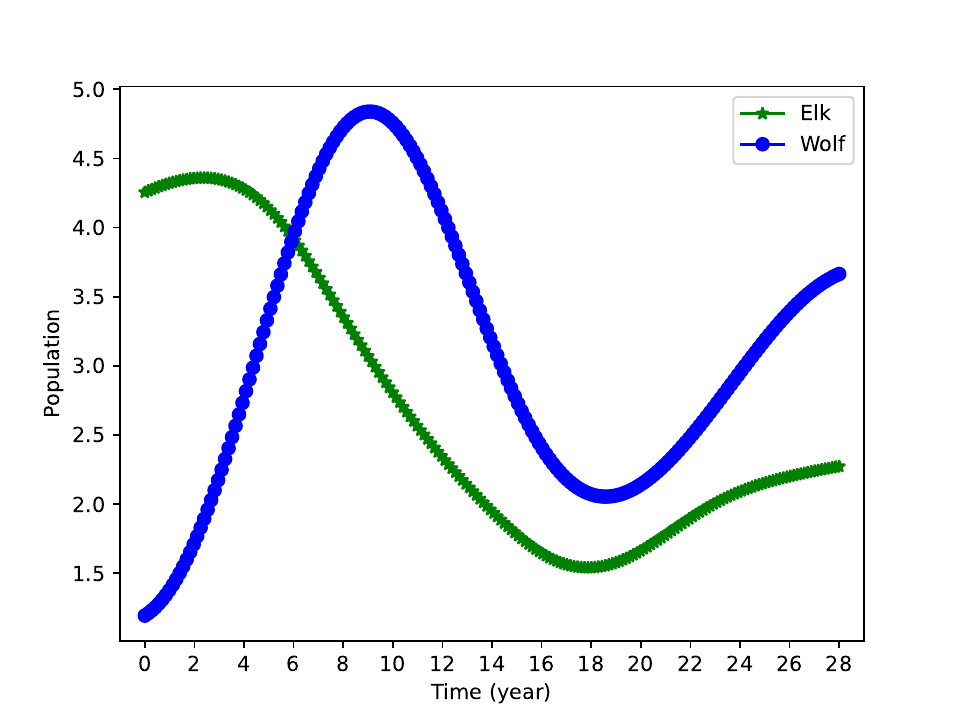}
		\caption{Standard deviation normalization of rescaled data.} 
		\label{Fig:7}
	\end{figure}
\end{minipage}

\vspace{1cm}

We assume that $X$ is the matrix which contains time series data of the species population in which first state variable is the elk population, which is denoted by $e$ and second state variable is the wolf population, which is denoted by $w$. The second-order finite difference approach is used to calculate the derivative $\dot{X}$ from the data. We then select potential functions that cover the terms included in the real model to create the experimental settings.  After running the E-SINDy algorithm we have got $117$ different mathematical models. As shown in Table \ref{table:3} (see in Appendix A), we have calculated AIC and BIC for each of the $117$ model systems. Based on the model selection criterion, we chose the best and the most
economical model using Table \ref{table:3}. The best model formulated from the data using the E-SINDy algorithm is given by

\begin{equation} \label{eq:2}
	\begin{aligned}
		\dot{e} &= 1.782 - 0.504e - 2.038w + 1.357ew - 0.175e^2w - 0.039ew^2 \\
		\dot{w} &= 5.366 - 1.586e - 5.744w + 3.478ew + 0.144w^2 - 0.405e^2w - 0.110ew^2 - 0.012w^3.
	\end{aligned}
\end{equation}

\begin{figure}[H]
	\centering
	\includegraphics[width=\textwidth, height=10cm]{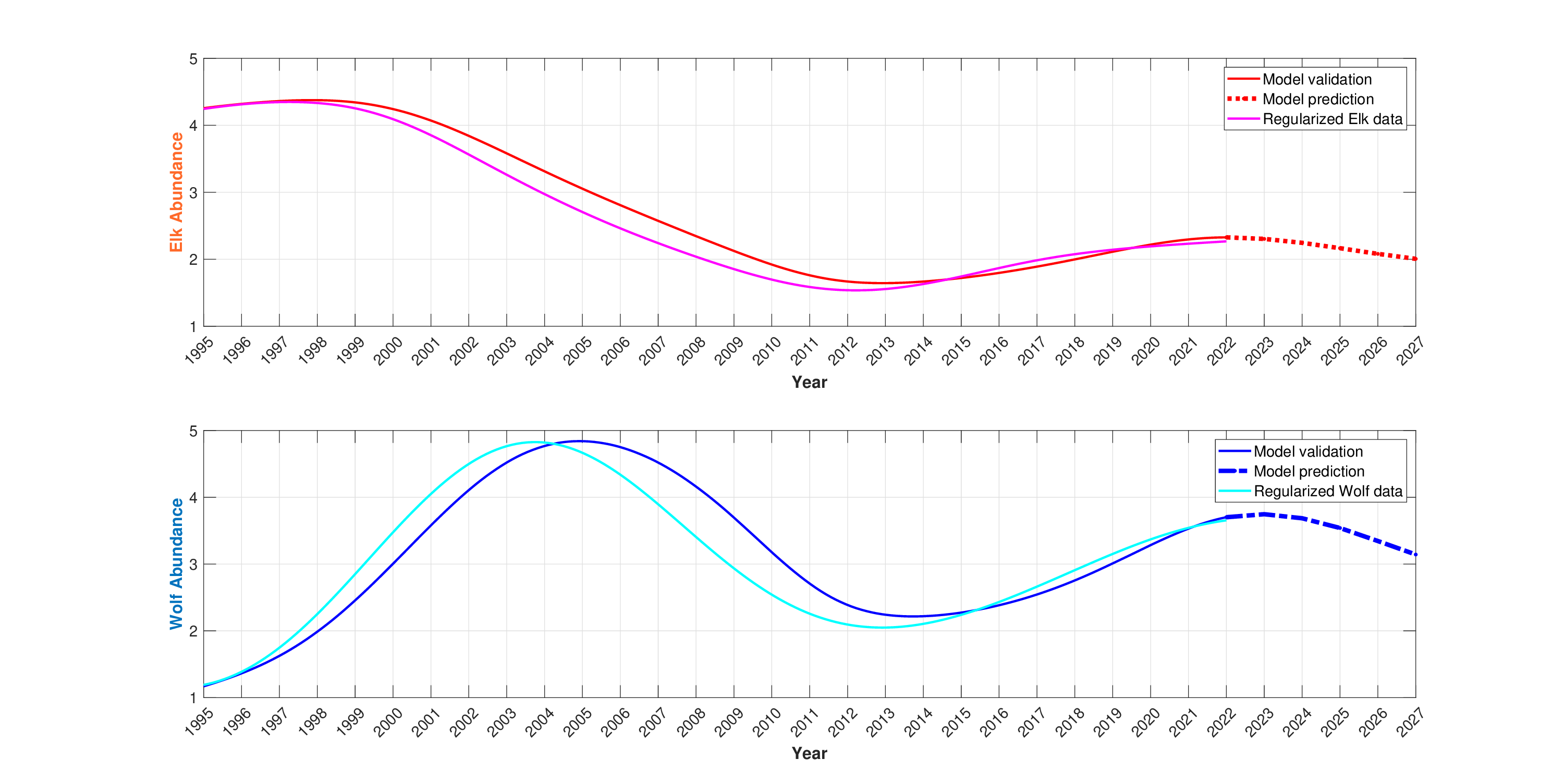}
	\caption{Validation and prediction of Elk (red line) and wolf (blue line) population using E-SINDy formulated model (\ref{eq:2}) with regularized data ($200$ points). The solid line represents model validation and dotted line shows results from model prediction.} 
	\label{Fig:8}
\end{figure}

We consider time $t$ from $0$ to $28$ with $200$ equally spaced time points. To perform the numerical simulation, the initial condition used is $[e_0,~ w_0] = [4.25493696,~ 1.17008188]$.
From Fig. \ref{Fig:8}, it can be seen that the E-SINDy model is able to capture the dynamics of the regularized data. There is a slight trade-off between model prediction and regularized data. Also, the model shows the prediction of the populations of elk and wolves for the next five years, i.e. from $2023$ to $2027$. 

In Fig. \ref{Fig:8}-\ref{Fig:bifurcations}, the populations of elk and wolf are observed in normalized space. For physical interpretation, these values can be rescaled back by multiplying with their respective standard deviations. Since the population scales of elk and wolf are different, we have used normalized space for better clarity.

\section{Nonlinear analysis}
In this section, we conduct a nonlinear analysis of the derived differential equations. We begin by defining the significance of each parameter in the model. Subsequently, we perform a detailed stability and bifurcation analysis.

\subsection{Ecological interpretation of the model} 
The earlier model, derived using the data-driven E-SINDy algorithm, effectively captures the dynamics of the elk–wolf interaction in Yellowstone National Park. However, short-term dynamics alone are insufficient to understand the full spectrum of ecological behavior, particularly the long-term patterns, nonlinear feedbacks, and parametric sensitivities inherent in prey-predator systems. To gain deeper insights into the ecological stability, potential bifurcations, and the influence of inter- and intra-species interactions, we reformulate system (\ref{eq:2}) by introducing ecologically meaningful parameters that reflect direct predation, population saturation, cooperative hunting behavior, and resource competition. These parameters, carefully interpreted in the context of elk–wolf dynamics, allow us to analyze both transient and asymptotic behaviors of the system. Therefore, we rewrite system (\ref{eq:2}) using the parameter set described in Table \ref{table:1}, aiming to facilitate a detailed exploration of the system’s nonlinear dynamics and ecological implications.
\begin{equation} \label{eq:4}
	\begin{aligned}
		\dot{e} &= a_0 - a_1 e - a_2 w + a_3 e w - a_4 e^2 w - a_5 e w^2 \\
		\dot{w} &= b_0 - b_1 e - b_2 w + b_3 e w + b_4 w^2 - b_5 e^2 w - b_6 e w^2 - b_7 w^3.
	\end{aligned}
\end{equation}


Several empirical and modeling studies of the Yellowstone elk-wolf system provide strong justification for incorporating ecologically meaningful parameters such as direct predation, population saturation, cooperative hunting behavior, and resources competition into our model.
MacNulty et. al demonstrated that wolf hunting success increases with pack size, particularly when targeting large prey like bison, thus supporting the inclusion of a cooperative hunting parameter \cite{macnulty2014influence}. Metz et. al quantified seasonal variation in wolf kill rates, highlighting how environmental factors such as snow depth influence prey vulnerability and necessitating seasonally driven saturation terms in the model \cite{metz2012seasonal}. Finally, studies following wolf reintroduction have shown shifts in elk dynamics from food limitation to predator limitation underscoring the importance of parameters for population saturation and intraspecific resource competition \cite{smith2004winter}. 

For the ecological interpretation of the parameters $a_2$ and $a_3$, we can rewrite the first equation of the model~(\ref{eq:4}) as $ \frac{de}{dt} = (a_3 e - a_2) w$,
assuming that the other parameters are zero. For a fixed number of wolves, the derivative becomes positive if $(a_3 e - a_2) w > 0 \Rightarrow e > \frac{a_2}{a_3} = e_c$, and negative if $e < \frac{a_2}{a_3}$. Thus, the fraction $e_c = \frac{a_2}{a_3}$ defines a critical threshold of the elk population.

\begin{figure}[htbp]
	\centering
	\includegraphics[width=\textwidth]{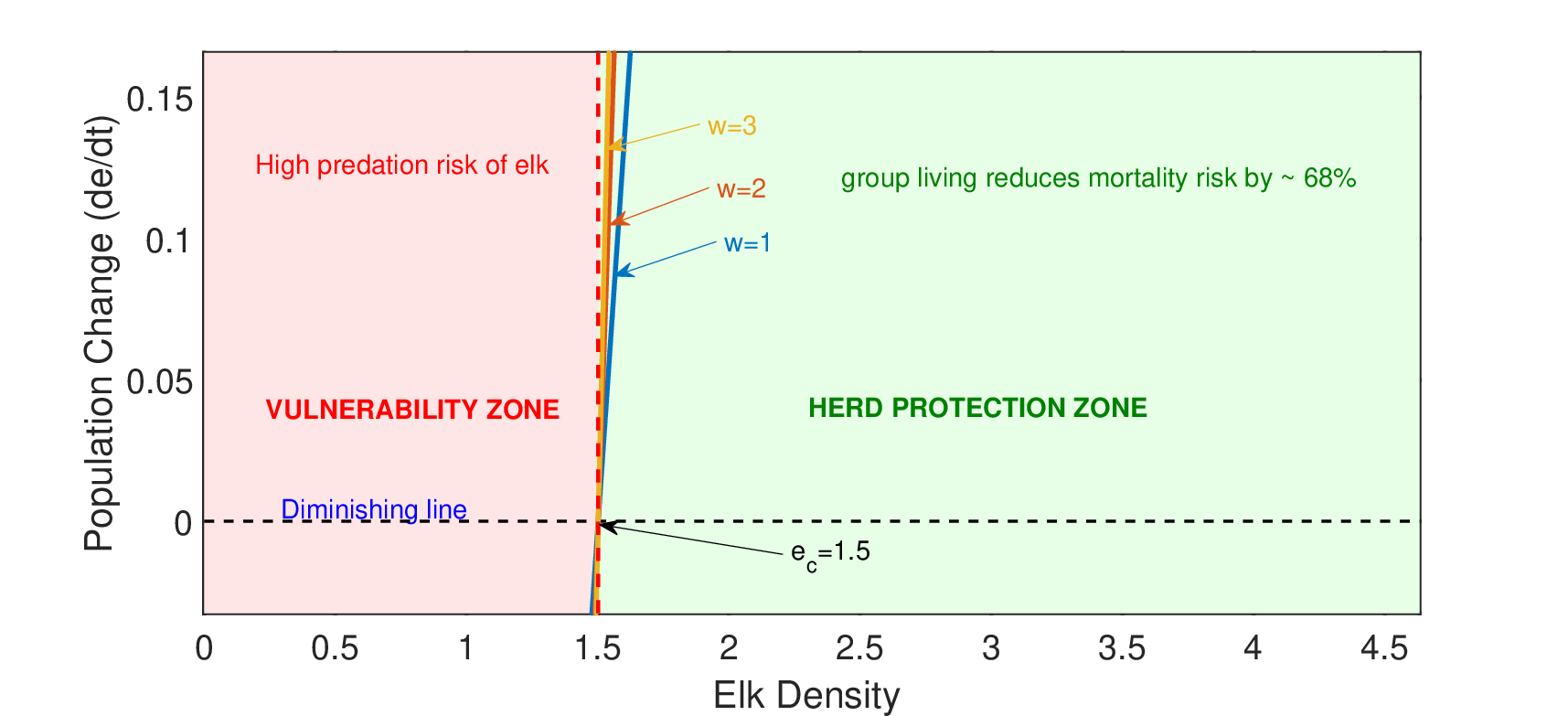}
	\caption{Elk population change $\frac{de}{dt} = (a_3 e - a_2) w$ in Yellowstone National Park. Critical threshold $e_c$ separates vulnerability (red) and herd protection (green) zones. Parameters: $a_2$ (baseline predation) and $a_3$ (herd effect of elk population) are from Table \ref{table:1}.}
	\label{fig:herd_effect}
\end{figure}

The model parameters $a_2=2.038$ and $a_3=1.357$ can be compared to these empirical patterns. The baseline predation rate $a_2=2.038$ is very close to the observed winter kill rate of $\approx 1.9$ elk per wolf-month ($\approx 22$ elk/year per wolf) \cite{smith2004winter}. In other words, Smith et al. \cite{smith2004winter} reported empirical kill rate is of the same order as $a_2$. The herd‐protection factor $a_3=1.357$ implies that group living reduces risk of elk mortality by about 0.681 (since $2.038-1.357=0.681$, roughly a $68\%$ drop). Thus, herd size reduces per-capita kill risk by $68\%$.
The rescaled critical value $e_c$ is approximately 4712, which defines the starting herd size of the elk population.
To illustrate this concept, we refer to Fig. \ref{fig:herd_effect}, where the red shaded area (below $e_c = a_2/a_3 = 1.5$) represents the vulnerability zone, and the green shaded area (above $e_c$) defines the herd protection zone. The vulnerability zone is where the elk predation is high due to their incapability to form group for protecting themselves. The herd protection zone is  where low predation of elk is seen due to herd behaviour. The diminishing protection line separates the positive and negative rate of change of elk population. For ecological validation and further insights into elk herd behavior, please refer to the research article \cite{creel2005responses}.

Collectively, these findings validate our model's design, which integrates direct predation coefficients, carrying capacity-driven saturation terms, pack size-dependent cooperative predation, and resource competition parameters to reflect the complex ecology of the Yellowstone elk-wolf system. The ecological interpretation of these parameters is based on the literature of Yellowstone National Park.

\renewcommand{\arraystretch}{0.76}
\begin{table}[H]
	\centering
	\caption{Descriptions of parameters involved in model (\ref{eq:4}).}
	\label{table:1}
	\begin{tabular}{|c|p{6.7cm}|p{3cm}|p{4.9cm}|}
		\hline
		\textbf{Parameter} & \textbf{Ecological Description} & \textbf{E-SINDy value}&\textbf{Range} \\ \hline
		$a_0$ & Constant growth rate of elk population & $1.782$& $[0.7497892, 2.3018789]$\\ \hline
		$a_1$ & Natural mortality of elk & $0.504$ & $[0.37553122, 0.5229818]$\\ \hline
		$a_2$ & Baseline predation coefficient: the fraction of elk lost per wolf, which does not necessarily imply that wolves consume the elk after killing them & $2.038$& $[1.8944218, 2.0522361]$\\ \hline
		$a_3$ & Herd protection coefficient for elk population& $1.357$ & $[1.389831, 1.3481359]$\\ \hline
		$a_4$ & Saturation coefficient of elk population & $0.175$& $[0.16729211, 0.18071395]$ \\ \hline
		$a_5$ & Wolf interference competition coefficient (higher wolf density decreases elk) & $0.039$& $[0.027339292, 0.43366126]$ \\ \hline
		$b_0$ & Growth rate, immigration, or constant reproduction of wolf population & $5.366$& $[4.4139183, 5.4476791]$ \\ \hline
		$b_1$ & Measures the Starvation in wolves due to low elk availability & $1.586$& $[1.5369023,1.8383628]$ \\ \hline
		$b_2$ & Natural mortality of wolf population & $5.744$& $[5.7079497, 6.0280868]$ \\ \hline
		$b_3$ & Predator-prey conversion coefficient (wolves benefit from eating elk) & $3.478$& $[3.4091542, 3.5012869]$ \\ \hline
		$b_4$ & Cooperative benefit coefficient; could indicate behaviors like wolf packs achieving higher hunting success rate & $0.144$& $[0.031872062, 0.16088238]$ \\ \hline
		$b_5$ & Predator decline due to elk overcrowding/ herd behavior of wolves & $0.405$& $[0.39023032,0.42255967]$ \\ \hline
		$b_6$ & Competition coefficient for limited elk (high wolf density causes higher competition for prey population) & $0.110$& $[0.0123105,0.13585081]$ \\ \hline
		$b_7$ & Intraspecific competition coefficient (higher wolf density leads to mortality due to fighting, starvation or disease) & $0.012$& $[0.0036901601, 0.070005091]$ \\ \hline
	\end{tabular}
\end{table}

\subsection{Co-existence through bifurcation structures}
%

Ecologically, each parameter plays a crucial role in determining system dynamics. While the E-SINDy derived model (\ref{eq:2}) can capture long-term dynamics, it cannot directly reveal how parameter perturbations affect these dynamics. To address this, we introduce a parameter-dependent formulation (\ref{eq:4}) that enables systematic bifurcation analysis across all system parameters. Before examining bifurcations, we must first establish the existence and stability of interior equilibrium states, where both populations coexist. This foundational analysis is essential for two reasons. Ecological relevance: only stable equilibria correspond to biologically sustainable population levels. Analytical prerequisite: stability determines how the system responds to perturbations. Accordingly, our investigation will proceed as follows: first, we identify all possible interior equilibrium states of system (\ref{eq:4}); next, we analyze their stability properties; and finally, we examine bifurcations arising from parameter variations. This approach ensures we properly characterize the system's behavior before exploring its parametric sensitivity.

For system (\ref{eq:4}), we focus exclusively on the interior equilibria where both the population coexist before conducting a deeper bifurcation and stability analysis. Let $[e,w]$ be an interior equilibrium state obtain by solving the right hand part of the system (\ref{eq:4}), we get 
$$e={\frac {-a_{{4}}w^2 \left( b_{{7}}w-b_{{4}} \right) +w \left( b_{{5}}a_{{2}}-a_{{4}}b_{{2}} \right) -b_{{5}}a_{{0}}+a_{{4}}b_{{0}}}{{w}^{2}
		\left( -b_{{5}}a_{{5}}+a_{{4}}b_{{6}} \right) +w \left( b_{{5}}a_{{3}}-a_{{4}}b_{{3}} \right) -b_{{5}}a_{{1}}+a_{{4}}b_{{1}}}}$$ and  the second component $w$ is the positive root of the equation
\begin{equation}\label{epoly}
	P_7\,{w}^{7}+P_6\,{w}^{6}+P_5\,{w}^{5}+P_4\,{w}^
	{4}+P_3\,{w}^{3}+P_2\,{w}^{2}+P_1 w+P_0=0,
\end{equation}
where $P_7={\it b_7}\, \left( {\it b_5}\,{{\it a_5}}^{2}+{{\it a_4}}^{2}{\it b_7}-{
	\it b_6}\,{\it a_5}\,{\it a_4} \right)$,\\ 
$P_6=-2\,{{\it a_4}}^{2}{\it b_4}\,{\it b_7}+{\it b_6}\,{\it a_3}\,{\it a_4}\,{
	\it b_7}+{\it b_3}\,{\it a_5}\,{\it a_4}\,{\it b_7}+{\it b_6}\,{\it a_5}\,{
	\it a_4}\,{\it b_4}-2\,{\it b_7}\,{\it b_5}\,{\it a_3}\,{\it a_5}-{\it b_4}\,
{\it b_5}\,{{\it a_5}}^{2},$\\ 
$P_5=-{\it b_6}\,{\it b_5}\,{\it a_5}\,{\it a_2}+{\it b_2}\,{\it b_5}\,{{\it a_5}}
^{2}-{\it b_1}\,{\it a_5}\,{\it a_4}\,{\it b_7}-{\it b_6}\,{\it a_5}\,{\it 
	a_4}\,{\it b_2}+2\,{\it b_4}\,{\it b_5}\,{\it a_3}\,{\it a_5}+{\it b_7}\,{
	\it b_5}\,{{\it a_3}}^{2}+{{\it b_6}}^{2}{\it a_4}\,{\it a_2}-{\it b_3}\,{
	\it a_5}\,{\it a_4}\,{\it b_4}+{{\it a_4}}^{2}{{\it b_4}}^{2}-{\it b_3}\,{
	\it a_3}\,{\it a_4}\,{\it b_7}-2\,{\it b_5}\,{\it a_2}\,{\it a_4}\,{\it b_7}-
{\it b_6}\,{\it a_1}\,{\it a_4}\,{\it b_7}-{\it b_6}\,{\it a_3}\,{\it a_4}\,{
	\it b_4}+2\,{{\it a_4}}^{2}{\it b_2}\,{\it b_7}+2\,{\it b_7}\,{\it b_5}\,{
	\it a_1}\,{\it a_5},$\\ 
$P_4=-2\,{\it b_4}\,{\it b_5}\,{\it a_1}\,{\it a_5}+2\,{\it b_5}\,{\it a_2}\,{
	\it a_4}\,{\it b_4}+{\it b_3}\,{\it b_5}\,{\it a_5}\,{\it a_2}+{\it b_6}\,{
	\it b_5}\,{\it a_3}\,{\it a_2}+{\it b_3}\,{\it a_3}\,{\it a_4}\,{\it b_4}-{
	\it b_4}\,{\it b_5}\,{{\it a_3}}^{2}-{{\it b_6}}^{2}{\it a_4}\,{\it a_0}+{
	\it b_6}\,{\it a_3}\,{\it a_4}\,{\it b_2}+{\it b_1}\,{\it a_3}\,{\it a_4}\,{
	\it b_7}+{\it b_1}\,{\it a_5}\,{\it a_4}\,{\it b_4}-2\,{\it b_7}\,{\it b_5}\,
{\it a_1}\,{\it a_3}+{\it b_3}\,{\it a_5}\,{\it a_4}\,{\it b_2}+{\it b_6}\,{
	\it a_5}\,{\it a_4}\,{\it b_0}-{\it b_0}\,{\it b_5}\,{{\it a_5}}^{2}-2\,{
	\it b_3}\,{\it a_4}\,{\it b_6}\,{\it a_2}-2\,{\it b_2}\,{\it b_5}\,{\it a_3}
\,{\it a_5}+{\it b_6}\,{\it b_5}\,{\it a_5}\,{\it a_0}-2\,{{\it a_4}}^{2}{
	\it b_0}\,{\it b_7}+{\it b_3}\,{\it a_1}\,{\it a_4}\,{\it b_7}+2\,{\it b_5}\,
{\it a_0}\,{\it a_4}\,{\it b_7}-2\,{{\it a_4}}^{2}{\it b_2}\,{\it b_4}+{\it 
	b_6}\,{\it a_1}\,{\it a_4}\,{\it b_4},$\\ 
$P_3=-{\it b_6}\,{\it b_5}\,{\it a_3}\,{\it a_0}-{\it b_3}\,{\it b_5}\,{\it a_5}\,
{\it a_0}-2\,{\it b_5}\,{\it a_0}\,{\it a_4}\,{\it b_4}+2\,{\it b_1}\,{\it 
	a_4}\,{\it b_6}\,{\it a_2}+2\,{\it b_0}\,{\it b_5}\,{\it a_3}\,{\it a_5}+2\,{
	\it b_3}\,{\it a_4}\,{\it b_6}\,{\it a_0}-2\,{\it b_5}\,{\it a_2}\,{\it a_4}
\,{\it b_2}-{\it b_1}\,{\it a_5}\,{\it a_4}\,{\it b_2}-{\it b_1}\,{\it b_5}\,
{\it a_5}\,{\it a_2}-{\it b_3}\,{\it b_5}\,{\it a_3}\,{\it a_2}+2\,{\it b_2}
\,{\it b_5}\,{\it a_1}\,{\it a_5}+2\,{\it b_4}\,{\it b_5}\,{\it a_1}\,{\it 
	a_3}-{\it b_6}\,{\it a_3}\,{\it a_4}\,{\it b_0}-{\it b_3}\,{\it a_3}\,{\it a_4
}\,{\it b_2}-{\it b_3}\,{\it a_5}\,{\it a_4}\,{\it b_0}-{\it b_6}\,{\it b_5}
\,{\it a_1}\,{\it a_2}-{\it b_3}\,{\it a_1}\,{\it a_4}\,{\it b_4}-{\it b_6}\,
{\it a_1}\,{\it a_4}\,{\it b_2}-{\it b_1}\,{\it a_1}\,{\it a_4}\,{\it b_7}-{
	\it b_1}\,{\it a_3}\,{\it a_4}\,{\it b_4}+{\it b_7}\,{\it b_5}\,{{\it a_1}}^{
	2}+{{\it b_3}}^{2}{\it a_4}\,{\it a_2}+{\it b_2}\,{\it b_5}\,{{\it a_3}}^{2}
+2\,{{\it a_4}}^{2}{\it b_0}\,{\it b_4}+{{\it a_4}}^{2}{{\it b_2}}^{2}+{{
		\it b_5}}^{2}{{\it a_2}}^{2},$\\ 
$P_2=-{\it b_0}\,{\it b_5}\,{{\it a_3}}^{2}+2\,{\it b_5}\,{\it a_2}\,{\it a_4}\,{
	\it b_0}+{\it b_3}\,{\it a_3}\,{\it a_4}\,{\it b_0}+{\it b_1}\,{\it a_5}\,{
	\it a_4}\,{\it b_0}+{\it b_1}\,{\it b_5}\,{\it a_3}\,{\it a_2}-{{\it b_3}}^{2
}{\it a_4}\,{\it a_0}+{\it b_6}\,{\it b_5}\,{\it a_1}\,{\it a_0}-2\,{\it b_1}
\,{\it a_4}\,{\it b_3}\,{\it a_2}-2\,{{\it b_5}}^{2}{\it a_0}\,{\it a_2}+{
	\it b_3}\,{\it b_5}\,{\it a_3}\,{\it a_0}+{\it b_1}\,{\it b_5}\,{\it a_5}\,{
	\it a_0}-2\,{\it b_2}\,{\it b_5}\,{\it a_1}\,{\it a_3}+{\it b_6}\,{\it a_1}\,
{\it a_4}\,{\it b_0}-2\,{\it b_1}\,{\it a_4}\,{\it b_6}\,{\it a_0}+{\it b_1}
\,{\it a_3}\,{\it a_4}\,{\it b_2}-2\,{{\it a_4}}^{2}{\it b_0}\,{\it b_2}+{
	\it b_3}\,{\it a_1}\,{\it a_4}\,{\it b_2}+{\it b_3}\,{\it b_5}\,{\it a_1}\,{
	\it a_2}-{\it b_4}\,{\it b_5}\,{{\it a_1}}^{2}-2\,{\it b_0}\,{\it b_5}\,{
	\it a_1}\,{\it a_5}+2\,{\it b_5}\,{\it a_0}\,{\it a_4}\,{\it b_2}+{\it b_1}\,
{\it a_1}\,{\it a_4}\,{\it b_4},$\\
$P_1=-{\it b_1}\,{\it b_5}\,{\it a_1}\,{\it a_2}-{\it b_1}\,{\it b_5}\,{\it a_3}\,
{\it a_0}+{{\it a_4}}^{2}{{\it b_0}}^{2}-{\it b_3}\,{\it b_5}\,{\it a_1}\,{
	\it a_0}+{{\it b_5}}^{2}{{\it a_0}}^{2}+{\it b_2}\,{\it b_5}\,{{\it a_1}}^{2
}-{\it b_1}\,{\it a_1}\,{\it a_4}\,{\it b_2}+2\,{\it b_1}\,{\it a_4}\,{\it 
	b_3}\,{\it a_0}-2\,{\it b_5}\,{\it a_0}\,{\it a_4}\,{\it b_0}-{\it b_1}\,{
	\it a_3}\,{\it a_4}\,{\it b_0}+{\it a_4}\,{{\it b_1}}^{2}{\it a_2}-{\it b_3}
\,{\it a_1}\,{\it a_4}\,{\it b_0}+2\,{\it b_0}\,{\it b_5}\,{\it a_1}\,{\it 
	a_3},$\\  
$P_0={\it b_1}\,{\it b_5}\,{\it a_1}\,{\it a_0}+{\it b_1}\,{\it a_1}\,{\it a_4}\,{
	\it b_0}-{\it a_4}\,{{\it b_1}}^{2}{\it a_0}-{\it b_0}\,{\it b_5}\,{{\it a_1}
}^{2}.$

Due to the presence of big expressions of $P_{i}$ for $i=0,1,2,\ldots 7$, an analytical finding of the positive solution of the equation (\ref{epoly}) is not possible. Therefore, we rely on the numerical set of the parameters and obtained that there exists a maximum of three positive interior equilibrium points for the system (\ref{eq:4}). 

\begin{figure}[H]
	\centering
	\includegraphics[width=\textwidth, height=7cm]{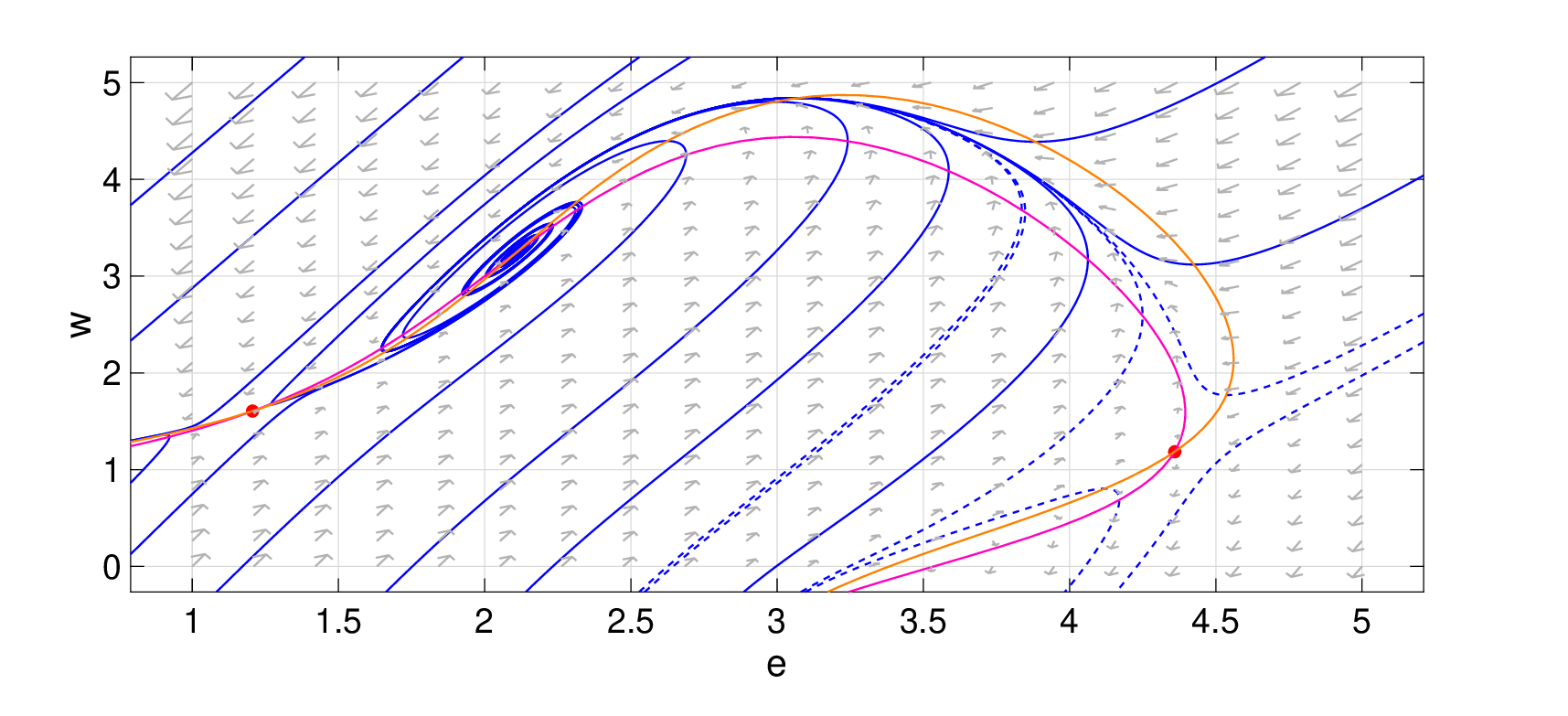}
	\caption{Dynamical behavior of all three positive equilibrium states of the system (\ref{eq:4}) with the parameters listed in the Table \ref{table:1}.} 
	\label{Fig:stability}
\end{figure}

The stability analysis of possible positive equilibrium states, using the original parameters from Table \ref{table:1}, reveals three distinct equilibria (Fig. \ref{Fig:stability}). The leftmost equilibria ([1.206091,~1.605300]) and rightmost equilibria ([4.359587,~1.184805]) exhibit saddle behavior, separated by a stable equilibrium ([2.107854,~3.243900]). Population trajectories originating in the first quadrant converge asymptotically toward this stable point. The right saddle point displays hyperbolic properties, represented by dotted divergence from the unstable manifold. Pink and yellow curves depict the nullclines of system (\ref{eq:4}), intersecting at equilibria where elk and wolf growth rates simultaneously vanish.

This stability structure reflects the ecological dynamics documented in Yellowstone National Park. The stable equilibrium represents a balanced state where elk and wolves coexist, sustained by regulated predation and available resources. The left saddle point may correspond to a scenario where elk densities fall too low, making the population vulnerable to extinction due to sustained wolf predation. For instance, in the Madison headwaters region, elk densities below approximately \(4\ \mathrm{elk/km}^2\)
  were linked to wolf pack dissolution or dispersal due to inefficient hunting \cite{becker2008wolf}. On the other hand, the right saddle equilibrium may reflect historical overabundance of elk, observed at densities exceeding \(20\ \mathrm{elk/km}^2\) in Yellowstone’s northern range \cite{houston1982northern,lemke1998winter}. Such overpopulation caused intensive winter browsing, leading to the decline of key vegetation species like willow and aspen, ultimately destabilizing the ecosystem \cite{ripple2012trophic}. The central stable equilibrium represents the observed prey-predator balance after wolf reintroduction \cite{ripple2012trophic,smith2003yellowstone}. 

\begin{figure}[H]
	\begin{subfigure}{0.32\textwidth} 
		\centering
		\includegraphics[width= \textwidth, height=3.3cm]{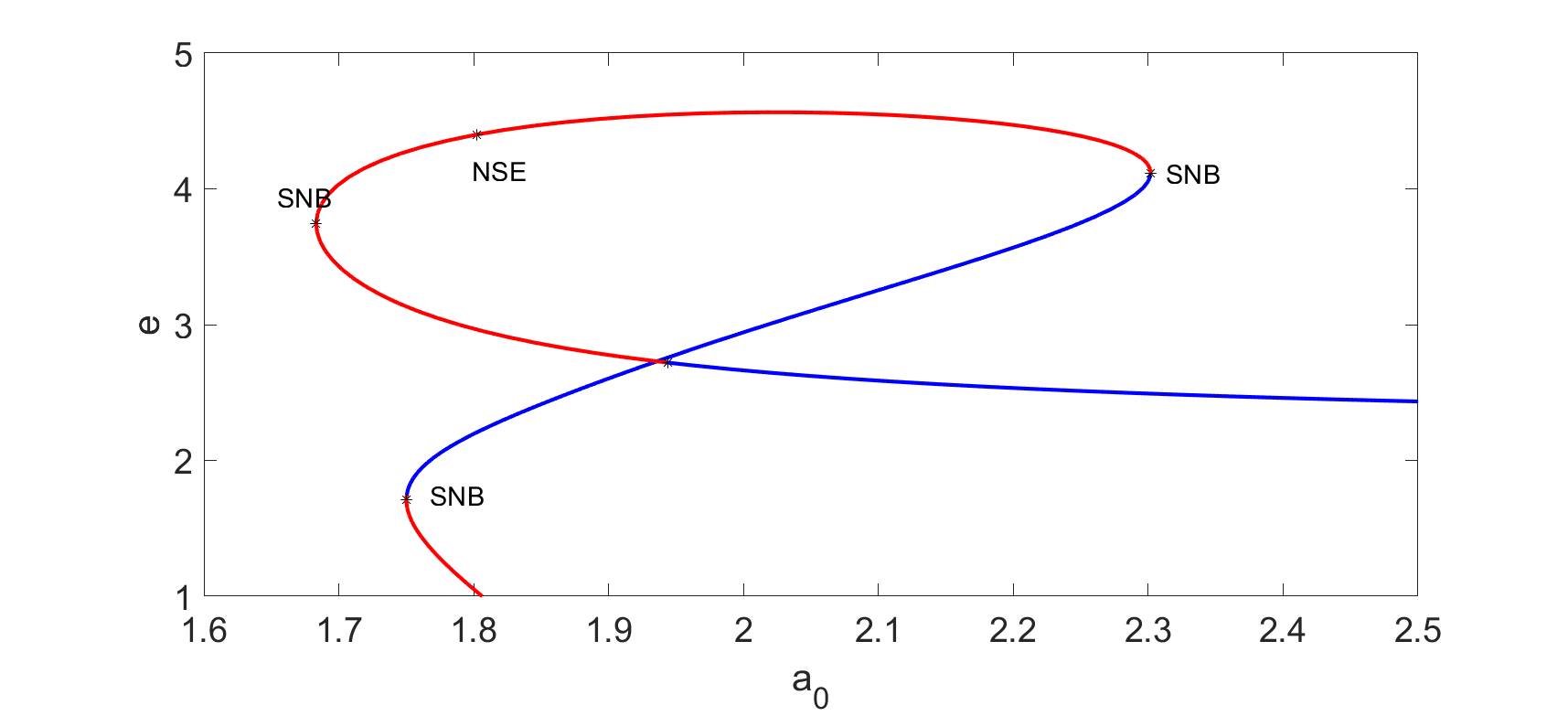} 
		\caption{\centering}
	\end{subfigure} 
	\begin{subfigure}{0.32\textwidth} 
		\centering
		\includegraphics[width= \textwidth, height=3.3cm]{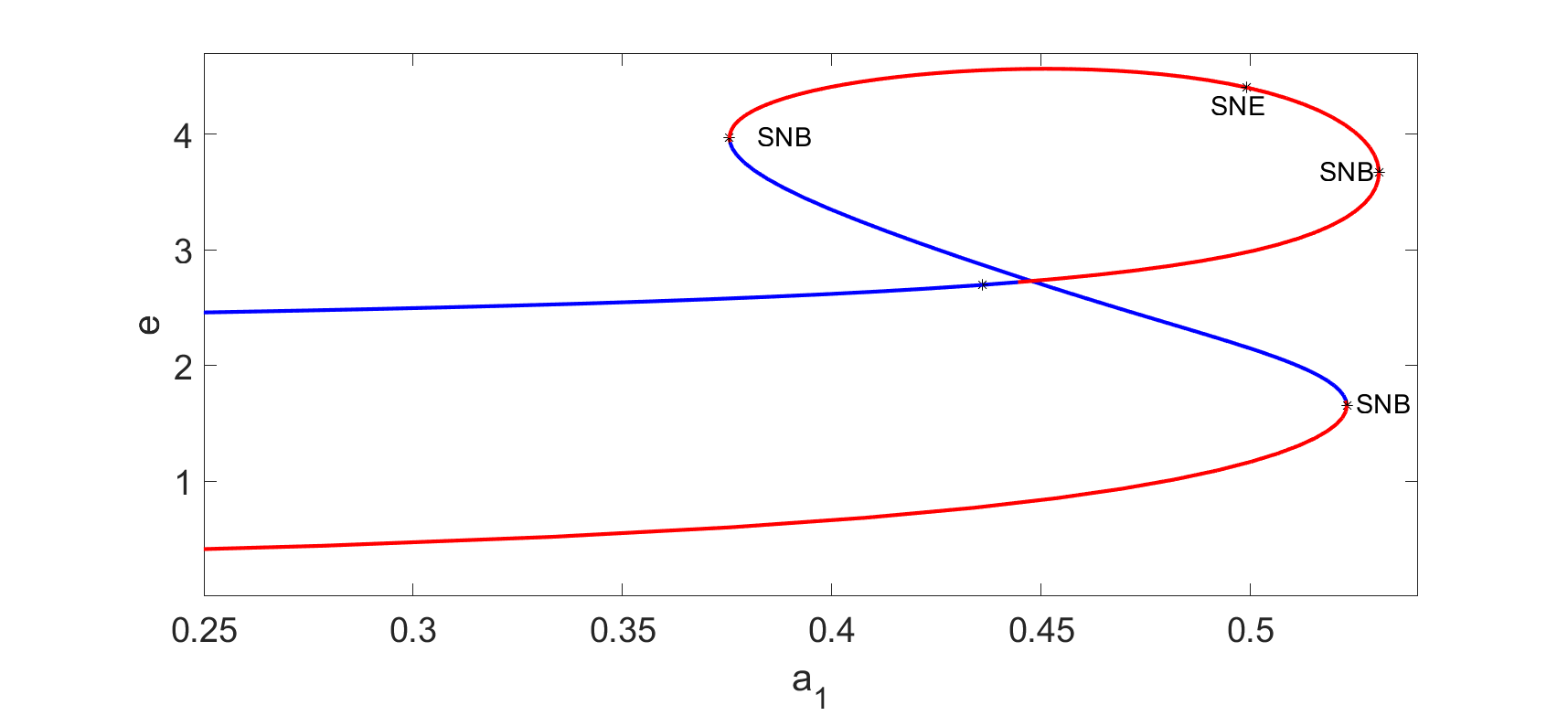} 
		\caption{\centering}
	\end{subfigure} 
	\begin{subfigure}{0.32\textwidth} 
		\centering
		\includegraphics[width= \textwidth, height=3.3cm]{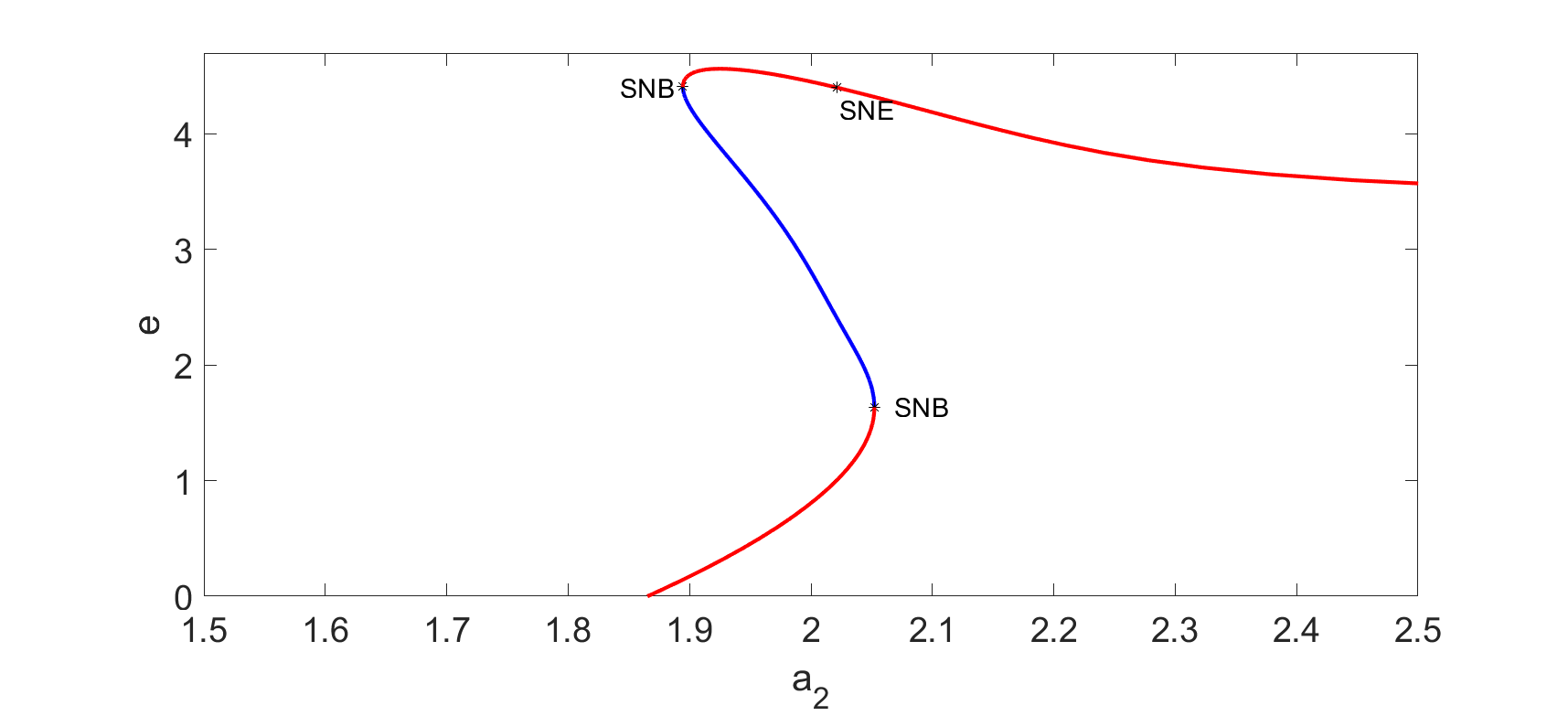} 
		\caption{\centering}
	\end{subfigure} 
	\begin{subfigure}{0.32\textwidth} 
		\centering
		\includegraphics[width= \textwidth, height=3.3cm]{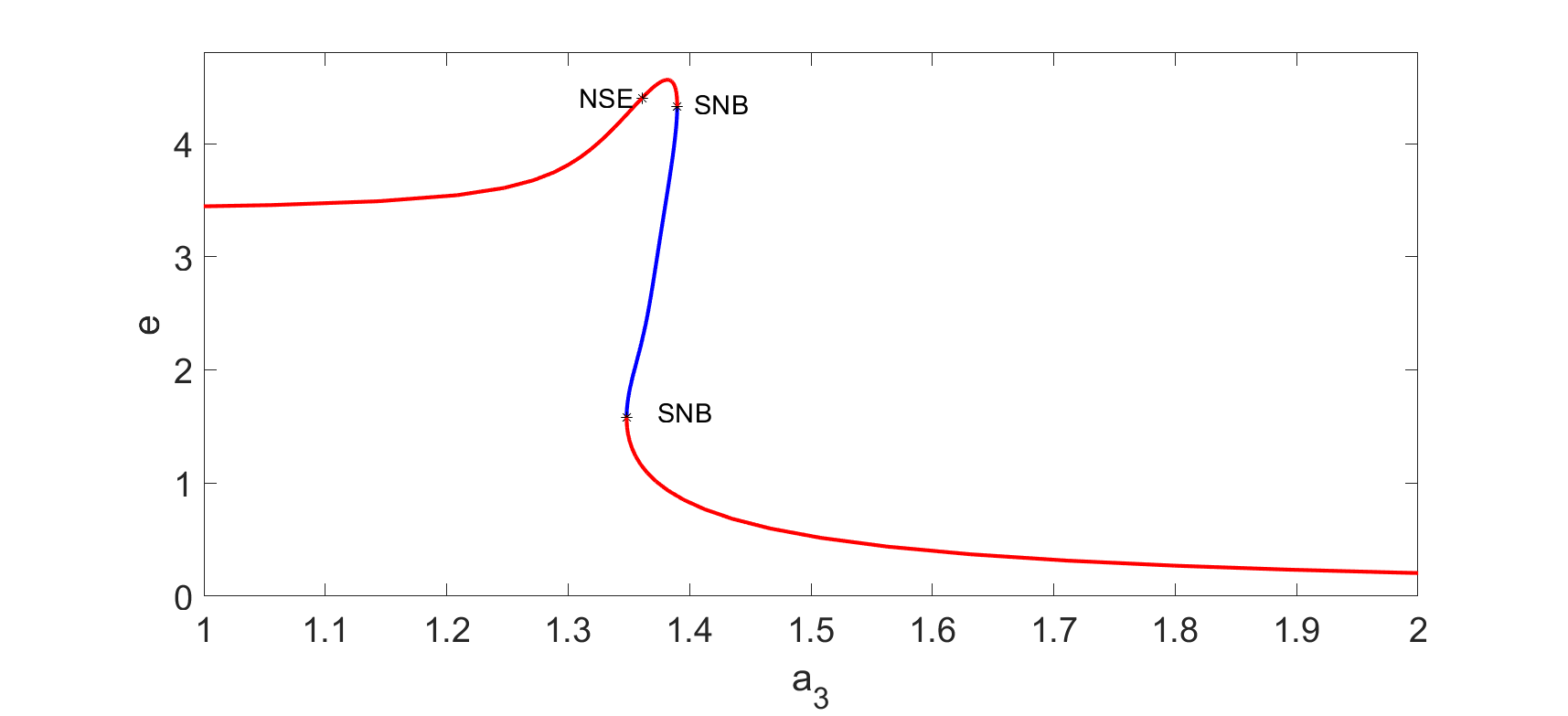} 
		\caption{\centering}
	\end{subfigure} 
	\begin{subfigure}{0.32\textwidth} 
		\centering
		\includegraphics[width= \textwidth, height=3.3cm]{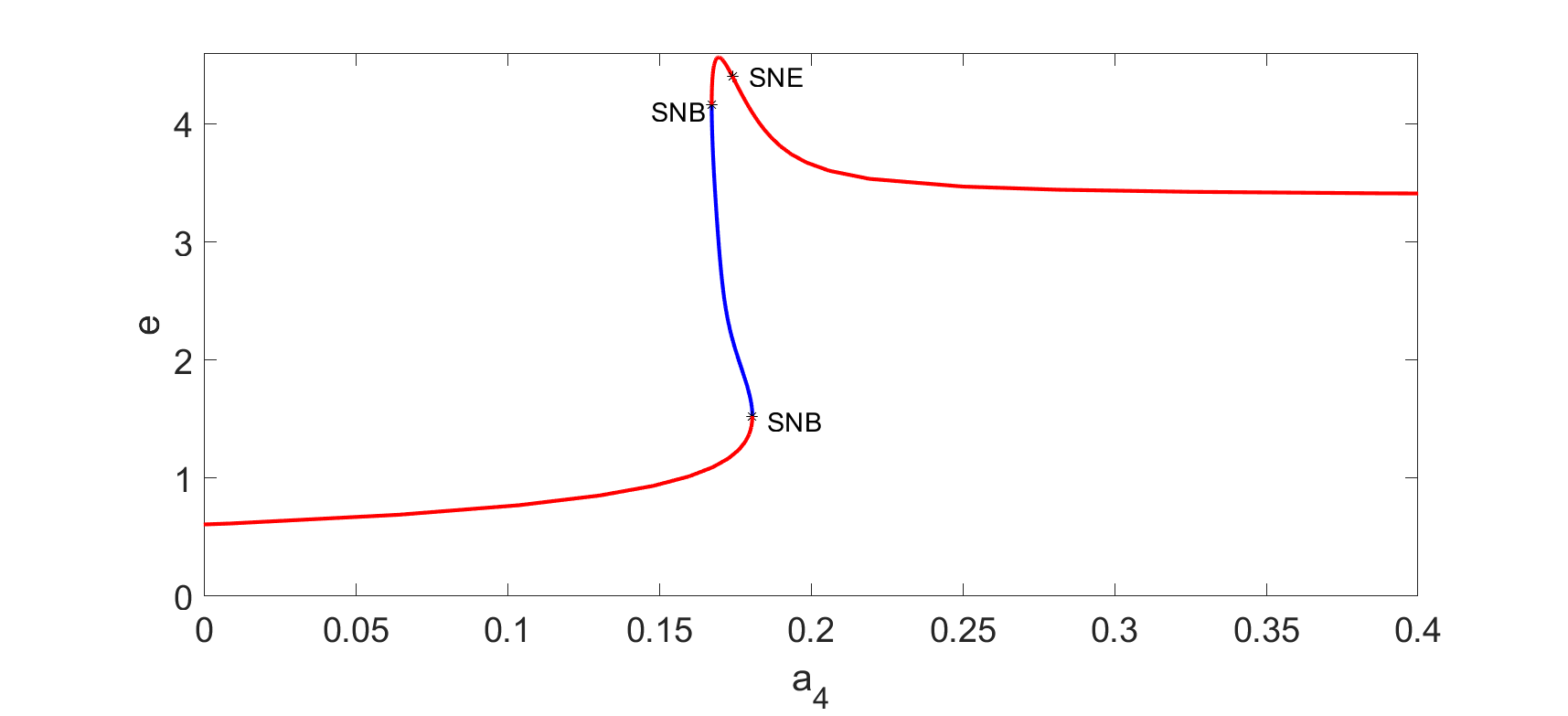} 
		\caption{\centering}
	\end{subfigure} 
	\begin{subfigure}{0.32\textwidth} 
		\centering
		\includegraphics[width= \textwidth, height=3.3cm]{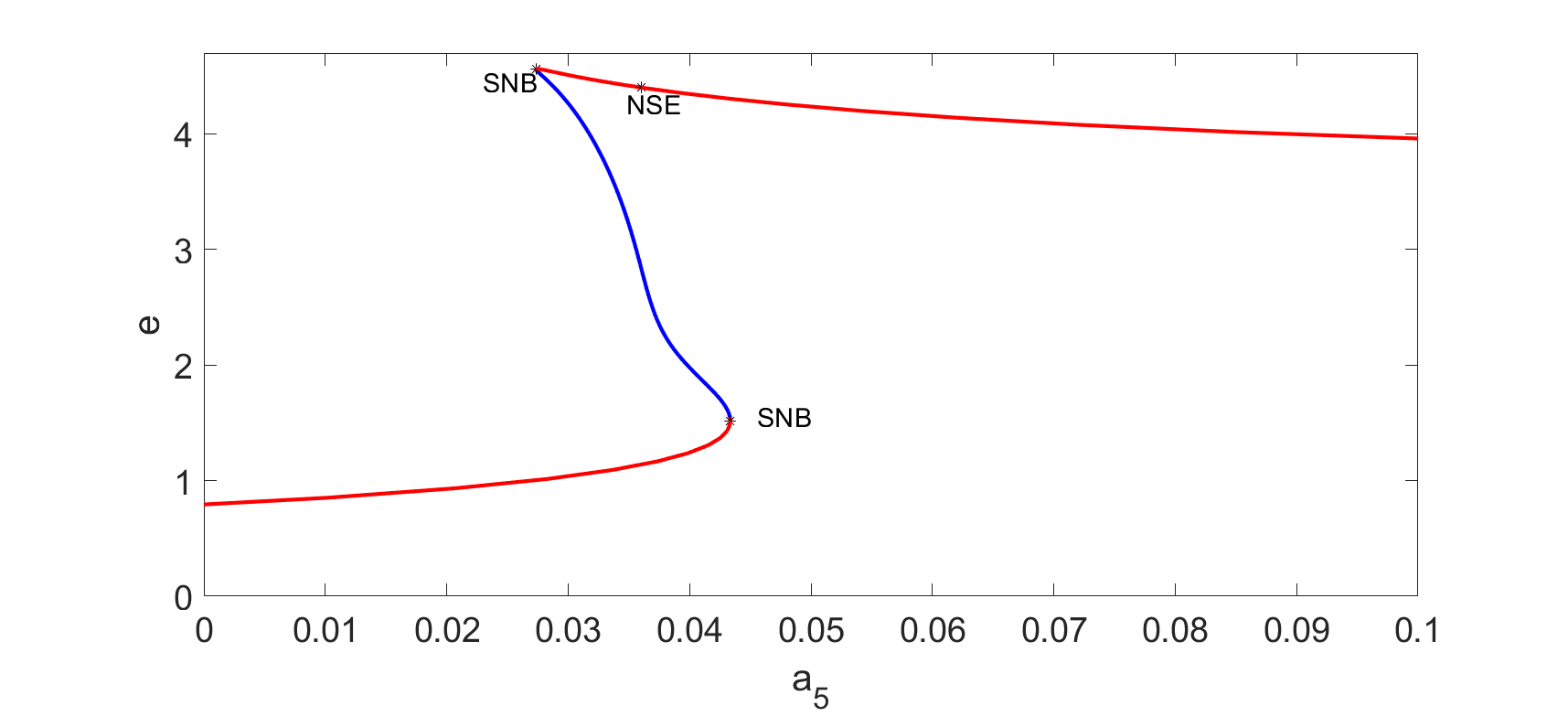} 
		\caption{\centering}
	\end{subfigure} 
	\begin{subfigure}{0.32\textwidth} 
		\centering
		\includegraphics[width= \textwidth, height=3.3cm]{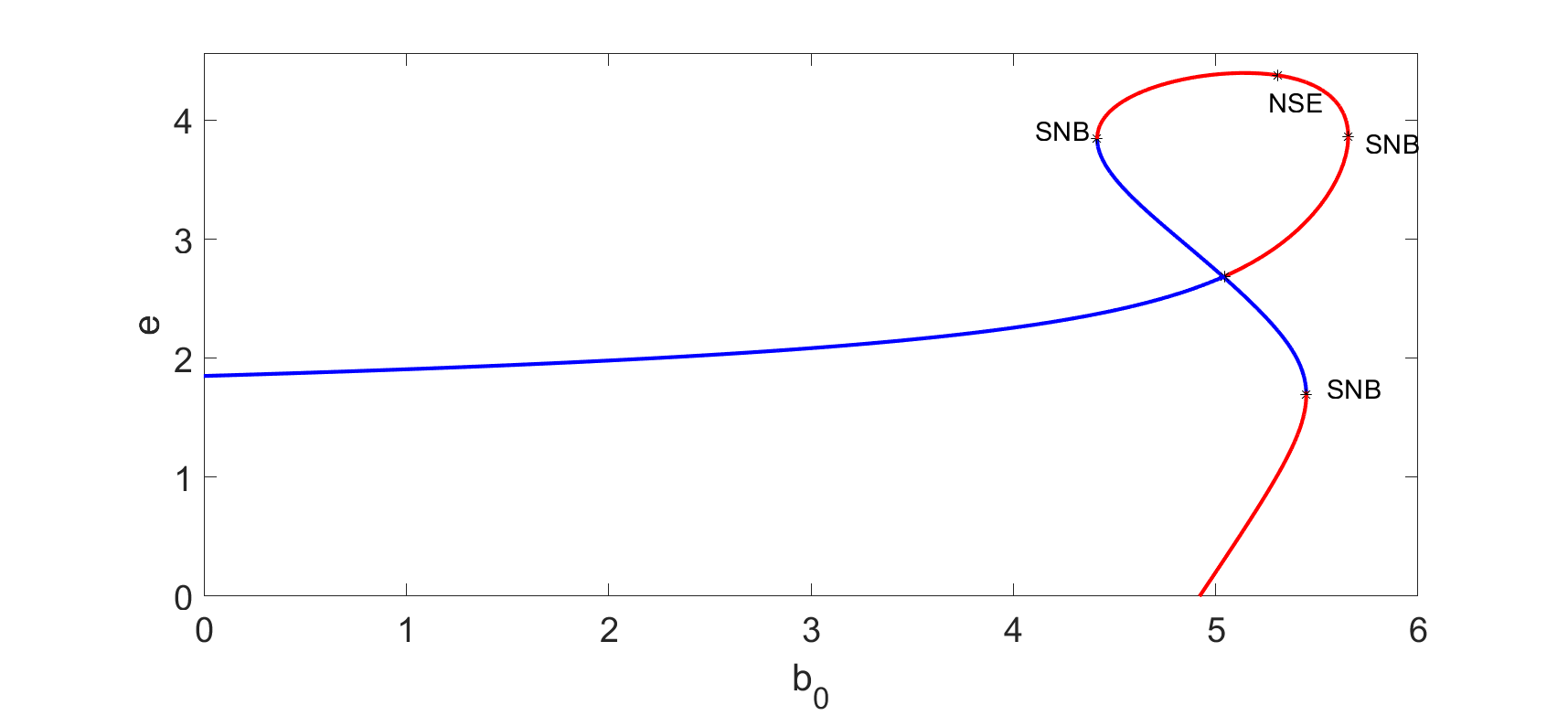} 
		\caption{\centering}
	\end{subfigure} 
	\begin{subfigure}{0.32\textwidth} 
		\centering
		\includegraphics[width= \textwidth, height=3.3cm]{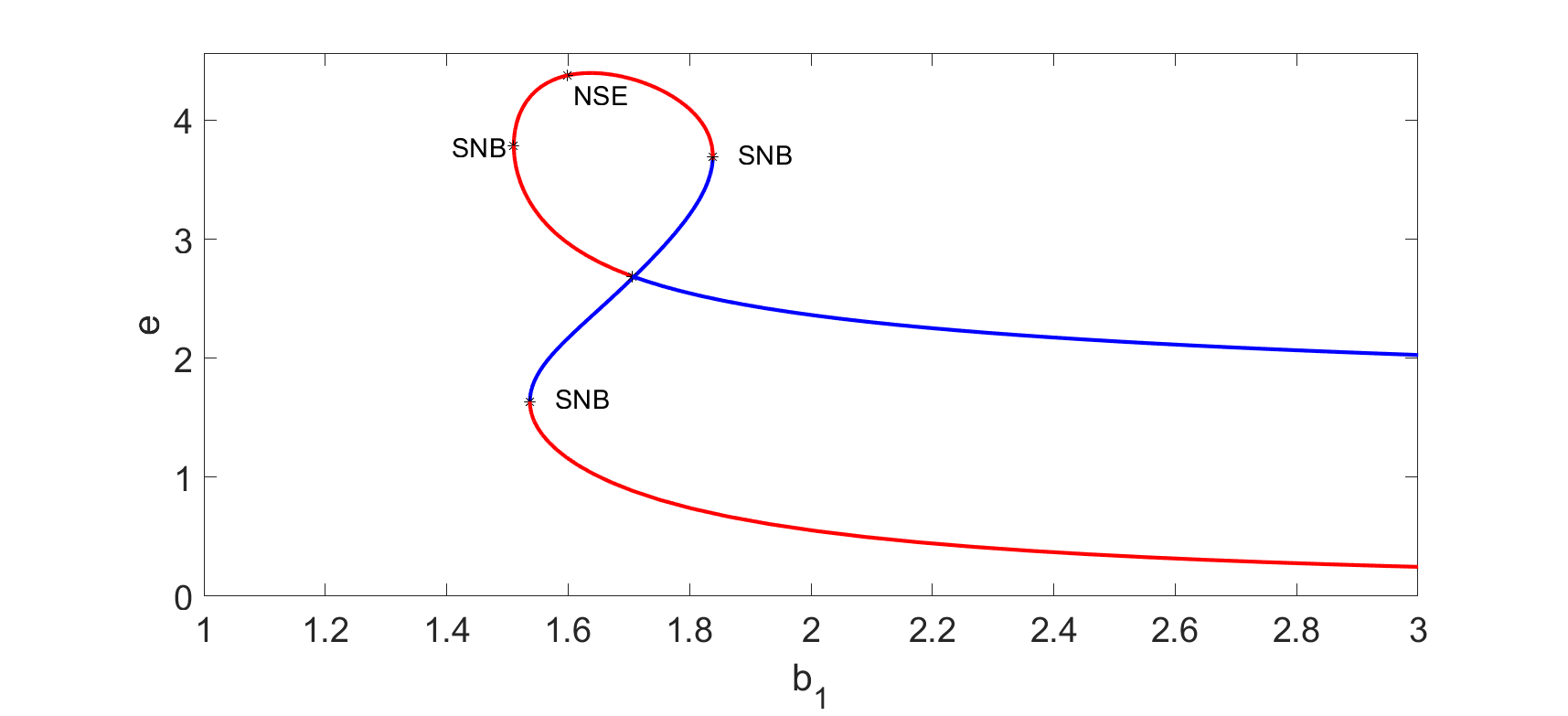} 
		\caption{\centering}
	\end{subfigure} 
	\begin{subfigure}{0.32\textwidth} 
		\centering
		\includegraphics[width= \textwidth, height=3.3cm]{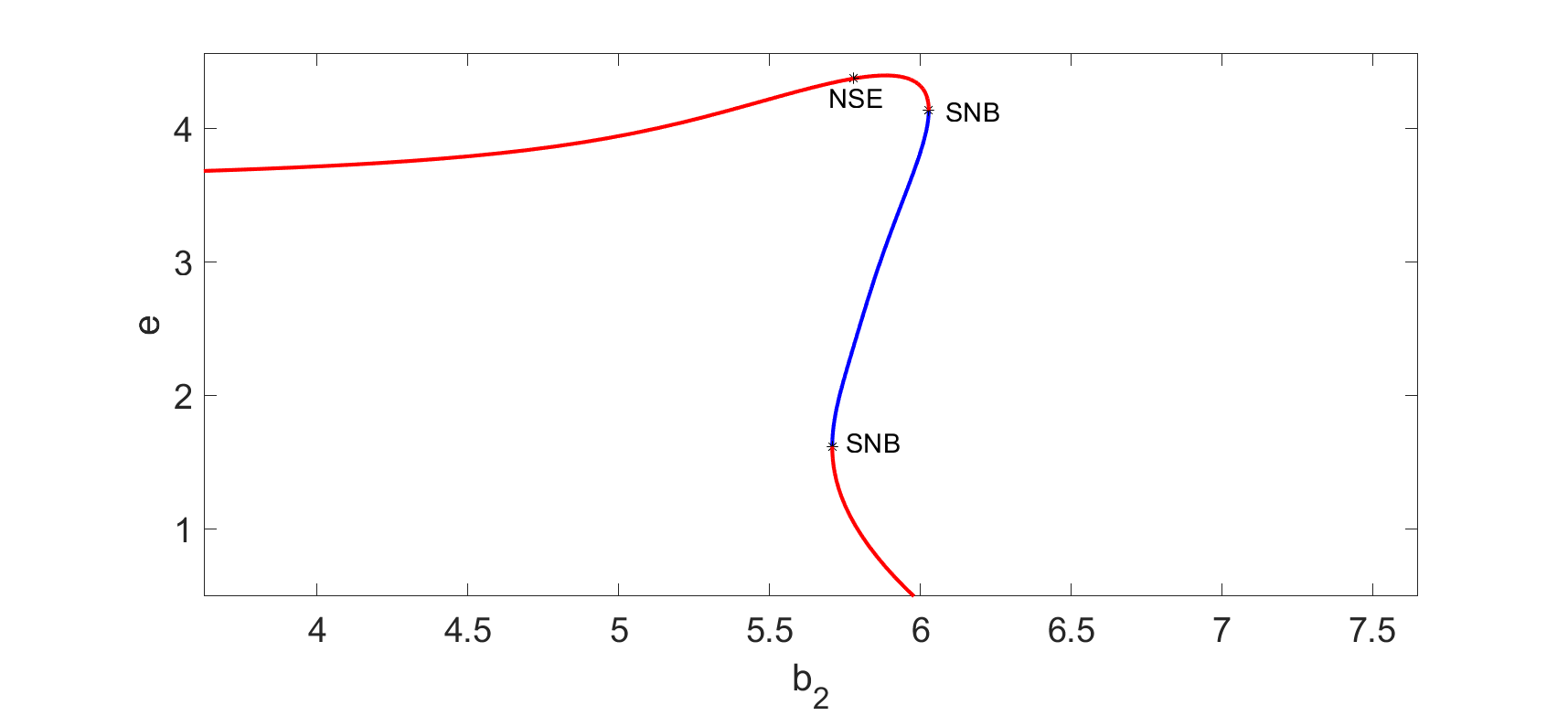} 
		\caption{\centering}
	\end{subfigure} 
	\begin{subfigure}{0.32\textwidth} 
		\centering
		\includegraphics[width= \textwidth, height=3.3cm]{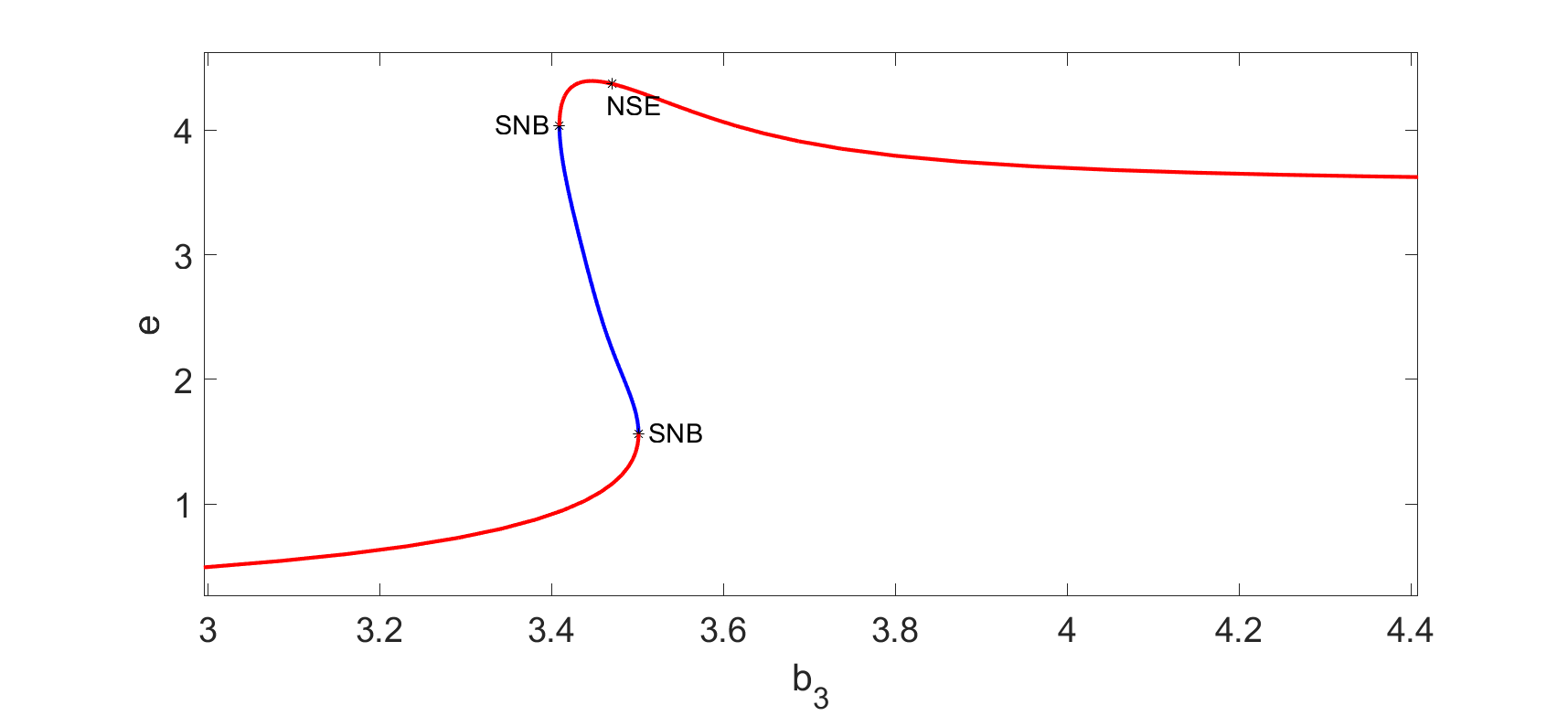} 
		\caption{\centering}
	\end{subfigure} 
	\begin{subfigure}{0.32\textwidth} 
		\centering
		\includegraphics[width= \textwidth, height=3.3cm]{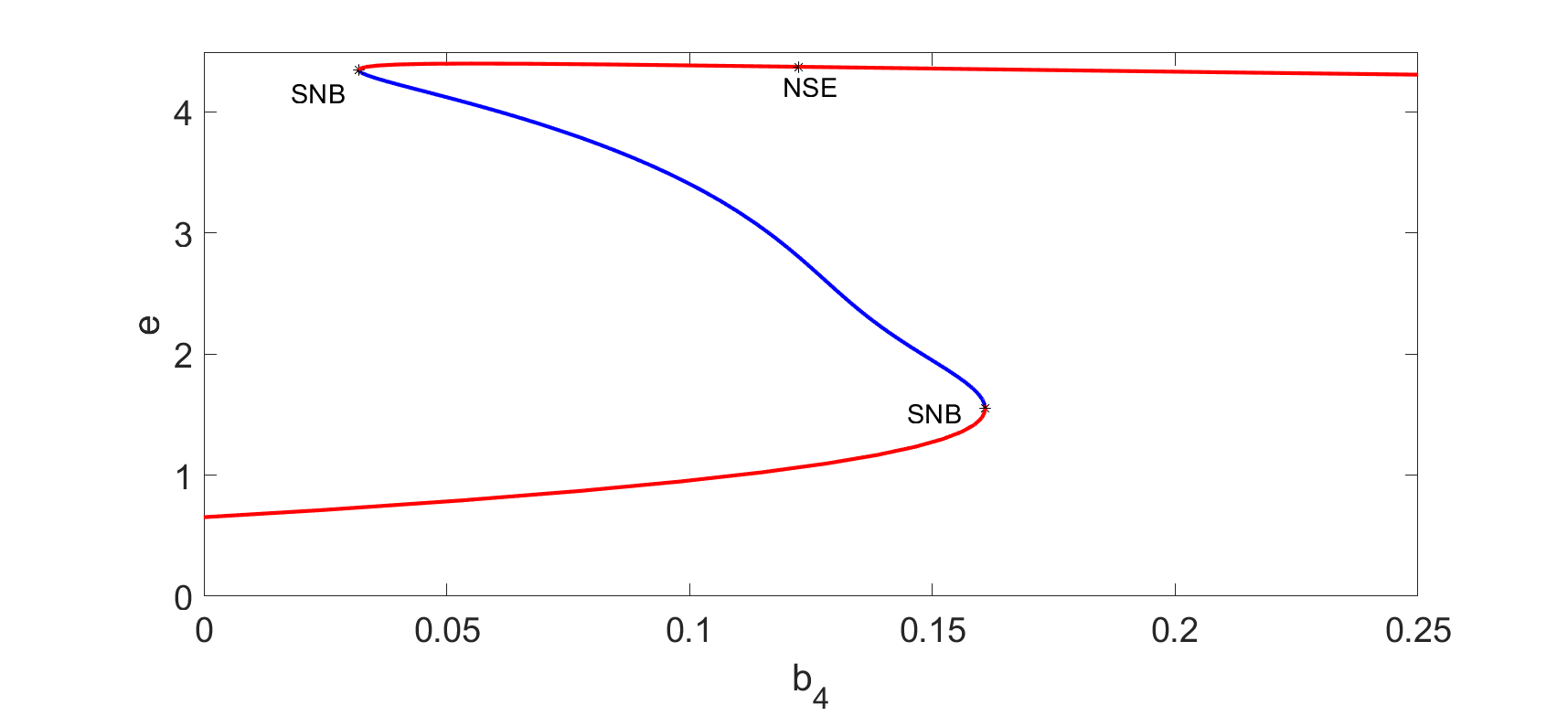} 
		\caption{\centering}
	\end{subfigure} 
	\begin{subfigure}{0.32\textwidth} 
		\centering
		\includegraphics[width= \textwidth, height=3.3cm]{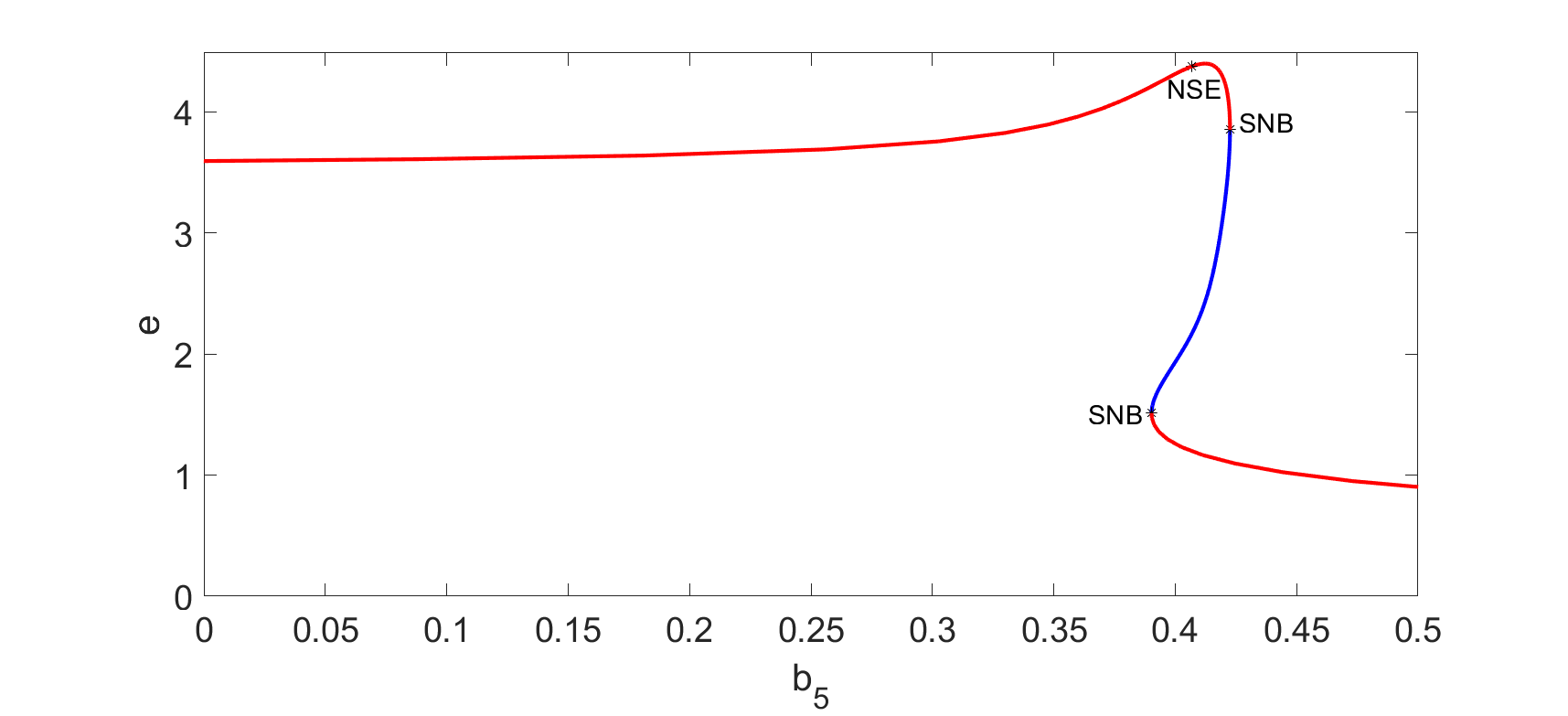} 
		\caption{\centering}
	\end{subfigure} 
	\begin{subfigure}{0.32\textwidth} 
		\centering
		\includegraphics[width= \textwidth, height=3.3cm]{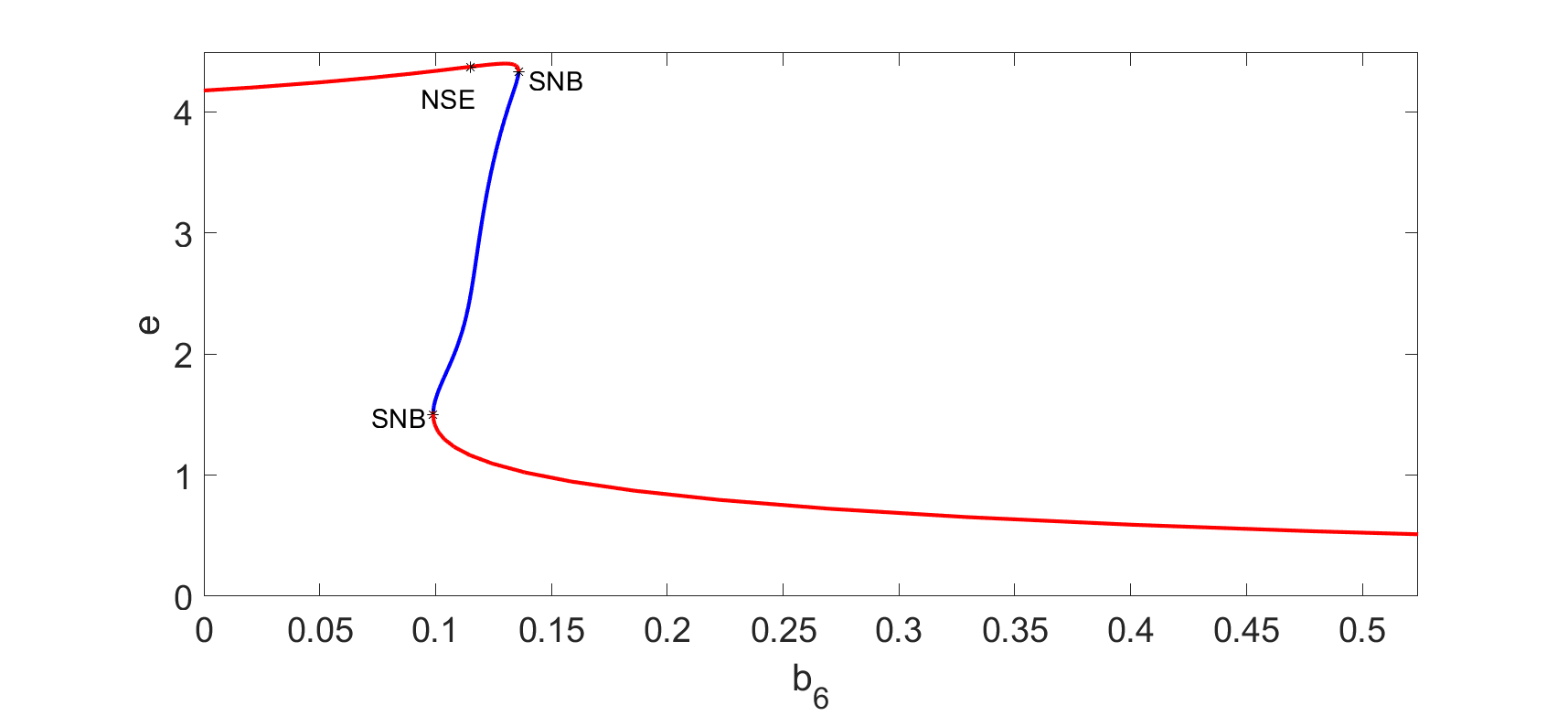} 
		\caption{\centering}
	\end{subfigure} 
	\begin{subfigure}{0.32\textwidth} 
		\centering
		\includegraphics[width= \textwidth, height=3.3cm]{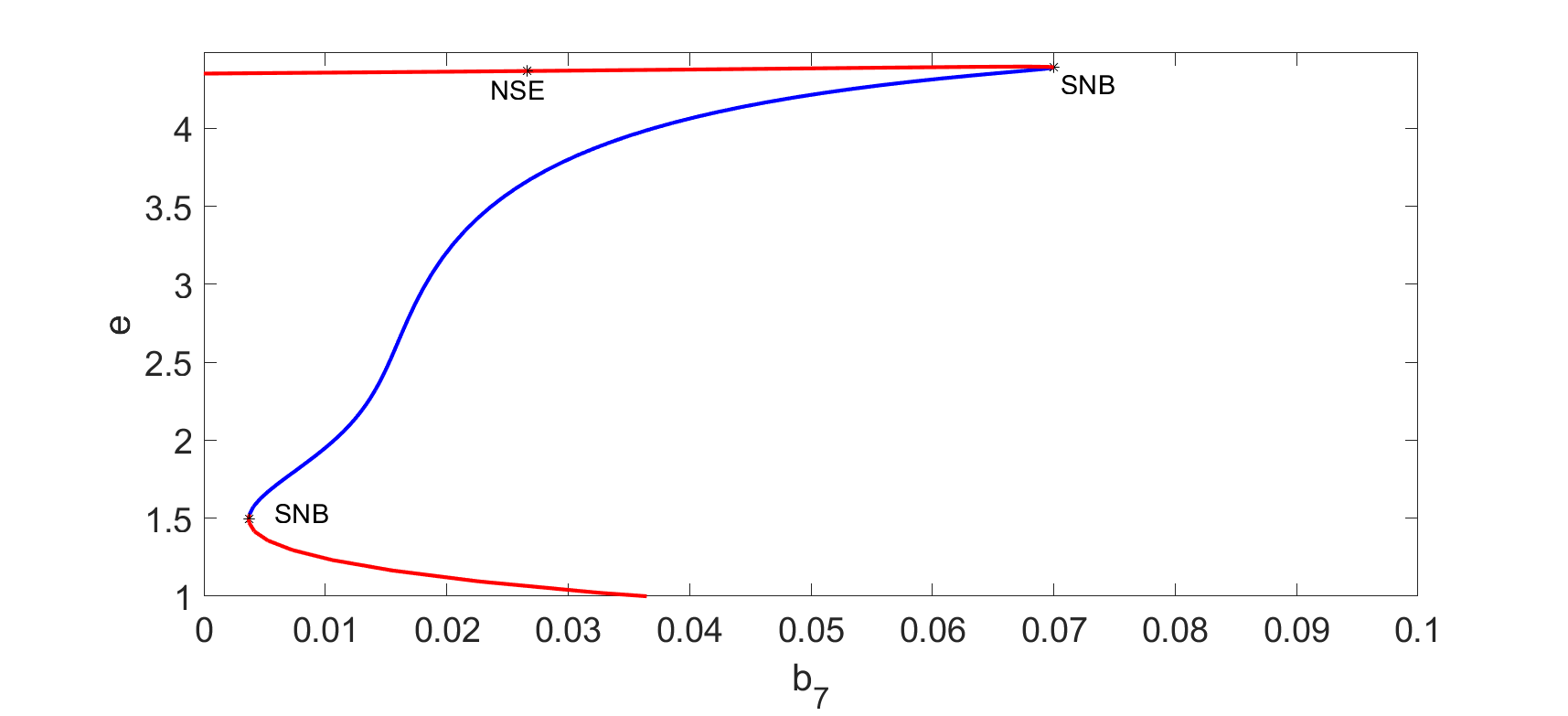} 
		\caption{\centering}
	\end{subfigure} 
	\caption{Co-dimension one bifurcation structures for each parameters of the model (\ref{eq:4}) with respect to elk population. In these figures SNB denotes the saddle node bifurcation, NSE stands for neutral saddle equilibrium.}
	\label{Fig:elk_bifurcation}   
\end{figure}

\begin{figure}[H]
	\begin{subfigure}{0.32\textwidth} 
		\centering
		\includegraphics[width=\textwidth, height=3.3cm]{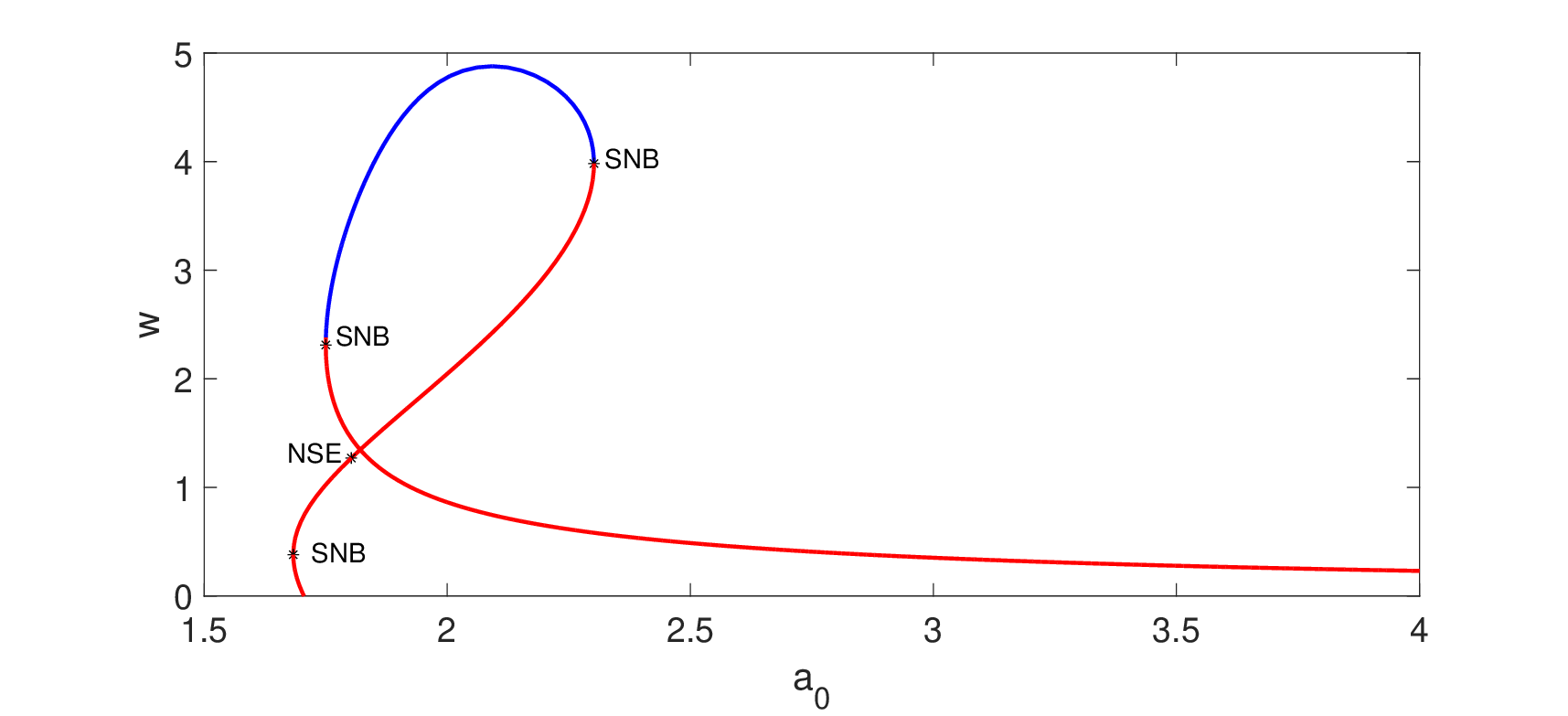} 
		\caption{\centering}
	\end{subfigure} 
	\begin{subfigure}{0.32\textwidth} 
		\centering
		\includegraphics[width= \textwidth, height=3.3cm]{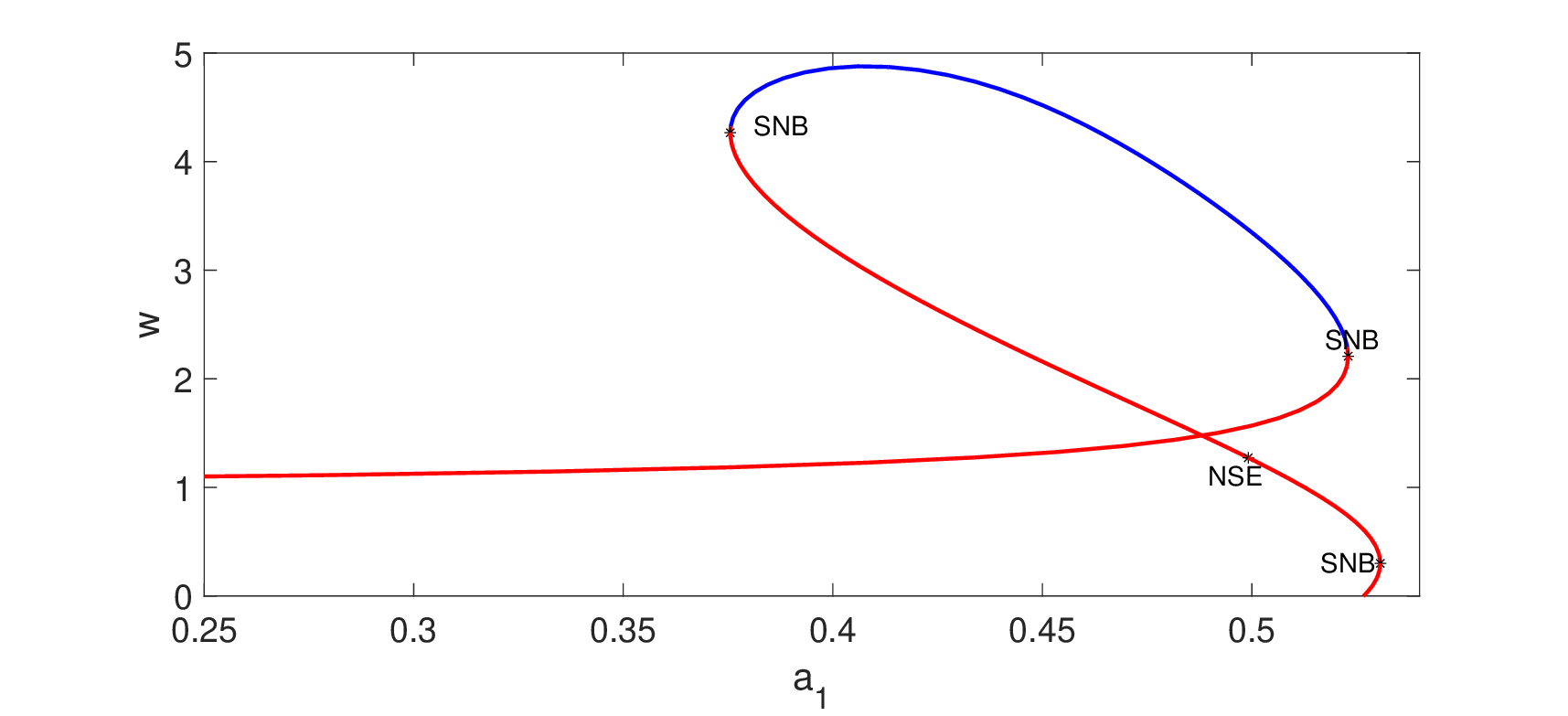} 
		\caption{\centering}
	\end{subfigure} 
	\begin{subfigure}{0.32\textwidth} 
		\centering
		\includegraphics[width= \textwidth, height=3.3cm]{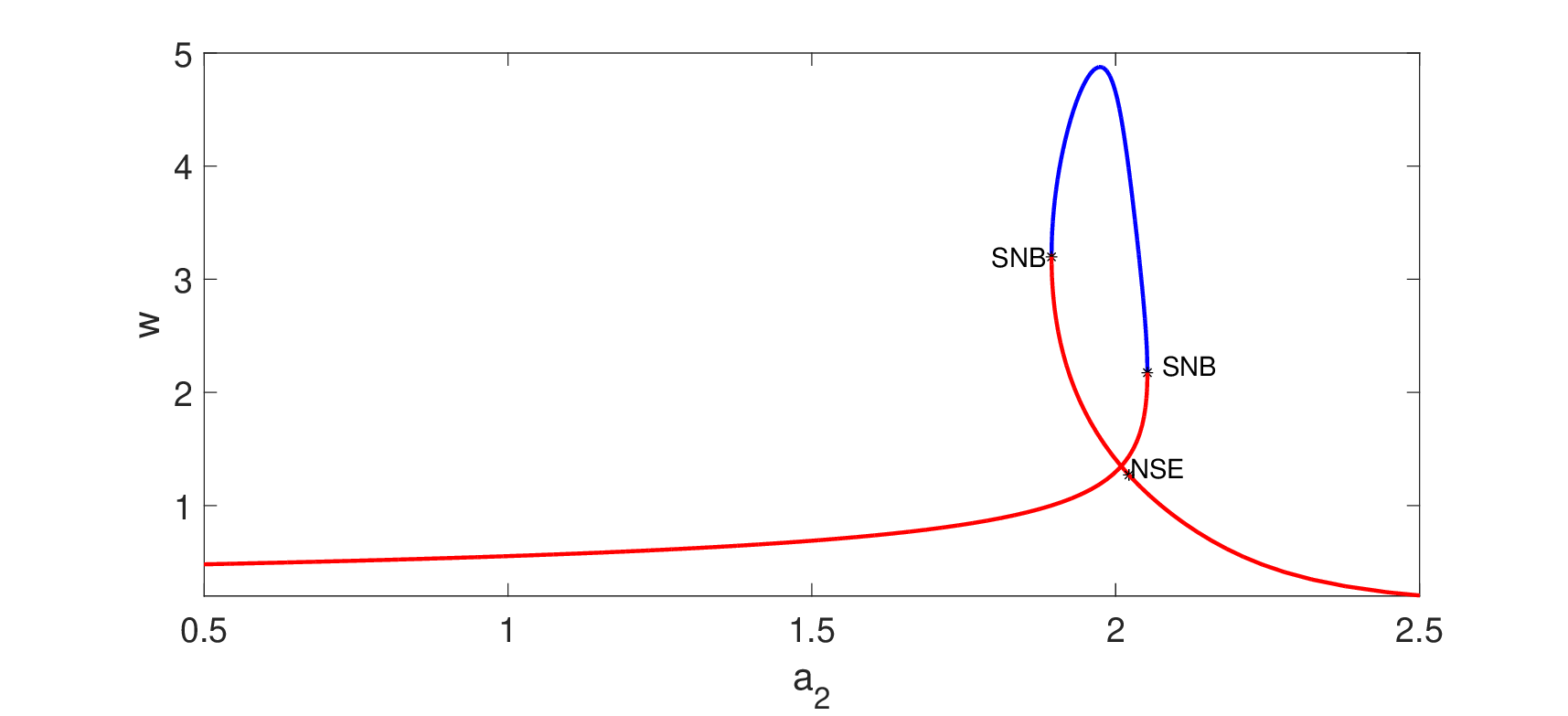} 
		\caption{\centering}
	\end{subfigure} 
	\begin{subfigure}{0.32\textwidth} 
		\centering
		\includegraphics[width= \textwidth, height=3.3cm]{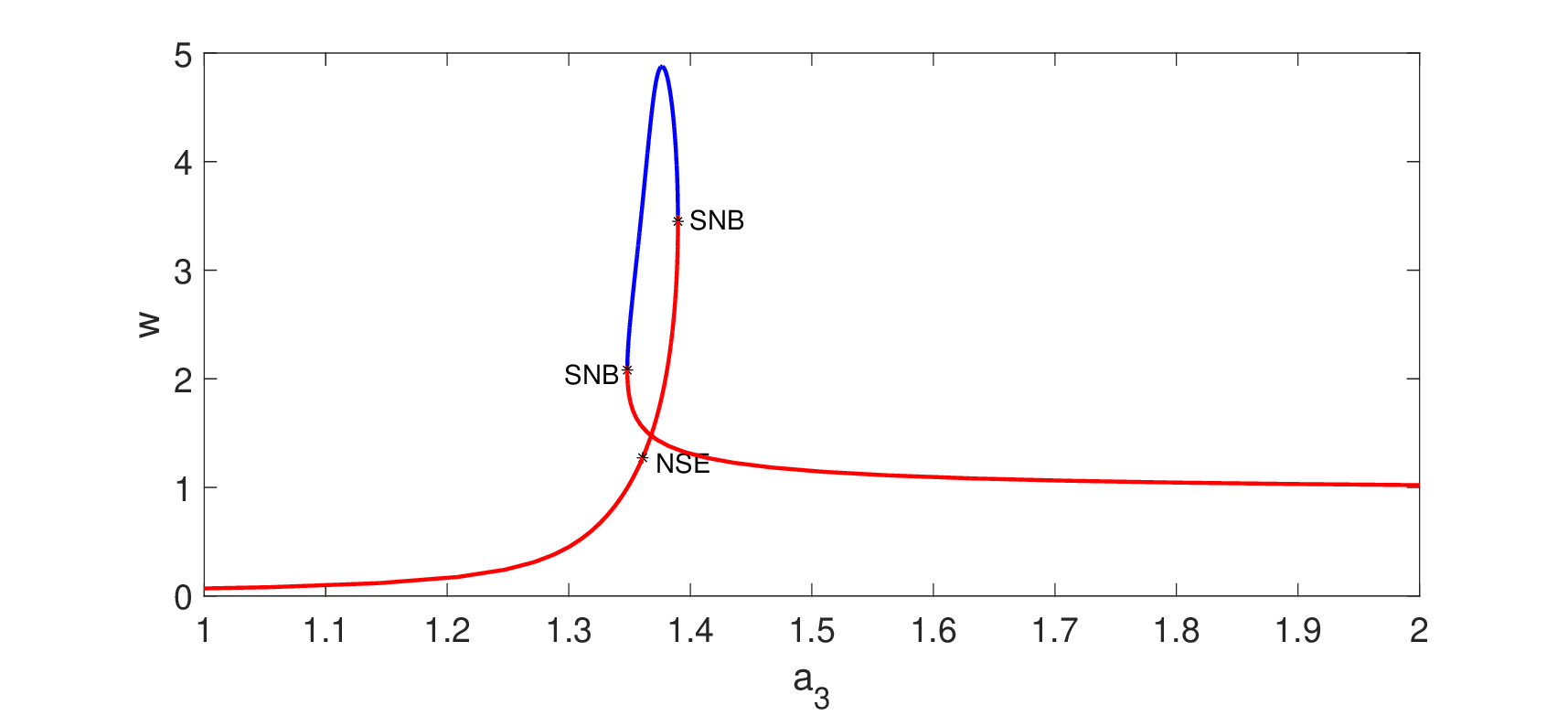} 
		\caption{\centering}
	\end{subfigure} 
	\begin{subfigure}{0.32\textwidth} 
		\centering
		\includegraphics[width= \textwidth, height=3.3cm]{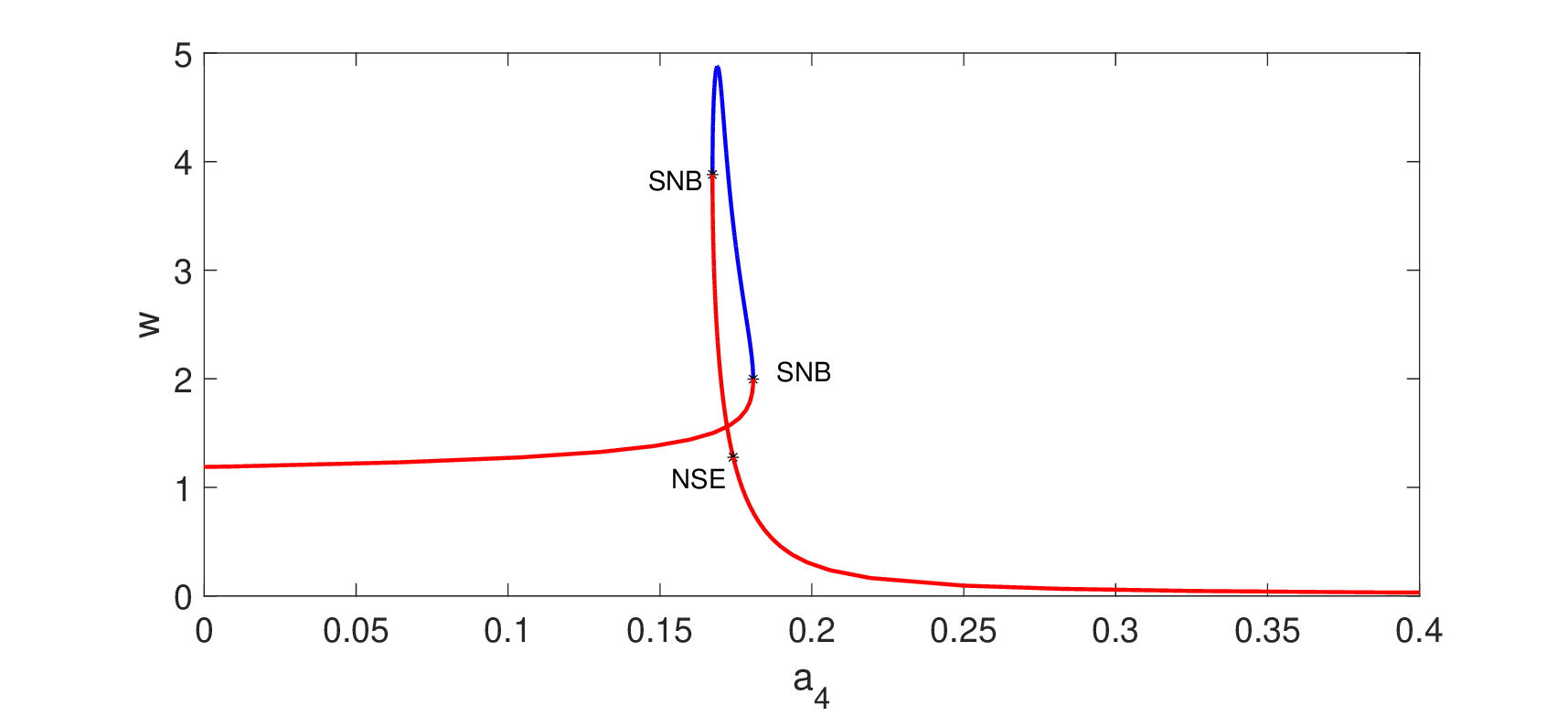} 
		\caption{\centering}
	\end{subfigure} 
	\begin{subfigure}{0.32\textwidth} 
		\centering
		\includegraphics[width= \textwidth, height=3.3cm]{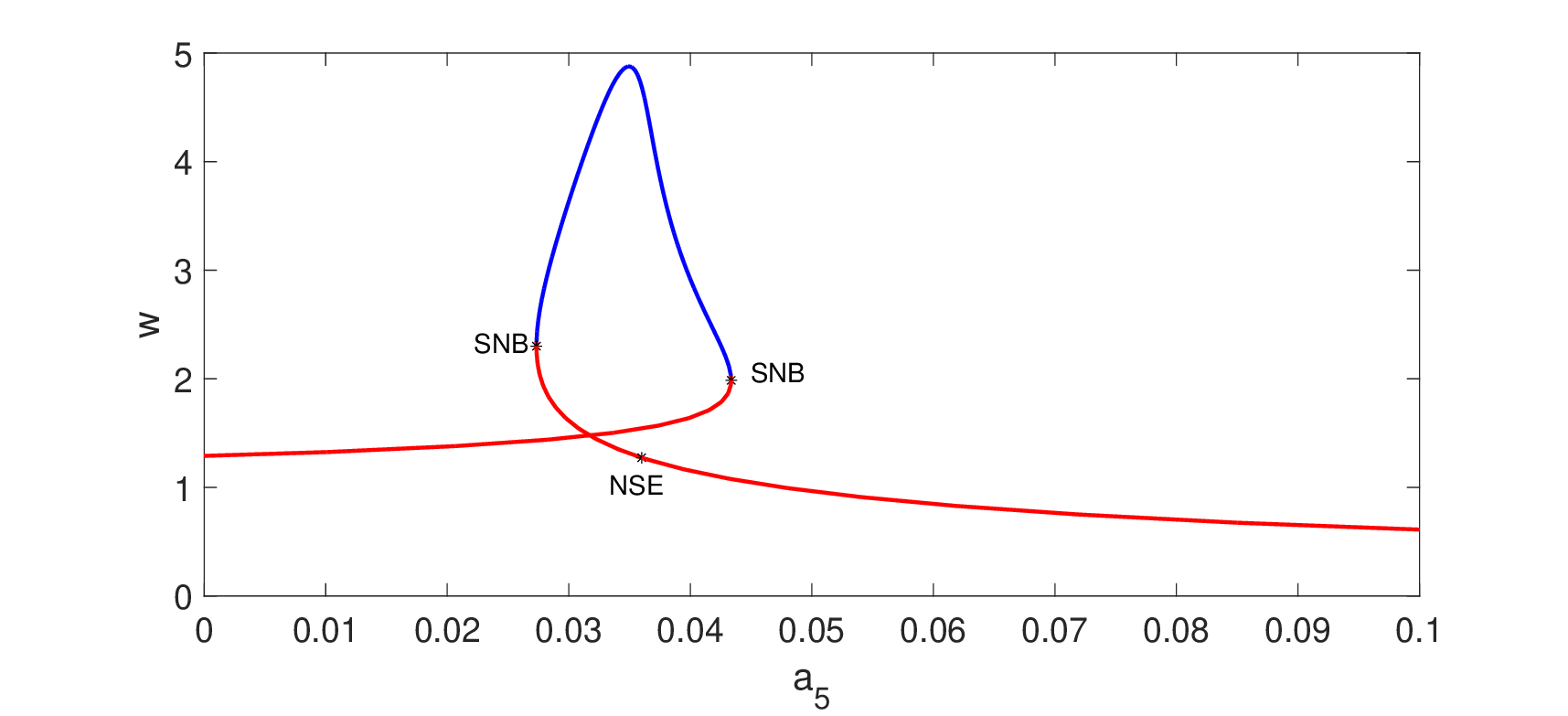} 
		\caption{\centering}
	\end{subfigure} 
	\begin{subfigure}{0.32\textwidth} 
		\centering
		\includegraphics[width= \textwidth, height=3.3cm]{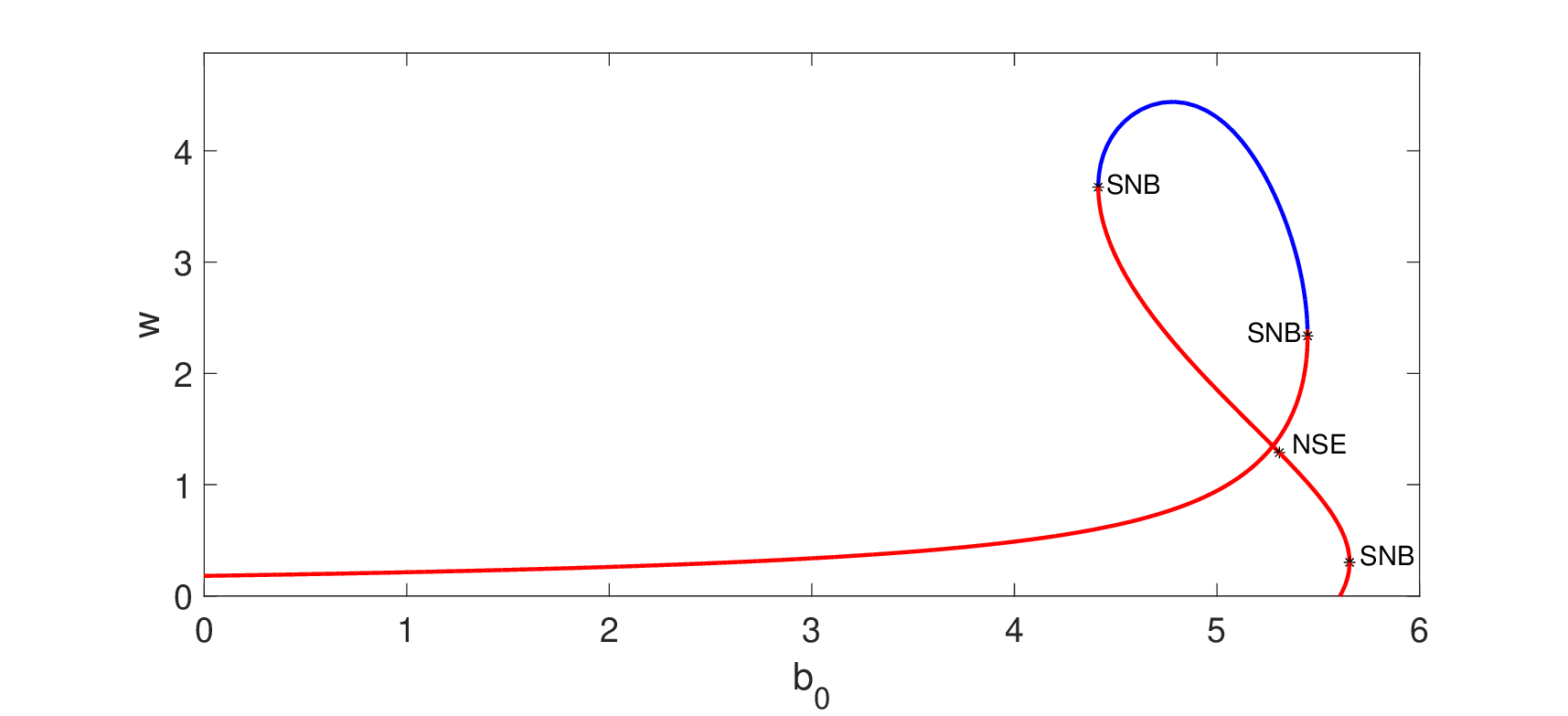} 
		\caption{\centering}
	\end{subfigure} 
	\begin{subfigure}{0.32\textwidth} 
		\centering
		\includegraphics[width= \textwidth, height=3.3cm]{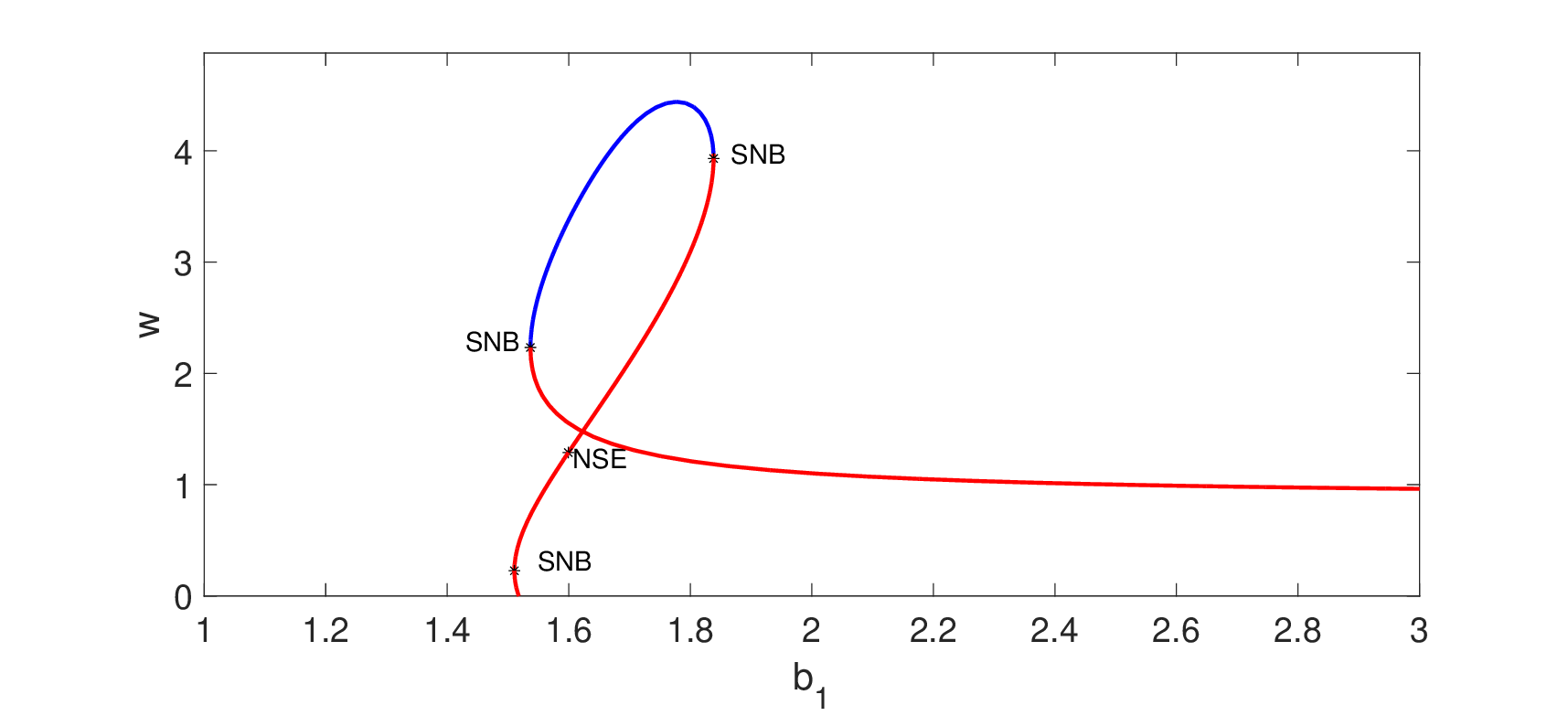} 
		\caption{\centering}
	\end{subfigure} 
	\begin{subfigure}{0.32\textwidth} 
		\centering
		\includegraphics[width= \textwidth, height=3.3cm]{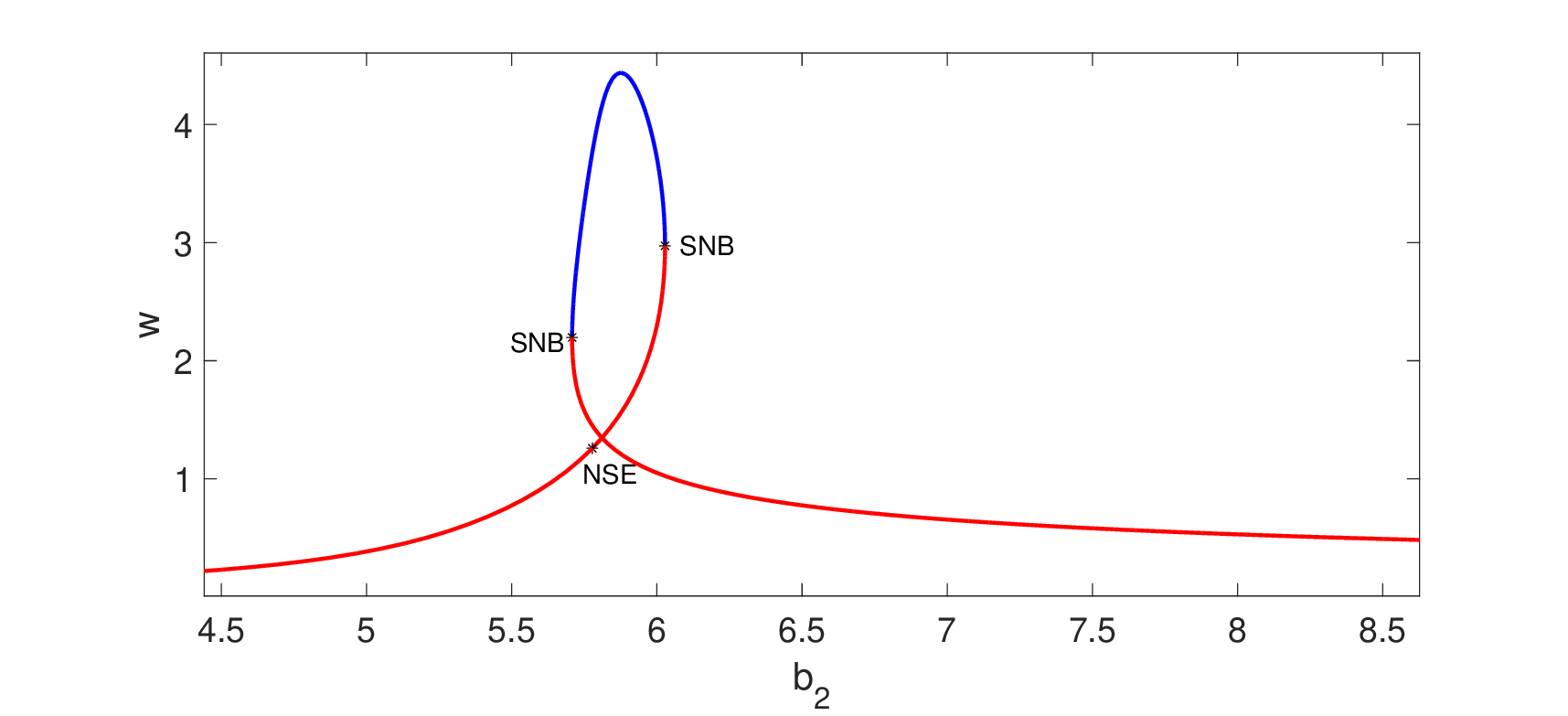} 
		\caption{\centering}
	\end{subfigure} 
	\begin{subfigure}{0.32\textwidth} 
		\centering
		\includegraphics[width= \textwidth, height=3.3cm]{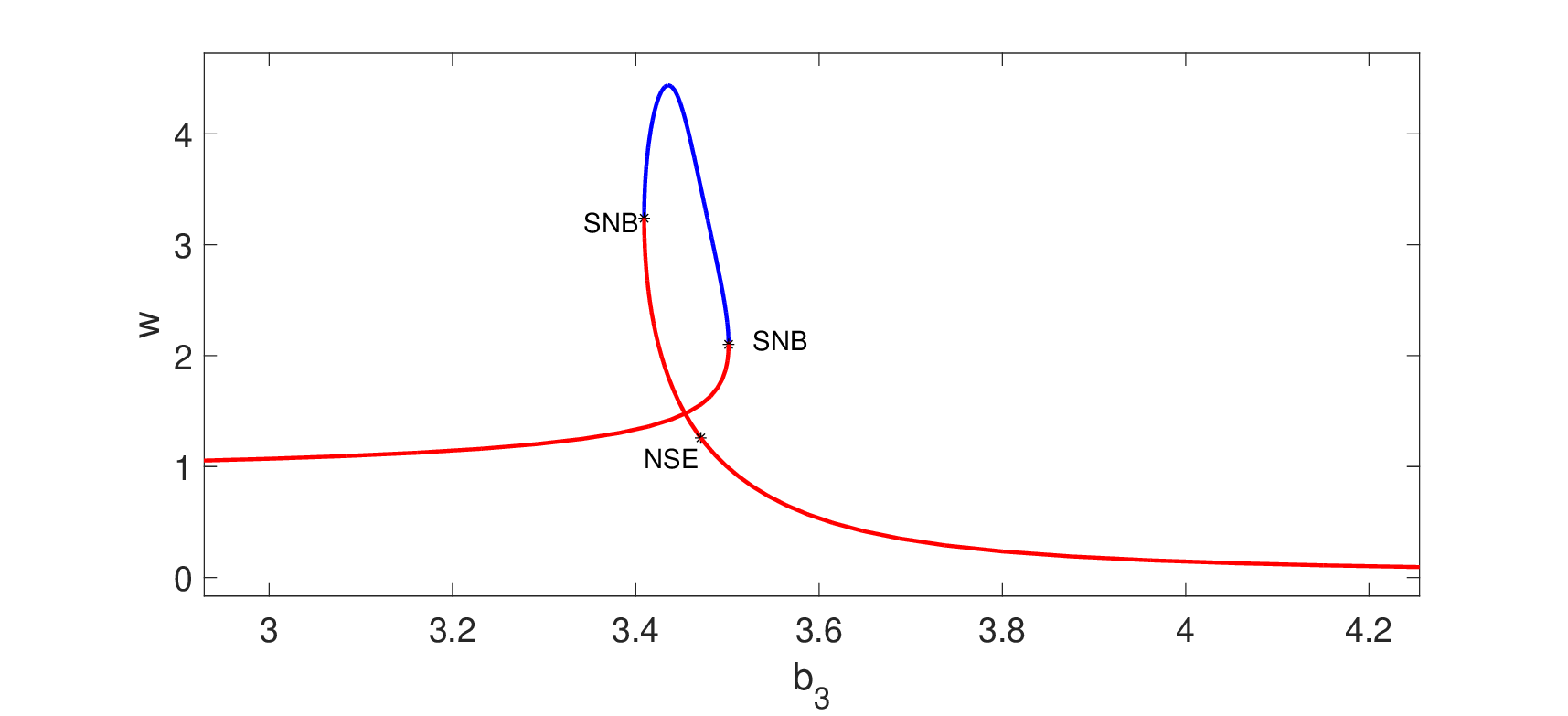} 
		\caption{\centering}
	\end{subfigure} 
	\begin{subfigure}{0.32\textwidth} 
		\centering
		\includegraphics[width= \textwidth, height=3.3cm]{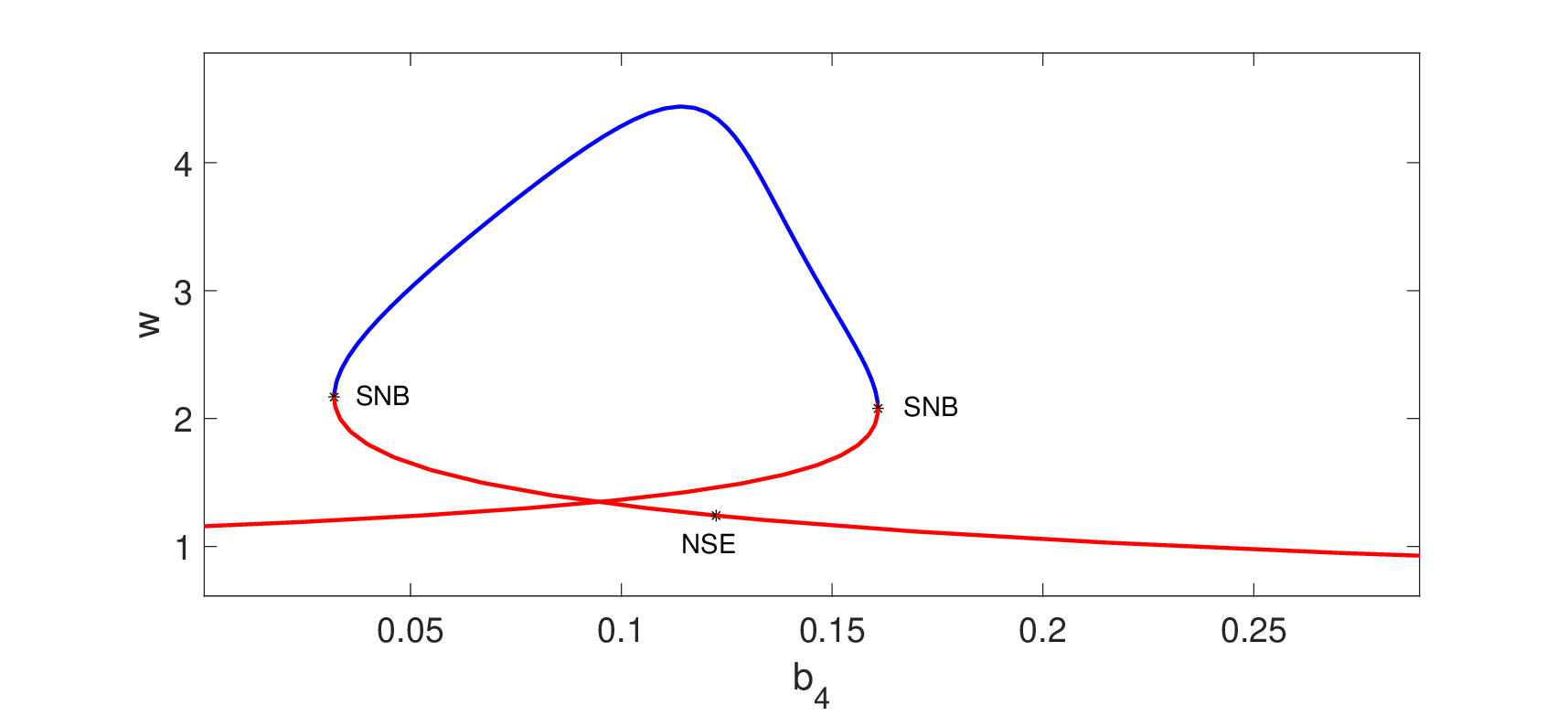} 
		\caption{\centering}
	\end{subfigure} 
	\begin{subfigure}{0.32\textwidth} 
		\centering
		\includegraphics[width= \textwidth, height=3.3cm]{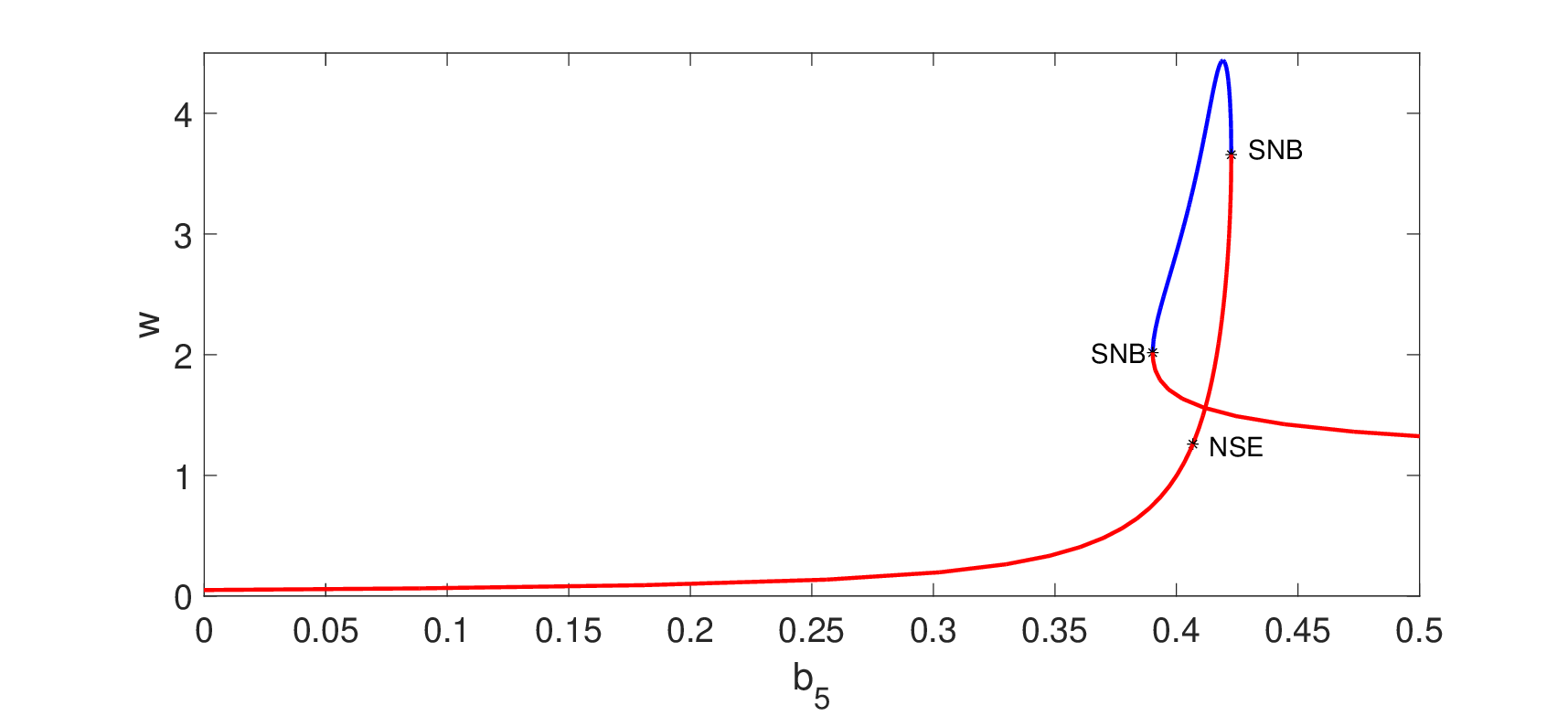} 
		\caption{\centering}
	\end{subfigure} 
	\begin{subfigure}{0.32\textwidth} 
		\centering
		\includegraphics[width= \textwidth, height=3.3cm]{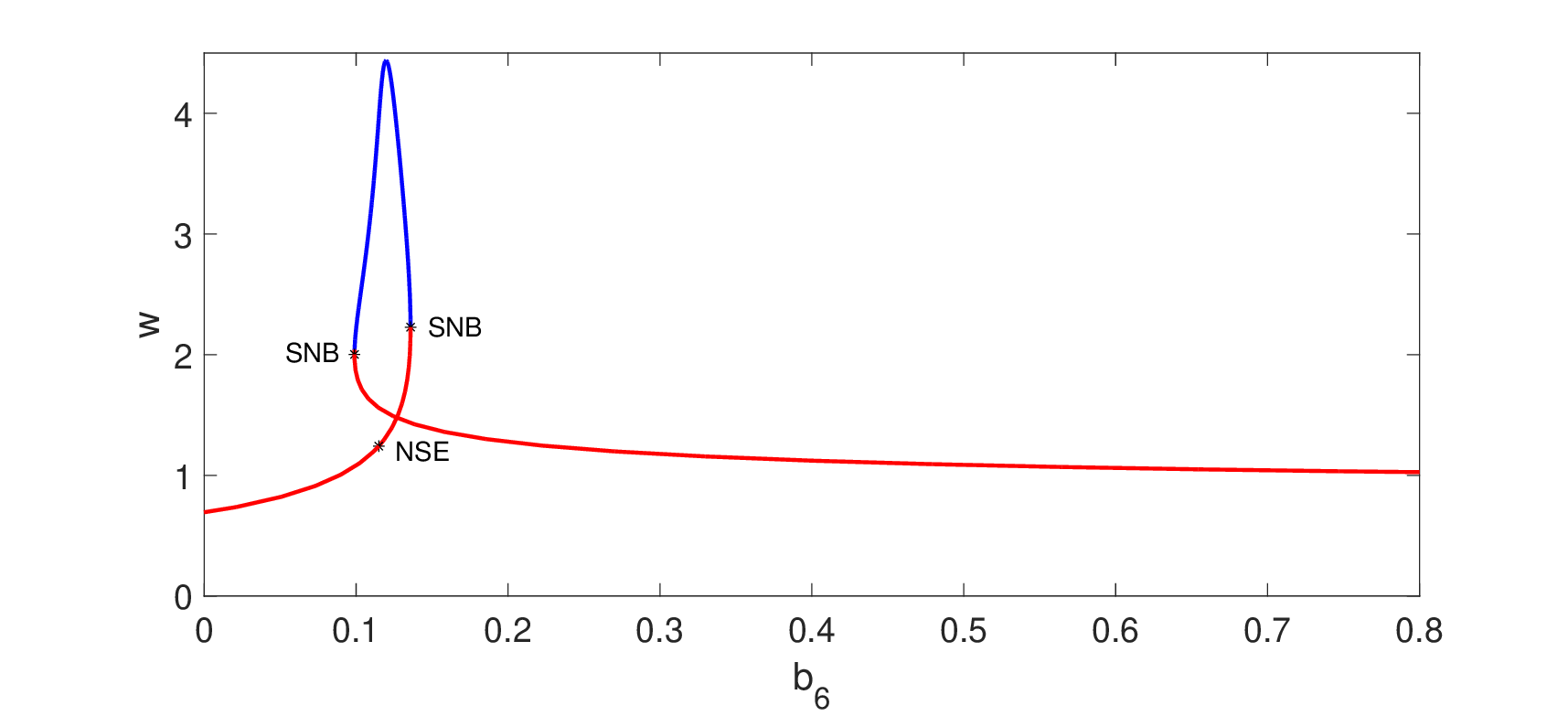} 
		\caption{\centering}
	\end{subfigure} 
	\begin{subfigure}{0.32\textwidth} 
		\centering
		\includegraphics[width= \textwidth, height=3.3cm]{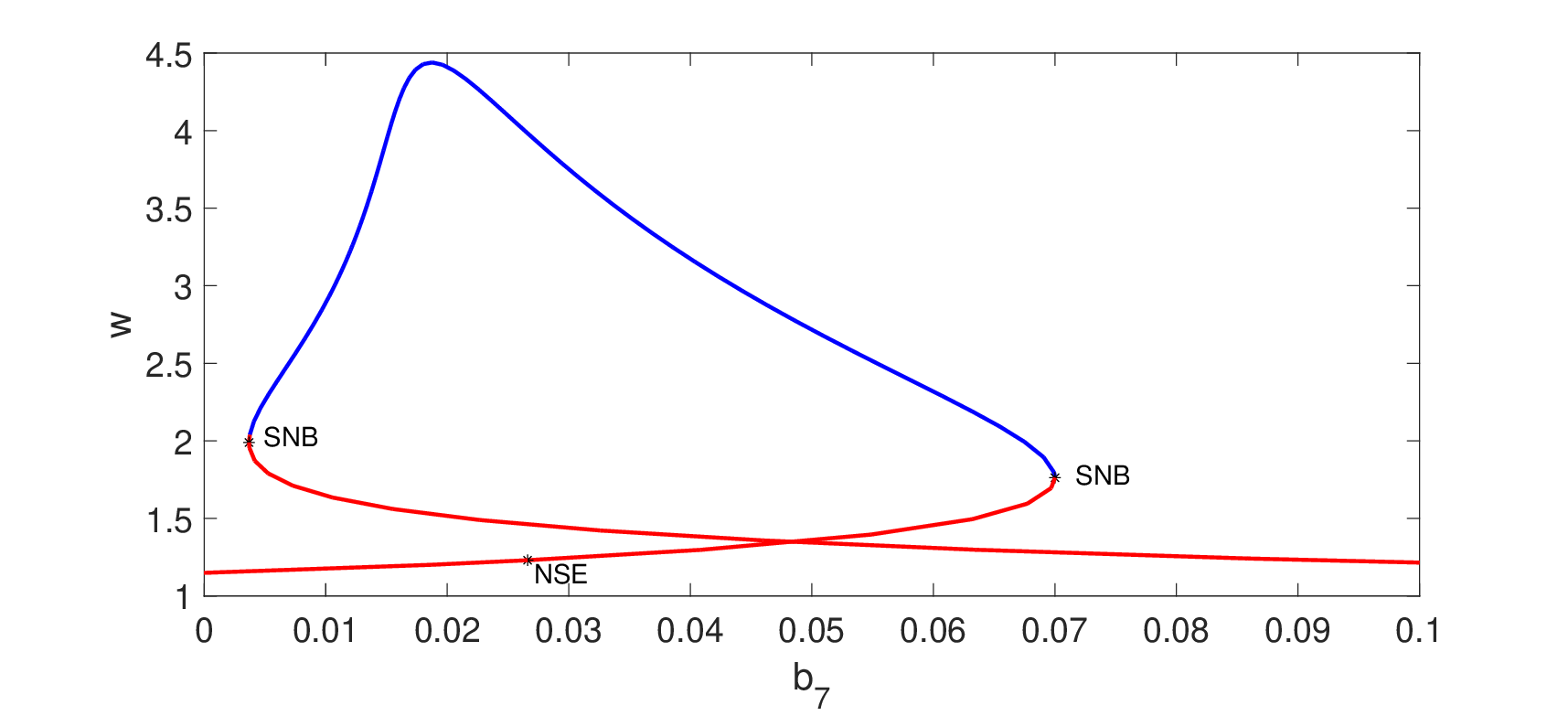} 
		\caption{\centering}
	\end{subfigure} 
	\caption{The possibilities of co-dimension one bifurcation for each parameters of the model (\ref{eq:4}) with respect to wolf population. In these figures SNB denotes the saddle node bifurcation, NSE stands for neutral saddle equilibrium.}
	\label{Fig:wolf_bifurcation}   
\end{figure}

\textbf{One parameteric bifurcation:} In Figs. \ref{Fig:elk_bifurcation} and \ref{Fig:wolf_bifurcation}, we present the one-parameter bifurcation diagrams of the data-driven model (\ref{eq:4}), using the parameter values listed in Table \ref{table:1}, which are informed by empirical data from northern Yellowstone National Park. In these figures, the blue curves denote stable equilibrium states, while the red curves represent unstable ones. Our analysis reveals the presence of saddle-node bifurcations for each parameter, with up to three critical threshold values. Among these, one bifurcation occurs between two unstable equilibria, while the other two mark transitions between stable and unstable branches. The saddle-node bifurcation points occur where the stable and unstable equilibrium branches merge, indicating critical thresholds in the ecological dynamics. Importantly, the region enclosed by the blue curves in both figures identifies the parameter ranges where coexistence of both elk and wolf populations is possible (the last column in the Table \ref{table:1}. Outside these thresholds, the system loses stability, leading to the potential extinction of one or both species. Thus, the presence of saddle-node bifurcations in the SINDy generated model delineates ecologically meaningful parameter intervals within which both elk and wolf populations can persist in the Yellowstone ecosystem. 

In the model~\eqref{eq:4}, evaluated near the numerical values of the parameters listed in Table~\ref{table:1}, we observe the occurrence of saddle-node bifurcations only, with no evidence of other bifurcation types. To explore richer dynamics, including periodic oscillations in species populations, we introduce small perturbations in the system parameters. Specifically, we modify the value of the predation related parameter $a_2=2.038$ to $a_2 = 2.138$ and examined the resulting bifurcation structure. Our analysis reveals that periodic solutions emerge only when the parameters $a_1$ (natural mortality of elk) and $a_3$ (herd protection coefficient for elk population) are varied, indicating the presence of Hopf bifurcations. Accordingly, in the following next part, we investigate two-parameter bifurcation diagrams to study the interaction between saddle-node and Hopf bifurcations, thereby capturing the transition to periodic behavior and complex ecological dynamics within the elk-wolf system.

\begin{figure}[H]
	\begin{subfigure}{0.5\textwidth} 
		\centering
		\includegraphics[width= \textwidth, height=5.5cm]{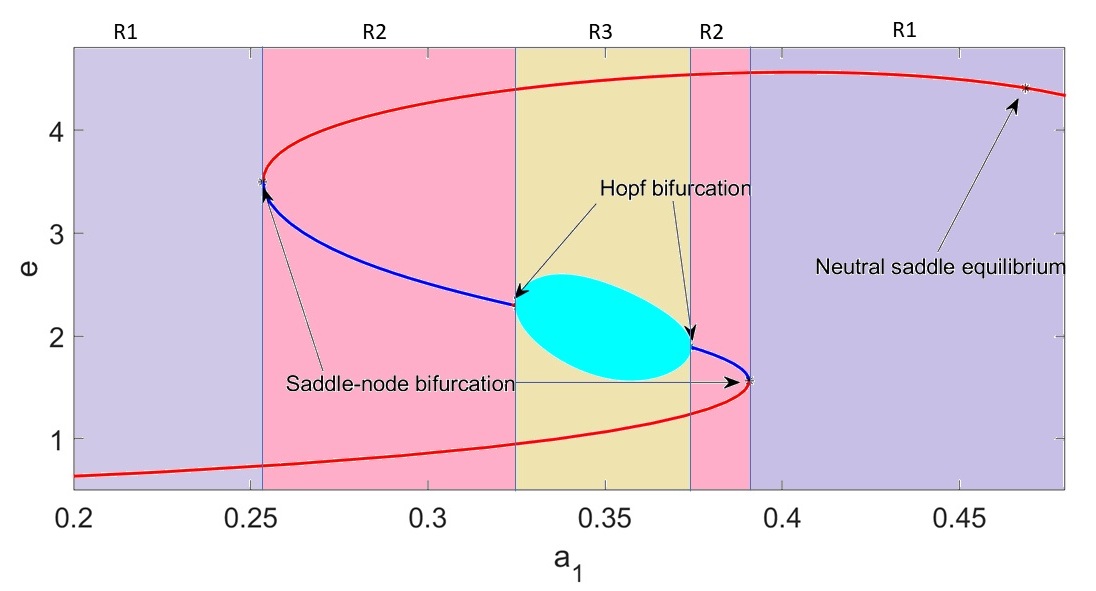} 
		\caption{\centering}
	\end{subfigure} 
	\begin{subfigure}{0.5\textwidth} 
		\centering
		\includegraphics[width= \textwidth, height=5.5cm]{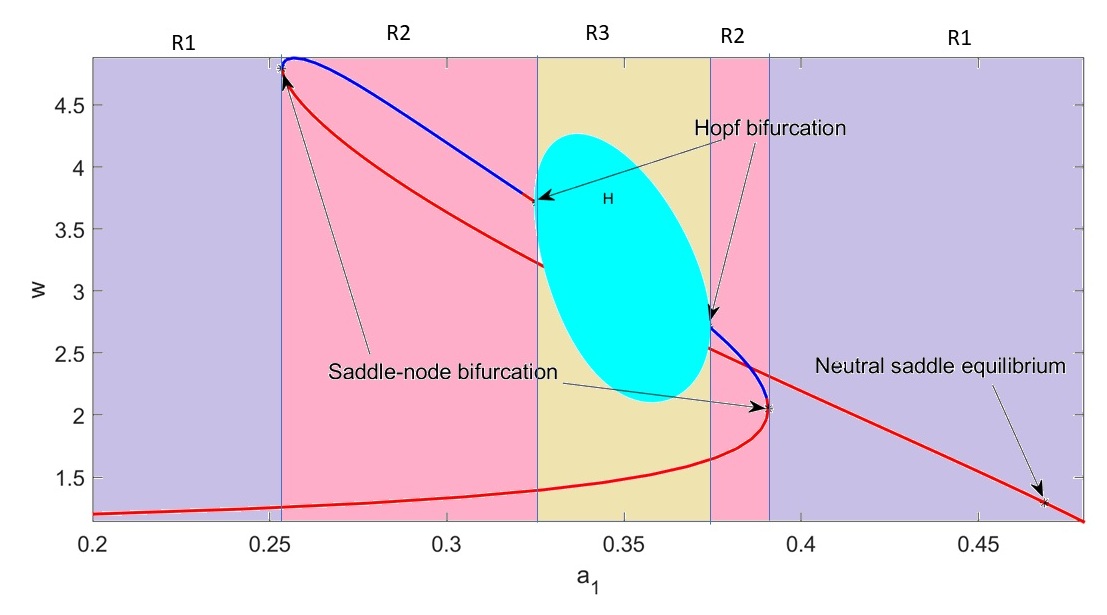} 
		\caption{\centering}
	\end{subfigure} 
	\caption{Bifurcation structure between the natural mortality of elk and the populations of elk and wolves.}
	\label{Fig:bifurcations_a1}   
\end{figure}

\begin{figure}[H]
	\begin{subfigure}{0.5\textwidth} 
		\centering
		\includegraphics[width= \textwidth, height=4.5cm]{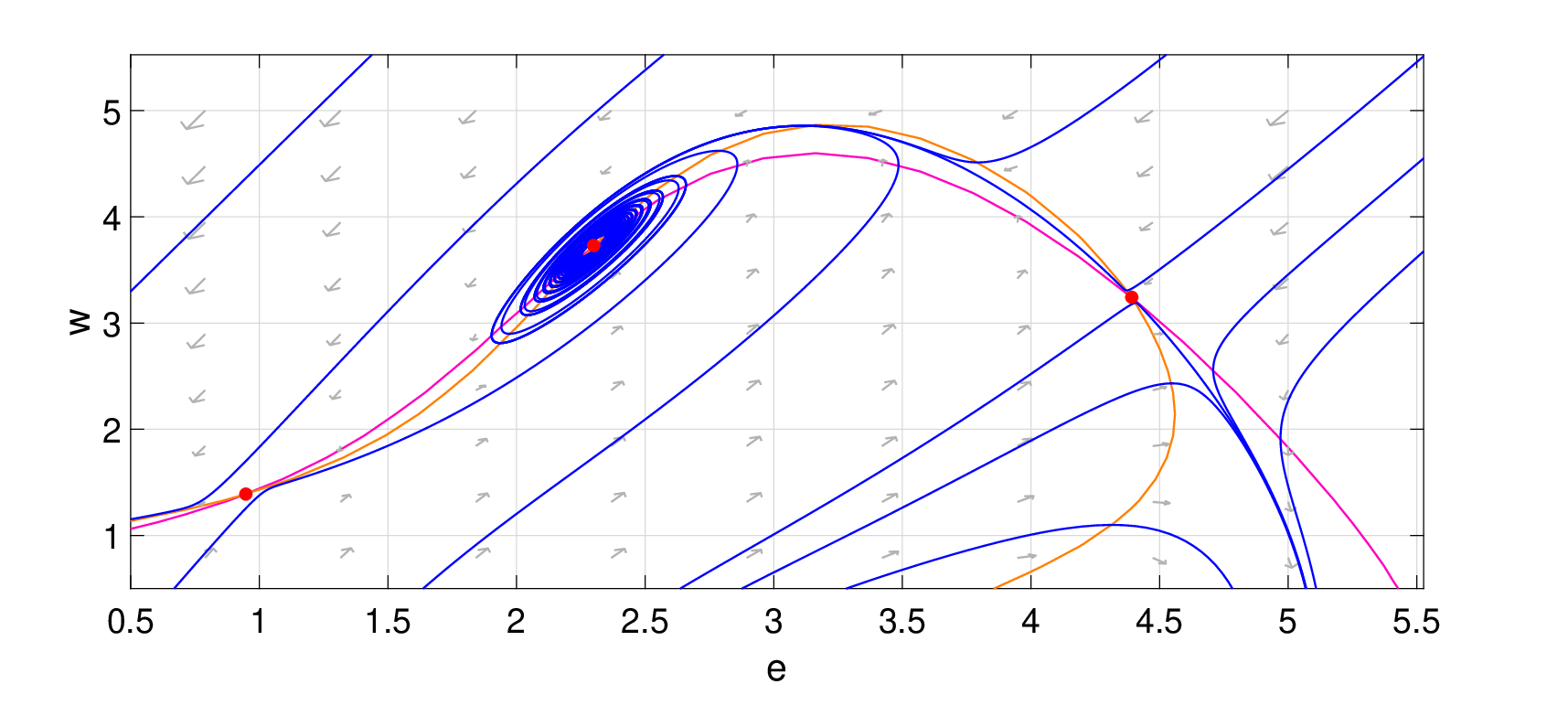} 
		\caption{\centering}
	\end{subfigure} 
	\begin{subfigure}{0.5\textwidth} 
		\centering
		\includegraphics[width= \textwidth, height=4.5cm]{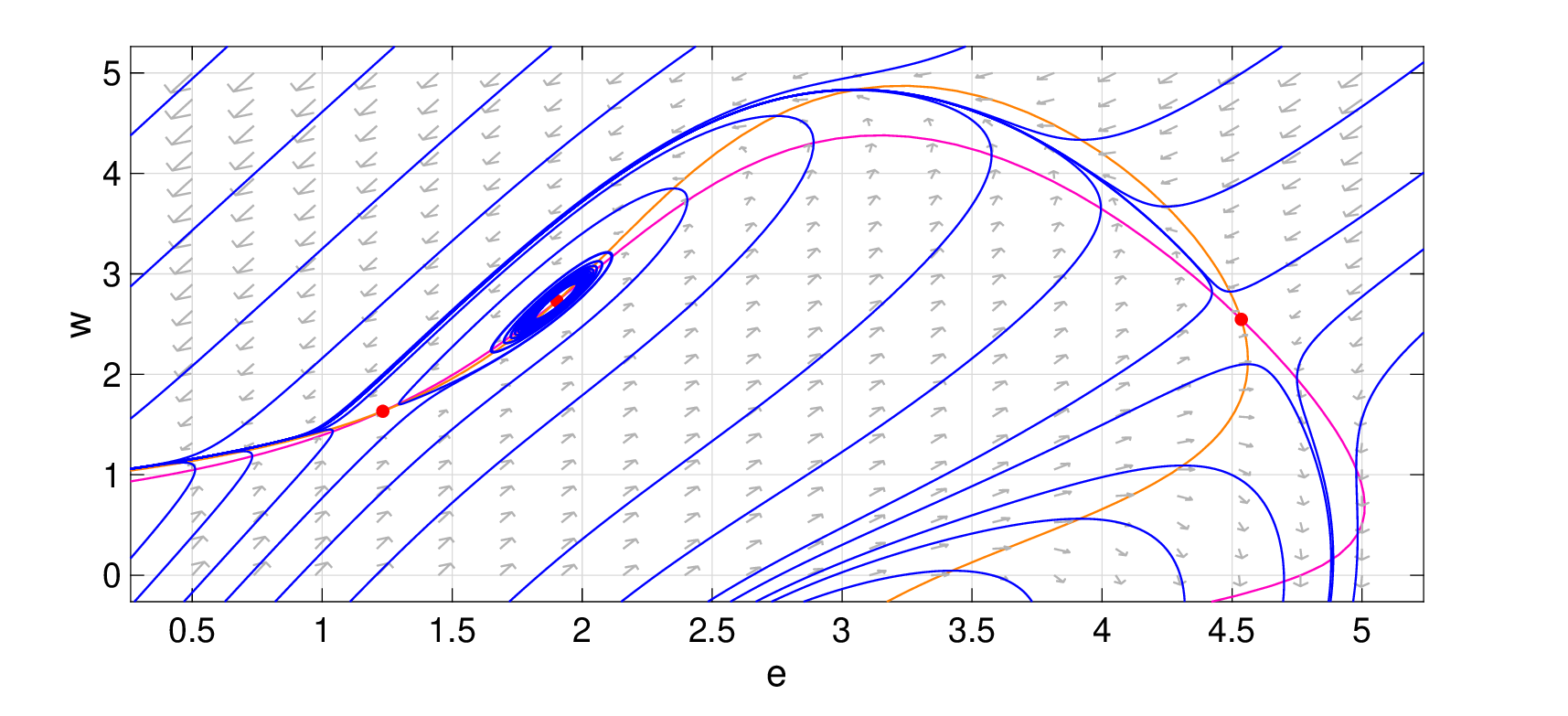} 
		\caption{\centering}
	\end{subfigure} 
	\begin{subfigure}{0.5\textwidth} 
		\centering
		\includegraphics[width= \textwidth, height=4.5cm]{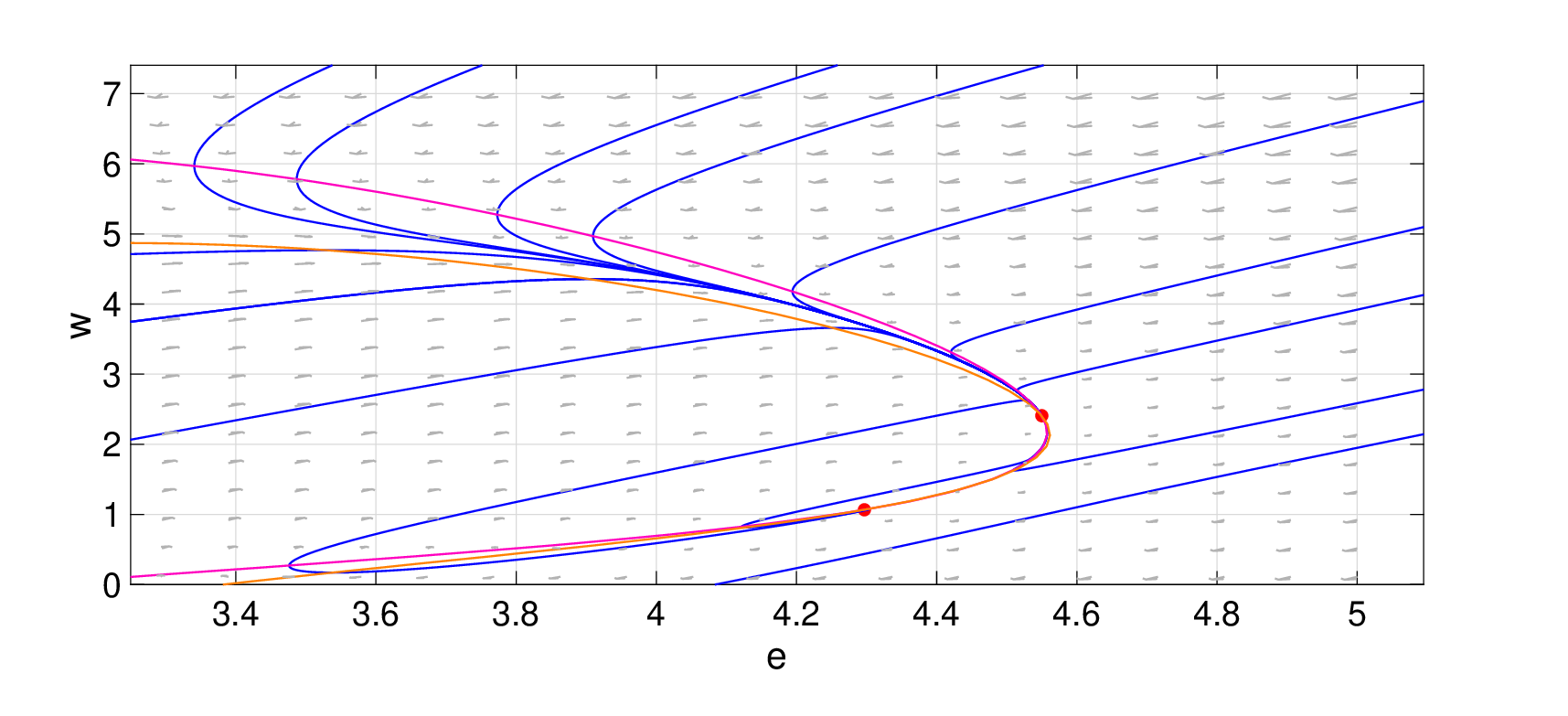} 
		\caption{\centering}
	\end{subfigure}
	\begin{subfigure}{0.5\textwidth} 
		\centering
		\includegraphics[width= \textwidth, height=4.5cm]{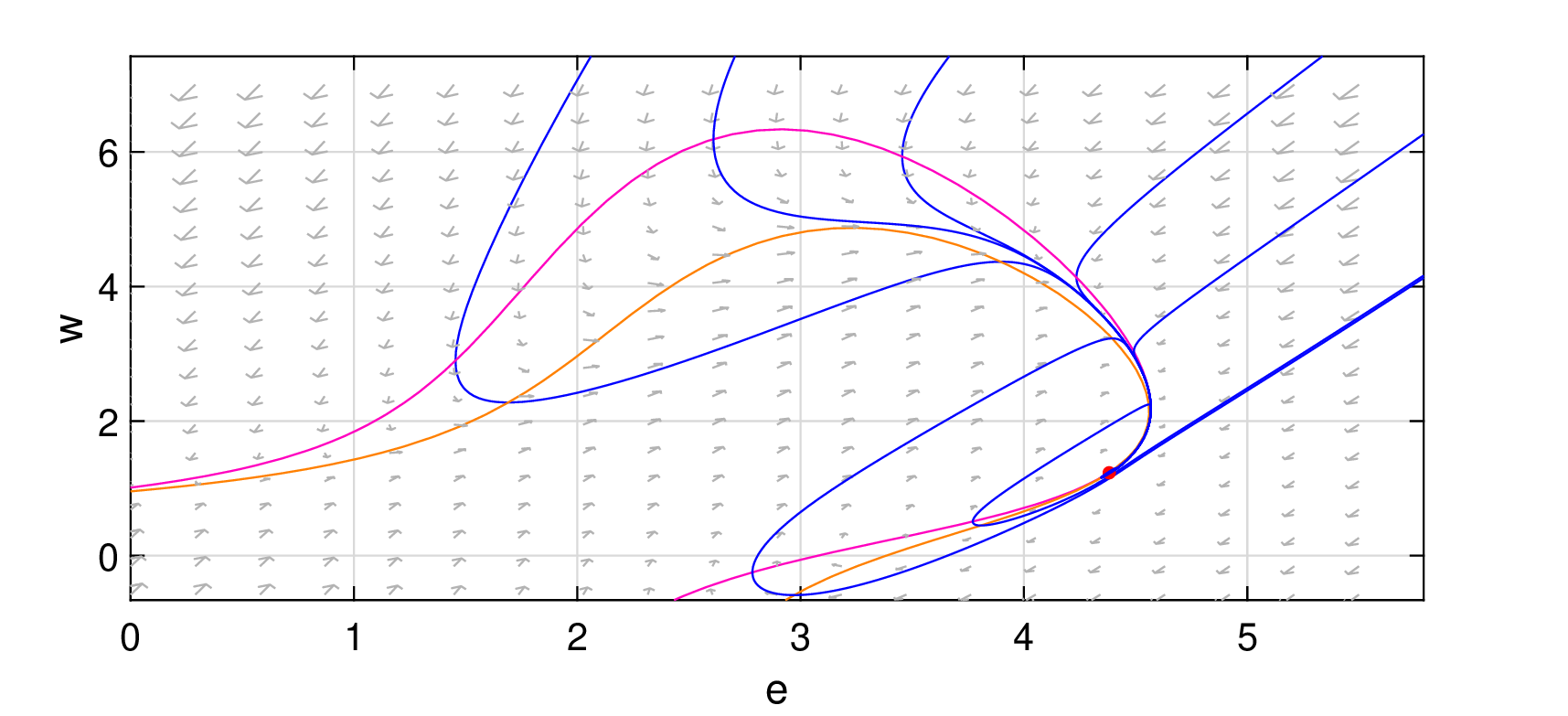} 
		\caption{\centering}
	\end{subfigure}
	\caption{Phase portraits representing at the first $a_1^{[h1]}$ and second $a_1^{[h2]}$ threshold values of the Hopf bifurcation for $a_1$ in (a) and (b). (c) At the Bogdanov-Takens bifurcation and (d) at Cusp bifurcation.}
	\label{Fig:pp_a1}   
\end{figure}

\textbf{Two parameteric bifurcation:} 
With $a_2 = 2.138$, the system (\ref{eq:4}) exhibits both periodic solutions and co-dimension two bifurcations, such as Bogdanov–Takens and Cusp bifurcations. In Fig.~\ref{Fig:bifurcations_a1}, we illustrate the possible co-dimension one bifurcations between the death rate of elk $(a_1)$ and the populations of elk and wolves. Initially, when $a_1 = 0.2534004 = a_1^{[\text{sn1}]}$, it is evident from the figure that two equilibrium branches—one stable (depicted in blue) and one unstable (depicted in red) merge. For values $a_1 < a_1^{[\text{sn1}]}$, only a single unstable equilibrium exists, whereas for $a_1 > a_1^{[\text{sn1}]}$, three equilibrium states emerge, of which only one is stable. This scenario indicates the occurrence of a saddle-node bifurcation. As $a_1$ increases further, another saddle-node bifurcation is observed at the threshold $a_1 = 0.39067856 = a_1^{[\text{sn2}]}$. Beyond this value ($a_1 > a_1^{[\text{sn2}]}$), the system (\ref{eq:4}) loses stability, making the coexistence of both elk and wolf populations unsustainable. Moreover, within the interval bounded by these two saddle-node bifurcation points, the system also exhibits periodic solutions through Hopf bifurcations \cite{kuznetsov1998elements}. These occur at two threshold values: $a_1 = 0.32514388 = a_1^{[\text{h1}]}$ and $a_1 = 0.3739507 = a_1^{[\text{h2}]}$. Therefore, for $a_1\in [a_1^{[\text{h1}]}, a_1^{[\text{h2}]}]$, Hopf bifurcations give rise to periodic solutions, forming a bubble structure in the bifurcation diagram (\ref{Fig:bifurcations_a1}). Inside the bubble both elk and wolf population coexist. Hence, a decrement in the death rate of elk form a stable limit cycle leads to the long term co-existence of both populations. At both threshold points, the system approaches a stable limit cycle, with the first Lyapunov numbers being $-0.143815482$ and $-0.255373648$ (Figs.~\ref{Fig:pp_a1}(a) and~\ref{Fig:pp_a1}(b)). 

Finally, when the natural mortality rate of elk changes, the system may pass through a neutral saddle equilibrium. In the ecological context of Yellowstone National Park, such a state implies a fragile balance: small changes in mortality may not significantly affect population dynamics immediately, but larger disturbances could destabilize the system and disrupt the long-term co-existence of elk and wolves. The colour description for Fig. \ref{Fig:bifurcations_a1} is given in Table \ref{table:2}.

\begin{table}[H]
	\centering
	\caption{Description of colored Regions in bifurcation diagram (\ref{Fig:bifurcations_a1})}
	\label{table:2}
	\renewcommand{\arraystretch}{1.4}
	\begin{tabular}{|c|c|p{13cm}|}
		\hline
		\rowcolor{gray!20}
		\textbf{Color} & \textbf{Region} & \textbf{Ecological/Bifurcation Interpretation} \\
		\hline
		\cellcolor{blue!15} & $R1$ & only one saddle equilibrium exists (not both population coexists) \\
		\hline
		\cellcolor{purple!30} & $R2$ & Presence of a saddle-node bifurcation, where two equilibria (one stable, one unstable) coalesce and annihilate each other along with third unstable equilibrium point. Indicates abrupt shifts in elk-wolf dynamics. \\
		\hline
		\cellcolor{yellow!30} & $R3$ & Region enclosing the Hopf bifurcation, indicating oscillatory behavior (limit cycles), likely corresponding to elk-wolf cycles in Yellowstone. \\
		\hline
	\end{tabular}
\end{table}

\newpage

\begin{figure}[H]
	\centering
	\includegraphics[width=\textwidth, height=8cm]{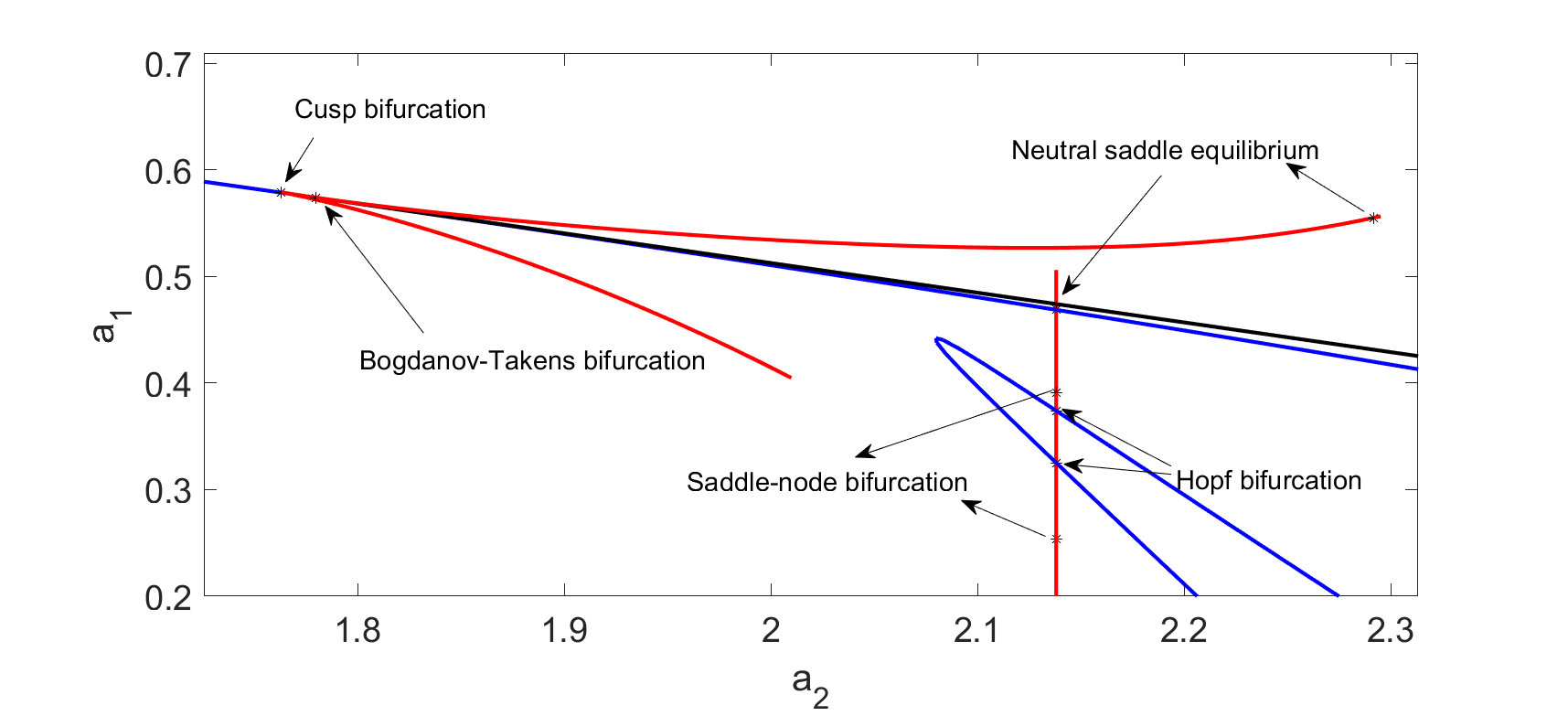}
	\caption{Co-dimension one and two bifurcations for the system (\ref{eq:4}) in the parametric space ($a_2$, $a_1$).} 
	\label{Fig:12}
\end{figure}

Furthermore, on the parametric surface defined by $a_1$–$a_2$, the system (\ref{eq:4}) exhibits co-dimension two bifurcations. Specifically, a Bogdanov–Takens bifurcation occurs at the point $(a_1, a_2) = (0.57383395,\ 1.7798861)$ (Fig. \ref{Fig:pp_a1}(c)) where the Hopf and saddle-node bifurcation curves intersect. Additionally, a cusp bifurcation appears at $(a_1, a_2) = (0.57888323,\ 1.7631351)$ (Fig. \ref{Fig:pp_a1}(d)) where two branches of saddle-node bifurcations merge. These bifurcations are illustrated in Fig.~\ref{Fig:12}. Ecologically, within the context of Yellowstone National Park, these bifurcations mark critical transitions in the dynamics of elk–wolf interactions. Here, $a_1$ denotes the natural mortality rate of elk, and $a_2$ represents the baseline predation coefficient that is, the fraction of elk lost per wolf, which may include elk killed by wolves but not necessarily consumed. The Bogdanov–Takens bifurcation indicates a delicate balance where small changes in elk mortality or wolf-induced elk loss may lead to dramatic shifts in population behavior from stable equilibria to oscillatory cycles or even extinction scenarios. The cusp bifurcation signifies a tipping point where the system transitions between having multiple equilibria (such as high and low elk densities) to a single equilibrium, indicating potential regime shifts in the ecosystem. Understanding such bifurcations is crucial for managing Yellowstone’s prey-predator balance, as even modest variations in predation intensity or elk vulnerability can drive the system toward drastically different ecological outcomes.

\begin{figure}[H]
	\begin{subfigure}{0.5\textwidth} 
		\centering
		\includegraphics[width= \textwidth, height=5cm]{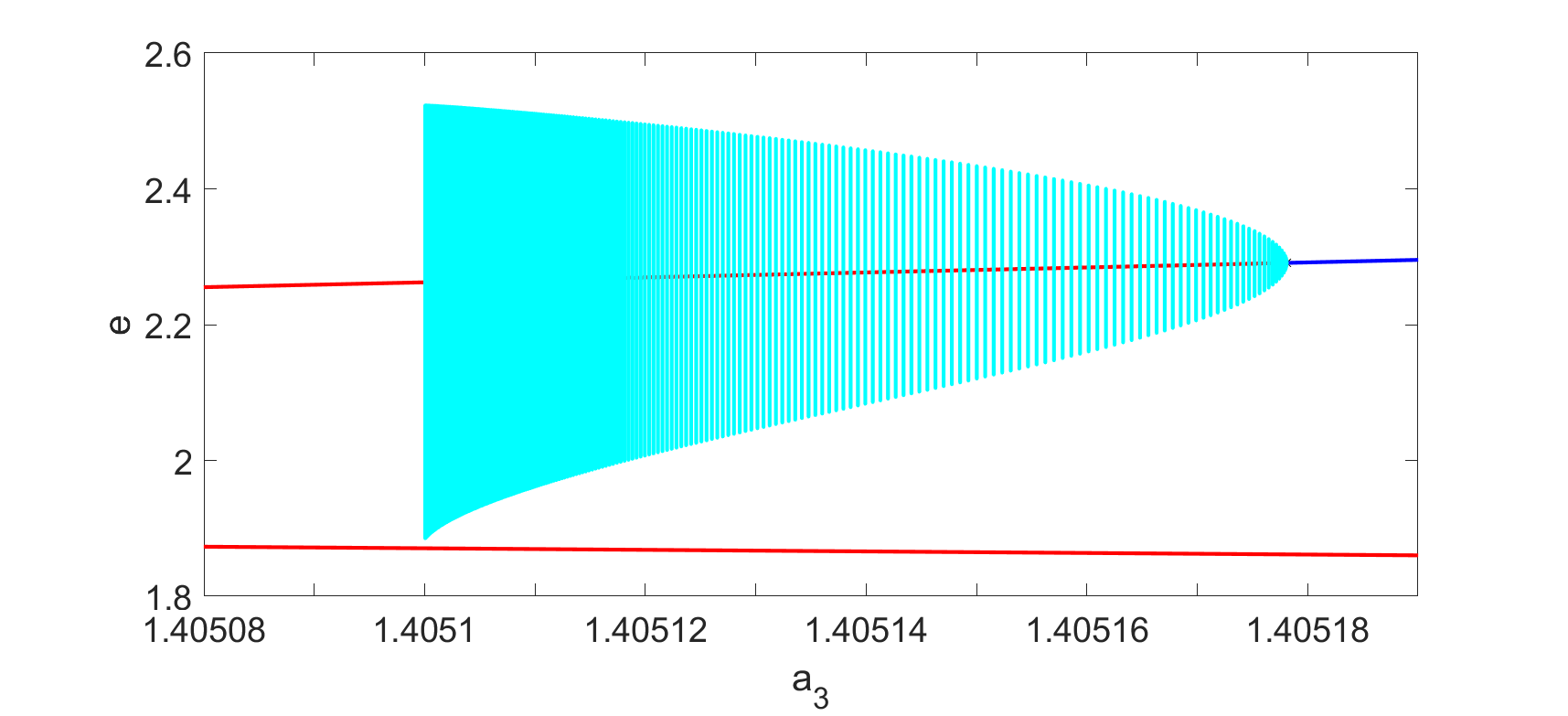} 
		\caption{\centering}
	\end{subfigure} 
	\begin{subfigure}{0.5\textwidth} 
		\centering
		\includegraphics[width= \textwidth, height=5cm]{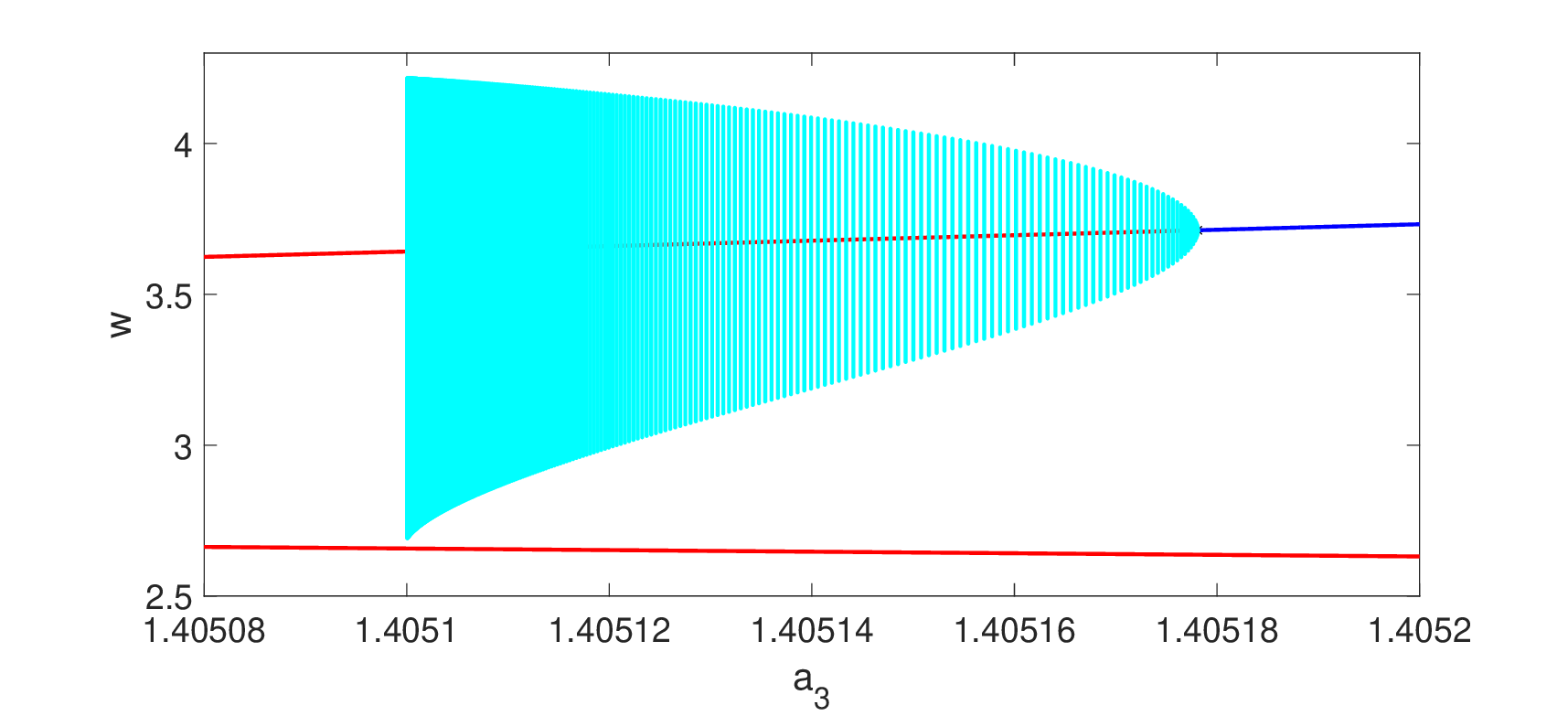} 
		\caption{\centering}
	\end{subfigure} 
	\caption{Periodic solution through Hopf bifurcation.}
	\label{Fig:bifurcations}   
\end{figure}

\begin{figure}[H]
	\centering
\begin{subfigure}{0.5\textwidth} 
		\includegraphics[width= \textwidth, height=5cm]{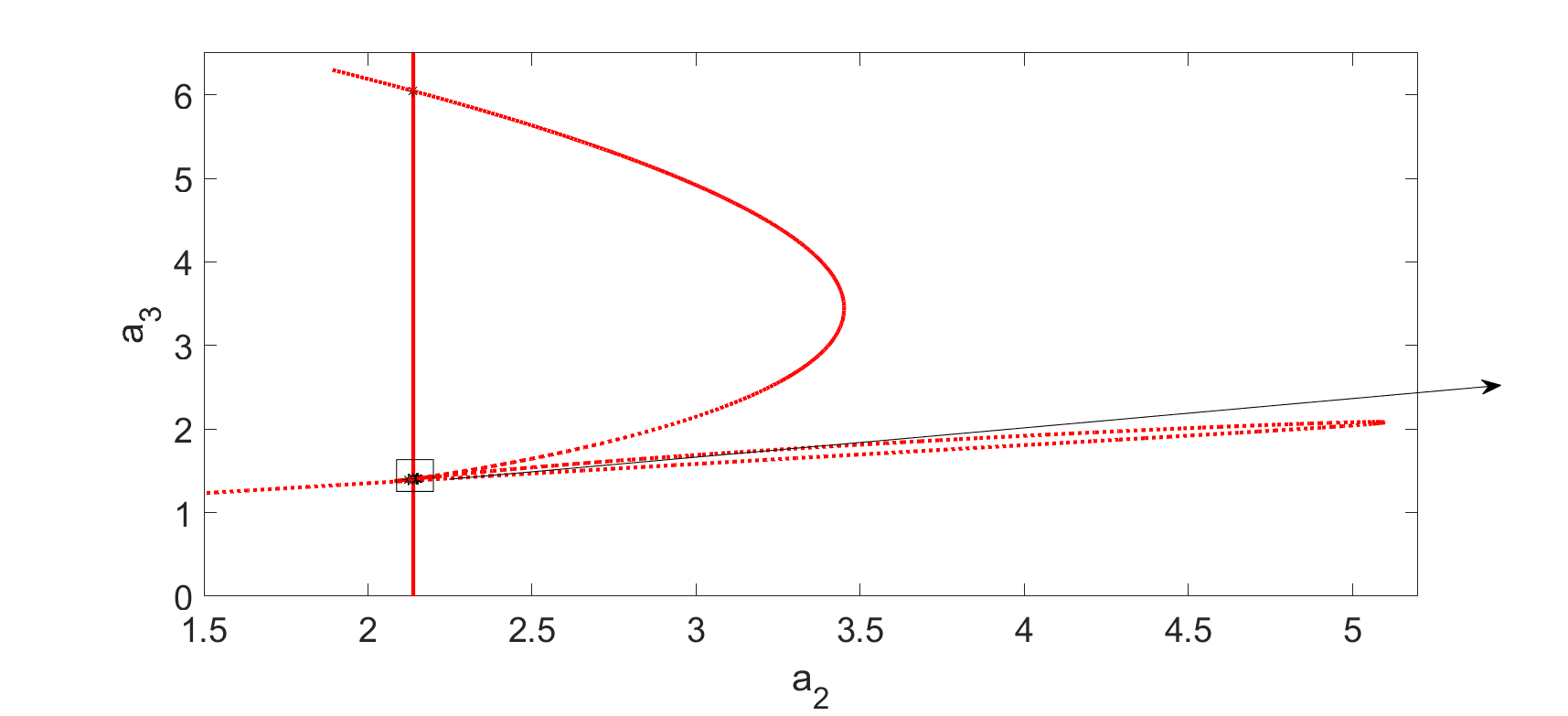} 
\end{subfigure} 
\hspace{-0.5cm}
\begin{subfigure}{0.5\textwidth} 
		\includegraphics[width= \textwidth, height=5cm]{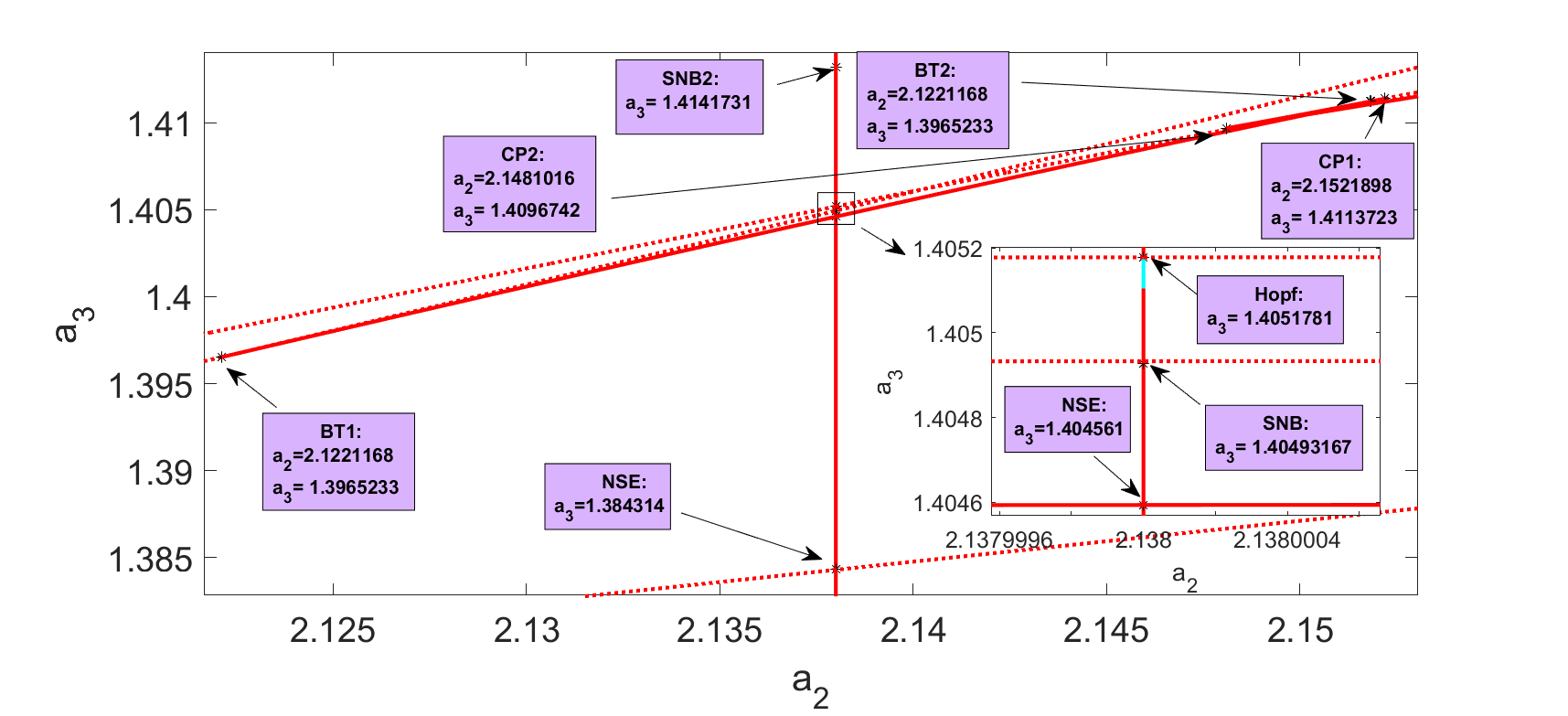} 
	\end{subfigure} \vspace{-0.1cm}
	\caption{Co-dimensions one and two bifurcations for the system (\ref{eq:4}) in the parametric space ($a_2$, $a_3$).} 
	\label{Fig:11}
\end{figure}

We have also examined the possibility of bifurcations with respect to other parameters and found that only the parameter representing the herd protection behavior of the elk population ($a_3$) is responsible for bifurcation. The bifurcation structure is illustrated in Fig. \ref{Fig:11}, where the system exhibits various co-dimension one and co-dimension two bifurcations. A periodic solution emerging from a Hopf bifurcation, associated with the parameter $a_3$, is shown in Fig. \ref{Fig:bifurcations}. An increase in herd behavior leads to periodic solutions over a very narrow range of $a_3$ values. In the context of Yellowstone National Park, herd protection behavior among elk plays a crucial role in their survival strategies, especially against predators such as wolves.

\section{Discussion and conclusion}
Gaining insight into the processes that shape ecosystem structure and function is a central objective in ecology. The predator-to-prey biomass ratio serves as a crucial indicator of trophic organization and community dynamics, offering valuable connections to various ecosystem functions and services. The elk and wolf populations play a vital role in sustaining the ecological balance and biodiversity of Yellowstone National Park. Prior to the reintroduction of wolves, the overabundance of elk had been a longstanding concern for ecologists due to overgrazing. Following the reintroduction of wolves to Yellowstone National Park in 1995 and 1996, significant ecological changes began to unfold. One of the most notable effects was the recovery of vegetation, particularly along riverbanks and riparian zones. With the return of this apex predator, elk populations were brought under better control, leading to reduced overgrazing pressure. As a result, willow, aspen, and cottonwood stands began to regenerate, contributing to the restoration of critical habitats for numerous other species. This trophic cascade highlighted the wolves role in helping to rebalance the park’s ecosystem and underscored the importance of top-down regulation in maintaining ecological integrity.

Maintaining ecological balance in the northern part of Yellowstone requires the co-existence of both elk and wolves. To effectively understand and analyze the dynamics of this prey-predator system, it is essential to develop a robust mathematical model whose predictions and implications closely reflect the real-world behavior observed in the park. Modeling real-world ecological systems presents numerous challenges due to their complexity and the influence of multiple interacting factors. Traditional mathematical modeling approaches often rely on parameter estimation techniques to align models with empirical data; however, their predictive accuracy and robustness are often limited. This is largely because such models are built upon numerous simplifying assumptions that may not fully capture the intricacies of natural systems.

In this study, we explored the data-driven mathematical model of the elk–wolf system in northern Yellowstone National Park, obtained using the E-SINDy framework, for its stability and bifurcation structures. The model parameters were interpreted ecologically to ensure their biological relevance to the Yellowstone ecosystem. We find the equilibrium points of the model and check the stability of the model around these equilibrium points. We also perform the bifurcation analysis of the model, which reveals interesting dynamics of the model. To the best of our knowledge, these results are reported for the first time in the literature.

Our analysis revealed that the system admits only three positive equilibrium states under the original parameter set. Among these, two are saddle points: one with a very low elk density and the other with a high elk density. These reflect ecological scenarios where either the elk population is too sparse to support wolves, or so abundant that their herd behavior leads to instability. Between these saddles, a single stable equilibrium exists, representing a biologically meaningful state where both elk and wolf populations coexist in balance a state reflective of the current dynamics in Yellowstone.

Bifurcation analysis showed the presence of co-dimension one and co-dimension two bifurcations. Under the original parameters, the system exhibits a saddle-node bifurcation. Ecologically, this indicates threshold values of key parameters (e.g., elk herd protection behavior) that separate regimes of extinction and co-existence. Below or above these thresholds, co-existence fails, highlighting the delicate balance required for maintaining prey-predator stability in the park. To further uncover complex dynamics, we perturbed the baseline predation coefficient and studied bifurcations by varying two ecologically relevant parameters: the natural mortality of elk and their herd protection coefficient. This led to the emergence of rich bifurcation structures, including Hopf, Bogdanov–Takens, and cusp bifurcations. These bifurcations capture critical transitions such as the onset of periodic cycles (e.g., population booms and crashes), sudden disappearance of equilibria, and shifts between ecological regimes.

From an ecological perspective, such bifurcations represent tipping points where small changes in environmental or biological parameters can lead to significant shifts in species dynamics such as oscillatory cycles of prey-predator abundance, abrupt elk population declines, or predator extinction. These insights underline the importance of parameter sensitivity in ecological systems and demonstrate how data-driven models can be used to predict and manage species interactions in protected habitats like northern Yellowstone National Park.

\section*{Authorship Contribution Statement}
\textbf{Anurag Singh}: Writing – original draft, Software, Methodology, Visualization, Conceptualization. \textbf{Nitu Kumari}: Writing – review \& editing, Visualization, Supervision, Conceptualization. \textbf{Arun Kumar}: Writing – original draft, Software, Methodology, Visualization, Conceptualization.

\section*{Acknowledgments}
The research of the corresponding author (Nitu Kumari) was supported by the Council of Scientific and Industrial Research (CSIR, India), under the research grant $\#$ $25\slash0326\slash23\slash\text{EMR}-\text{II}$. We thank Dr. Matt Metz, Biologist, Yellowstone center for resources for kindly providing the data file and related literature that supported this study.

\section*{Conflict of Interest}
The authors declare no conflict of interest.

\section*{Data Availability}
The empirical data used in this study were obtained through direct communication with the Yellowstone center for resources, U.S. national park service, upon request. 

\bibliographystyle{unsrt}
\bibliography{ref}

\section*{Appendix} \label{appendix}

\begin{table}[H]
\centering
\caption{Model selection results. Model selected as the best ones (Model \# $53$) are indicated
by bold text. $k$ is the total number of terms in each model system.}
\label{table:3}
\begin{minipage}[t]{0.49\textwidth}
\scriptsize
\centering
\begin{tabular}{p{1cm}p{1cm}p{1cm}p{1cm}p{1cm}}
\toprule
\textbf{Model} & \textbf{k} & \textbf{SSE} & \textbf{AIC} & \textbf{BIC} \\
\midrule
1 & 18 & 4.64 & -716.46 & -657.09 \\
2 & 19 & 16.14 & -465.32 & -402.65 \\
3 & 19 & 7.77 & -611.51 & -548.84 \\
4 & 19 & 4.67 & -713.24 & -650.57 \\
5 & 18 & 5.74 & -673.90 & -614.53 \\
6 & 17 & 8.08 & -607.58 & -551.51 \\
7 & 18 & 4.73 & -712.63 & -653.26 \\
8 & 18 & 7.83 & -611.97 & -552.60 \\
9 & 17 & 8.06 & -608.13 & -552.06 \\
10 & 17 & 16.45 & -465.55 & -409.48 \\
11 & 16 & 20.25 & -425.99 & -373.21 \\
12 & 16 & 10.08 & -565.43 & -512.65 \\
13 & 18 & 5.46 & -684.08 & -624.71 \\
14 & 16 & 7.61 & -621.71 & -568.93 \\
15 & 15 & 7.53 & -625.64 & -576.16 \\
16 & 15 & 15.42 & -482.50 & -433.02 \\
17 & 16 & 20.89 & -419.80 & -367.02 \\
18 & 17 & 14.73 & -487.56 & -431.49 \\
19 & 15 & 18.39 & -447.24 & -397.77 \\
20 & 14 & 3.09 & -805.92 & -759.74 \\
21 & 14 & 221.26 & 48.20 & 94.38 \\
22 & 15 & 11.71 & -537.56 & -488.08 \\
23 & 13 & 1.14 &  -1007.3 & -964.52 \\
24 & 16 & 15.19 & -483.48 & -430.71 \\
25 & 13 & 4.18 & -747.18 & -704.30 \\
26 & 14 & 12.56 & -525.55 & -479.37 \\
27 & 13 & 15.79 & -481.70 & -438.83 \\
28 & 16 & 14.06 & -498.85 & -446.08 \\
29 & 13 & 1.21 & -994.22 & -951.34 \\
30 & 13 & 5.41e+09 & 3448.98 & 3491.86 \\
31 & 13 & 5422.45 & 686.00 & 728.87 \\
32 & 11 & 48.18 & -262.64 & -226.35 \\
33 & 15 & 146049.35 & 1348.67 & 1398.15 \\
34 & 12 & 5512.24 & 687.28 & 726.86 \\
35 & 12 & 3.17e+16 & 6563.70 & 6603.28 \\
36 & 12 & 1.00e+33 & 14163.07 & 14202.65 \\
37 & 13 & 1.00e+33 & 14165.07 & 14207.95 \\
38 & 19 & 3.16 & -791.28 & -728.61 \\
39 & 17 & 18.48 & -442.30 & -386.22 \\
40 & 17 & 8.07 & -608.02 & -551.95 \\
41 & 16 & 10.42 & -558.83 & -506.05 \\
42 & 16 & 16.89 & -462.20 & -409.43 \\
43 & 18 & 9.71 & -569.00 & -509.63 \\
44 & 16 & 13.87 & -501.58 & -448.80 \\
\bottomrule
\end{tabular}
\end{minipage}
\hfill
\begin{minipage}[t]{0.49\textwidth}
\scriptsize
\centering
\begin{tabular}{p{1cm}p{1cm}p{1cm}p{1cm}p{1cm}}
\toprule
\textbf{Model} & \textbf{k} & \textbf{SSE} & \textbf{AIC} & \textbf{BIC} \\
\midrule
45 & 17 & 11.04 & -545.30 & -489.23 \\
46 & 15 & 19.89 & -431.54 & -382.06 \\
47 & 15 & 1.99 & -891.45 & -841.97 \\
48 & 15 & 25.00 & -385.84 & -336.36 \\
49 & 16 & 1.24 & -984.25 & -931.48 \\
50 & 14 & 2.01 & -891.75 & -845.57 \\
51 & 15 & 9.43 & -580.86 & -531.38 \\
52 & 14 & 25.19 & -386.32 & -340.15 \\
\textbf{53} & \textbf{14} & \textbf{0.91} & \textbf{-1049.6} & \textbf{-1003.5} \\
54 & 10 & 60.00 & -220.76 & -187.78 \\
55 & 11 & 279.97 & 89.27 & 125.55 \\
56 & 13 & 18003.86 & 926.00 & 968.88 \\
57 & 11 & 281.59 & 90.42 & 126.70 \\
58 & 12 & 64.81 & -201.34 & -161.76 \\
59 & 12 & 1.28e+16 & 6383.49 & 6423.07 \\
60 & 11 & 117.97 & -83.56 & -47.28 \\
61 & 10 & 153.98 & -32.29 & 0.68 \\
62 & 9 & 3.99e+17 & 7064.07 & 7093.76 \\
63 & 9 & 274.35 & 81.22 & 110.90 \\
64 & 9 & 216.75 & 34.08 & 63.77 \\
65 & 8 & 220.63 & 35.63 & 62.02 \\
66 & 10 & 1.55e+16 & 6417.21 & 6450.19 \\
67 & 9 & 1.37e+17 & 6851.21 & 6880.90 \\
68 & 20 & 10.52 & -548.94 & -482.97 \\
69 & 19 & 3.61 & -764.41 & -701.74 \\
70 & 18 & 4.29 & -731.98 & -672.61 \\
71 & 18 & 31.35 & -334.60 & -275.23 \\
72 & 18 & 41.42 & -278.88 & -219.51 \\
73 & 16 & 45.80 & -262.79 & -210.02 \\
74 & 17 & 25.77 & -375.79 & -319.72 \\
75 & 17 & 51.55 & -237.14 & -181.07 \\
76 & 14 & 38.54 & -301.29 & -255.11 \\
77 & 14 & 34.88 & -321.25 & -275.08 \\
78 & 16 & 89.15 & -129.58 & -76.80 \\
79 & 14 & 47.30 & -260.34 & -214.16 \\
80 & 15 & 54.66 & -229.40 & -179.93 \\
81 & 14 & 36.45 & -312.43 & -266.25 \\
82 & 14 & 35.00 & -320.54 & -274.37 \\
83 & 13 & 43.91 & -277.23 & -234.35 \\
84 & 14 & 125407.89 & 1316.20 & 1362.37 \\
85 & 14 & 2.37e+10 & 3746.84 & 3793.01 \\
86 & 9 & 3149.97 & 569.36 & 599.05 \\
87 & 12 & 40.88 & -293.52 & -253.94 \\
88 & 14 & 50094.30 & 1132.66 & 1178.84 \\
\bottomrule
\end{tabular}
\end{minipage}
\end{table}

\begin{table}[H]
\centering
\begin{minipage}[t]{0.48\textwidth}
\scriptsize
\centering
\begin{tabular}{p{1cm}p{1cm}p{1cm}p{1cm}p{1cm}}
\toprule
\textbf{Model} & \textbf{k} & \textbf{SSE} & \textbf{AIC} & \textbf{BIC} \\
\midrule
89 & 11 & 155.24 & -28.65 & 7.62 \\
90 & 10 & 3234.00 & 576.63 & 609.61 \\
91 & 13 & 141015.21 & 1337.66 & 1380.53 \\
92 & 10 & 2720.37 & 542.04 & 575.02 \\
93 & 11 & 215.86 & 37.26 & 73.55 \\
94 & 8 & 240.31 & 52.72 & 79.11 \\
95 & 8 & 239.28 & 51.86 & 78.25 \\
96 & 9 & 218.44 & 35.64 & 65.32 \\
97 & 9 & 229.81 & 45.79 & 75.48 \\
98 & 17 & 4.70 & -715.86 & -659.79 \\
99 & 17 & 4.70e+13 & 5270.61 & 5326.68 \\
100 & 18 & 4.22 & -735.23 & -675.87 \\
101 & 17 & 4.27 & -735.11 & -679.04 \\
102 & 16 & 3.33e+10 & 3818.20 & 3870.97 \\
 &        &         &         &          \\
\bottomrule
\end{tabular}
\end{minipage}
\hfill
\begin{minipage}[t]{0.48\textwidth}
\scriptsize
\centering
\begin{tabular}{p{1cm}p{1cm}p{1cm}p{1cm}p{1cm}}
\toprule
\textbf{Model} & \textbf{k} & \textbf{SSE} & \textbf{AIC} & \textbf{BIC} \\
\midrule
103 & 17 & 12.34 & -522.93 & -466.86 \\
104 & 15 & 5.91 & -674.20 & -624.73 \\
105 & 16 & 2.89 & -814.82 & -762.05 \\
106 & 12 & 4.64e+14 & 5718.55 & 5758.13 \\
107 & 13 & 4.81e+10 & 3885.67 & 3928.54 \\
108 & 13 & 4.84e+10 & 3887.28 & 3930.16 \\
109 & 10 & 5.28e+15 & 6201.24 & 6234.22 \\
110 & 12 & 1.81e+11 & 4149.15 & 4188.73 \\
111 & 11 & 3.13e+11 & 4256.58 & 4292.86 \\
112 & 10 & 71.19 & -186.58 & -153.60 \\
113 & 9 & 144.33 & -47.23 & -17.55 \\
114 & 10 & 78.08 & -168.09 & -135.11 \\
115 & 7 & 241.33 & 51.57 & 74.66 \\
116 & 8 & 238.10 & 50.88 & 77.26 \\
117 & 6 & 1015.48 & 336.96 & 356.75 \\
\bottomrule
\end{tabular}
\end{minipage}
\end{table}

\end{document}